\documentclass[a4paper]{amsart}
\pdfoutput=1
\usepackage[headings]{fullpage}
\usepackage{xcolor}
\usepackage{booktabs}
\usepackage{hyperref}
\hypersetup{colorlinks=true,
	linkcolor=blue!50!black,
	citecolor=blue!50!black,
	urlcolor=blue!50!black,
	filecolor=blue!50!black,
	pdfborder={0 0 0},
}

\usepackage[a4paper,margin=0.9in]{geometry}
\usepackage{amssymb,amsfonts,amsmath,amsthm}
\usepackage{latexsym,mathtools,stmaryrd}
\usepackage{tikz}
\usetikzlibrary{decorations.pathmorphing}
\tikzset{snake it/.style={decorate, decoration=snake}}
\usetikzlibrary{patterns, shadings}
\usetikzlibrary{hobby}
\usetikzlibrary{decorations.markings}

\newtheorem{theorem}{Theorem}[section]
\theoremstyle{definition}
\newtheorem{proposition}[theorem]{Proposition}
\newtheorem{lemma}[theorem]{Lemma}
\newtheorem{definition}[theorem]{Definition}
\newtheorem{remark}[theorem]{Remark}
\newtheorem{corollary}[theorem]{Corollary}
\newtheorem{conjecture}[theorem]{Conjecture}

\newtheorem{example}[theorem]{Example}

\usepackage{amssymb,amsmath,mathtools,bbm,tikz,tikz-cd,color,relsize,multirow}
\usepackage{diagbox}

\newcommand*\smat[4]{\begin{smallmatrix}#1&#2\\#3&#4\end{smallmatrix}}

\newcommand{\ol}[1]{\overline{#1}}

\newcommand*\ang[1]{\langle#1\rangle}

\newcommand{\sm}{\smallsetminus}

\newcommand{\CC}{\mathbb{C}}

\newcommand{\QQ}{\mathbb{Q}}
\newcommand{\RR}{\mathbb{R}}

\newcommand{\ZZ}{\mathbb{Z}}
\newcommand{\calC}{\mathcal{C}}
\newcommand{\Li}{\mathrm{Li}}

\newcommand{\rk}{\operatorname{rk}}
\newcommand{\tr}{\operatorname{tr}}

\newcommand{\eps}{\varepsilon}
\renewcommand{\phi}{\varphi}

\newcommand{\fa}{\mathfrak{a}}
\newcommand{\fb}{\mathfrak{b}}
\newcommand{\fc}{\mathfrak{c}}
\newcommand{\fd}{\mathfrak{d}}

\newcommand{\SL}{\operatorname{SL}}
\newcommand{\GL}{\operatorname{GL}}

\newcommand{\Hom}{\operatorname{Hom}}

\newcommand{\ee}{\mathbf{e}}

\newcommand{\tq}{\widetilde{q}}
\def\Res#1{\underset{#1}{\mathrm{Res}}}

\newcommand{\BZ}{\mathbb{Z}}

\newcommand{\BQ}{\mathbb{Q}}
\newcommand{\BR}{\mathbb{R}}
\newcommand{\BC}{\mathbb{C}}
\newcommand{\BK}{\mathbb{K}}
\newcommand{\ga}{\gamma}
\newcommand{\calZ}{\mathcal{Z}}
\newcommand{\ve}{\varepsilon}
\newcommand{\calO}{\mathcal{O}}
\newcommand{\calT}{\mathcal{T}}
\newcommand{\ZTV}{\mathrm{Z}}
\newcommand{\ZTVred}{\mathrm{Z}^{\mathrm{red}}}
\newcommand{\dimq}{\mathrm{dim}_q}
\newcommand{\End}{\mathrm{End}}
\newcommand{\Aut}{\mathrm{Aut}}
\newcommand{\HF}{\mathrm{HF}}
\renewcommand{\rk}{\mathrm{rk}}
\newcommand{\Sk}{\mathrm{Sk}}

\begin{document}

\title[Real quadratic fields and finite quantum dilogarithms II]{
Real quadratic fields and finite quantum dilogarithms II:\\ 
traces of quantum modular cocycles}

\author{Stavros Garoufalidis}
\address[S. Garoufalidis]{
  International Center for Mathematics, Department of Mathematics \\
  Southern University of Science and Technology \\
  Shenzhen, China \newline
  {\tt \url{http://people.mpim-bonn.mpg.de/stavros}}}
\email{stavros@mpim-bonn.mpg.de}

\author{Danylo Radchenko}
\address[D. Radchenko]{Institut des Hautes \'Etudes Scientifiques, 
CNRS, Laboratoire Alexander Grothendieck,
	35 route de Chartres, Bures-sur-Yvette 91440, France\newline
  {\tt \url{https://www.danrad.net}}}
\email{danradchenko@gmail.com}

\author{Campbell Wheeler}
\address[C. Wheeler]{Coburg 3058, Australia\newline
  {\tt \url{https://www.ihes.fr/~wheeler/}}}
\email{wheeler@ihes.fr}

\date{9 November 2026}
\subjclass[2010]{}
\keywords{Real quadratic fields, quantum dilogarithm, finite 
quantum dilogarithm, quantum modular cocycle, 
$\mathrm{SL}_2(\mathbb{Z})$, TQFT, fusion categories, 
3-manifolds, linear $q$-difference equations, 
state-integrals, algebraic units.}

\begin{abstract}
We study traces of quantum modular cocycles at fixed points of 
hyperbolic elements $\gamma\in\SL_2(\BZ)$, formulate
an algebraicity conjecture for these traces, and prove it for 
the cocycles associated to 1-dimensional state-integrals. In 
addition, we analyse a finite TQFT associated with $\gamma$ 
and conjecture that its invariants for 3-manifolds agree 
with the traces of their $\mathfrak{sl}_2$ quantum modular
cocycles. The finite TQFTs predict a yet-to-be-defined motive
with cohomology equipped with a $\BZ/2\BZ$-degree 
grading and a $\BZ$-weight grading coming from a version of the 
Riemann 
hypothesis, which suggests a similar  grading should exist for 
an extension of the Kauffman bracket skein module,
for the framed $\SL_2(\BC)$ Floer cohomology of 
Abouzaid--Manolescu, and for the quantum modular cocycle.
This leads to  surprising additive identities among algebraic 
units, which we confirm   for the complement of the simplest 
hyperbolic knot $4_1$, as well as for $S^1 \times \Sigma_g$
for a closed genus $g$ surface.
\end{abstract}

\maketitle

{\footnotesize
\tableofcontents
}

\section{Introduction}

\subsection{Quantum modular cocycles}
\label{sub.modular}

Matrix-valued holomorphic quantum modular forms were introduced 
in~\cite{GZ23,GZ24} as a way of unifying the intertwined properties 
of exact and perturbative complex Chern--Simons theory.
The concept that emerged, a quantum modular cocycle, was the 
result of a gradual evolution of ideas and computations over many 
years building on an earlier notion of a quantum modular form 
introduced by Zagier~\cite{Zag10}. 
By definition, a quantum modular cocycle is a collection 
$\Omega_\gamma$ of matrix-valued holomorphic functions on 
the cut-plane 
$\BC_\gamma=\BC\sm\{\tau \, | \, c\tau+d \leq 0\}$ that is a 
cocycle of $\SL_2(\BZ)$, i.e., it satisfies
\begin{equation}
\label{Omcoc}
\Omega_{\gamma \gamma'}(\tau)
= \Omega_{\gamma}(\gamma'(\tau))  \Omega_{\gamma'}(\tau) 
\end{equation}
for all $\gamma, \gamma' \in \SL_2(\BZ)$.

One of the main conjectures of~\cite{GZ23} states that complex 
Chern--Simons  theory associates a quantum
modular cocycle $\Omega_M$ to every compact, oriented 3-manifold 
$M$, possibly with torus boundary components. This conjecture is 
currently known only for a handful of 
hyperbolic knot complements (such as the three simplest hyperbolic 
knots~\cite{GZ23,Ga:237}) and the integer homology sphere obtained 
by $1/2$-surgery on the $4_1$-knot~\cite{Wheeler:QMC}.

Once introduced, quantum modular cocycles can be studied without 
knowledge of 3-dimensional topology and quantum field theory, 
and they are very interesting 
objects from many points of view, including their relation to 
$q$-series, $p$-adic analysis, and arithmetic geometry. What is
more, they are easy to 
construct using state-integrals. The latter are finite dimensional 
integrals of
products of a special function, the Faddeev quantum dilogarithm, 
evaluated at
integral linear forms times an exponential of a quadratic form. The 
integrands
of these state-integrals, and hence the corresponding cocycle, are 
determined 
by discrete combinatorial data. An example is the cocycle 
$\Omega_{A,B}$ defined below that depends on two integers~$A,B$.

The origin of this project was a search for algebraicity for 
the entries of
$\Omega_{\gamma}(\tau^*)$ where $\tau^* \in \BR$ is a fixed point 
of $\gamma$, i.e., satisfies $\gamma \tau^*=\tau^*$. 
Whereas we failed to recognise the entries of
this matrix explicitly, numerical experiments that predated 
the work of~\cite{RW} confirmed that the trace 
$\tr(\Omega_{\gamma}(\tau^*))$ is an 
algebraic number, and more generally, that the characteristic 
polynomial of 
$\Omega_{\gamma}(\tau^*)$ has algebraic coefficients.
It is important to note that the characteristic polynomial of 
$\Omega_\gamma(\tau^*)$ makes sense, since at these points 
$\gamma\tau^*=\tau^*$ implies that 
$\Omega_\gamma(\tau^*)$ behaves like an endomorphism 
(i.e., is a ``square matrix'').
These characteristic polynomials are expected to have coefficients 
in the abelian number field $\BK_\gamma$ generated by the values 
of the finite quantum dilogarithm, and some roots of unity 
(explicitly discussed below). 
We further found experimentally that 
$\mathrm{tr}(\Omega_{A,B,\gamma}(\tau^*))$ was periodic 
in $A$ and growing exponentially in $B$.

These experimental observations can be proved for the 
quantum modular cocycles associated with 1-dimensional 
state-integrals using a residue 
calculation reminiscent of the fundamental lemma
~\cite[Lem.2.1]{GK15} and of the proof of the functional 
equations satisfied by the finite quantum dilogarithms 
(see~\cite[Thm.~2]{RW}).

To state our main conjecture and our results, we begin with an 
elementary observation that computes the traces of modular 
cocycles evaluated at fixed points 
$\tau^*$ of $\gamma$.
Let $\Lambda_{\Omega,\gamma}$ denote the set of eigenvalues 
of the matrix $\Omega_{\gamma}(\tau^*)$.
Using the cocycle property, we have the following identity: 
\begin{equation}
\label{omg}
\tr(\Omega_{\gamma^n}(\tau^*)) = 
\sum_{\lambda \in \Lambda_{\Omega,\gamma}} \lambda^n \,.
\end{equation}
Moreover, both sides depend only on the conjugacy class of 
$\gamma$ in $\SL_2(\BZ)$. 
We can now formulate an algebraicity conjecture for these 
traces.

\begin{conjecture}
\label{conj.tralg}  
For every quantum modular cocycle $\Omega$, and every 
hyperbolic $\gamma$,  
\begin{itemize}
\item 
$\tr(\Omega_{\gamma^n}(\tau^*)) \in \calO_{\BK_\gamma}$,
hence (see Corollary~\ref{cor.AB}) 
$\Lambda_{\Omega,\gamma}$ is a set of units. 
\item 
If $\sigma\in\mathrm{Gal}(\overline{\QQ}/\QQ)$ satisfies 
$\sigma\eps=\eps^{-1}$ and 
$\lambda \in \Lambda_{\Omega,\gamma}$, then
\begin{equation}
\label{wep}
|\sigma\lambda|\in\eps^{\frac{1}{2}\ZZ}\,.
\end{equation} 
\end{itemize}
\end{conjecture}

In all cases that we know, the first part of the conjecture follows 
from explicit formulas for the trace given by sums over the group 
$G_\gamma=\BZ^2/(I-\gamma)\BZ^2$ of special values of the Faddeev 
quantum dilogarithm.
These formulas parallel the associated $q$-hypergeometric formulas 
with the $q$-Pochhammer symbol being replaced by the finite quantum 
dilogarithm $F_\gamma$ of~\cite{RW} (we take $F_\gamma:=F_\gamma^+$ 
in that notation).
Importantly, there are boundary terms that appear, which come from 
the failure of the five-term relation for $F_\gamma$ at lattice 
points.

\subsection{Examples from state-integrals}
\label{sub.13}

We will now give evidence for our conjecture by considering 
a family of cocycles $\Omega_{A,B}$ 
parametrised by a pair of integers $A$ and $B$. These cocycles
are expressed in terms of a special function, the Faddeev 
quantum dilogarithm~\cite{Fad95}, and are closely related to 
the $q$-hypergeometric series
\begin{equation}
\label{eq:fAB}
f_{A,B}(q)
=
\sum_{k=0}^{\infty}(-1)^{Ak}\frac{q^{Ak(k+1)/2}}{(q;q)_k^B}
\begin{pmatrix}
1\\    
q^k\\
\vdots\\
q^{(\max(A,B)-1)k}
\end{pmatrix}
\end{equation}
(also known as Nahm sums) which are 
vector-valued quantum modular forms for all $A,B\in\ZZ_{\geq 1}$. 
(This was proved in~\cite{Wheeler:QMC} for even $A$ and $B=1$, 
and the methods can be extended to prove these all 
of these  cases.)
To describe the associated cocycle, we define matrices $Q,R$ of 
size $A\times A$ 
and $(B-A)\times (B-A)$ so that 
$Q_{ij}=(-1)^{i+j-A-1}\binom{B}{i+j-A-1}$ and 
$R_{ij}=(-1)^{B-A-i-j+1}\binom{B}{B-A-i-j+1}$.
Then we take the matrix $P=(\smat{Q}{0}{0}{R})$, where we consider 
$Q$ 
as empty 
when $A\leq 0$ and similarly for $R$ when $B-A\leq 0$.
Finally, we define the rank-one matrix 
$D(z;\tau)_{ij}=\ee(jz+(i-1)z/\tau)$, where
\[
    \ee(x):=\exp(2\pi i x)\,,\qquad\ee_N(x):=\ee(x/N)\,.
\]
The values of $f_{A,B}$ at $q=\ee(\tau)$ and $\tq=\ee(-1/\tau)$ 
are related by 
$$
f_{A,B}(q)=\Omega_{A,B,S}(\tau)f_{A,B}(\tq)
$$
where the quantum modular cocycle $\Omega_{A,B,S}$ is a 
matrix-valued state integral given by
\[
\begin{aligned}
  \Omega_{A,B,S}(\tau)
  =
  \int_{\calC}
  \ee\Big(\frac{Az(z-\tau-1)}{2\tau}\Big)\;
  \Phi_{S,1,0}(z;\tau)^B
  \;D(z;\tau)\;P\;dz \,,
\end{aligned}
\]
where, for any $\gamma=(\smat{a}{b}{c}{d})\in\SL_2(\ZZ)$, we define 
$\Phi_{\gamma,m,n}$ to be the meromorphic continuation of
\begin{equation}
\label{Phimn}
\Phi_{\gamma,m,n}(z;\tau)
=
\frac{(\ee(z+m\tau);\ee(\tau))_\infty}{
(\ee(\tfrac{z}{c\tau+d}+n\gamma\tau);\ee(\gamma\tau))_\infty}\,.
\end{equation}
This cocycle makes sense for any $A,B$ with $B>0$.
There are associated quantum modular forms, however, when $A<0$, 
the formulas are not as neatly presented as that in 
Equation~\eqref{eq:fAB}.
Using the symmetrisation relation, which relates 
$\Phi_{\gamma,m,n}(z;\tau)^{-1}$ and $\Phi_{\gamma,m,n}(-z;\tau)$, 
these cover all cases with 
$B\neq0$ (which are trivially just Gaussian integrals and can be 
separately 
treated with $\theta$-functions).

With this explicit formula, one can then compute the trace to 
find that
\[
\tr(\Omega_{A,B,S}(\tau))
=
\int_{\calC}
  \ee\Big(\frac{Az(z-\tau-1)}{2\tau}\Big)\;
  \Phi_{S,1,0}(z;\tau)^B
  \;\frac{(1-\ee(z/\tau))^B\ee(Az)-(1-\ee(z))^B\ee(Az/\tau)}{1-\ee(z(1
  -\tau)/
  \tau))}\;dz\,.
\]
In Theorem~\ref{thm:cocycle} below we will show that 
$\Omega_{A,B,S}$ extends
to a quantum modular cocycle $\Omega_{A,B}$. Using the analogous 
formula for general $\gamma$ and the same contour 
manipulation in~\cite{RW}, 
we will prove the following, which confirms the first part of the 
main Conjecture~\ref{conj.tralg} 
for the quantum modular cocycle $\Omega_{A,B}$:

\begin{theorem}
\label{thm.nahm}
For all integers $A$ and $B$ with $AB(A-B)\neq0$ and all 
hyperbolic $\gamma$, we have (with the notation of 
Section~\ref{sub.gammaQDL}):
\begin{equation}
\label{trom} 
\tr(\Omega_{A,B,\gamma}(\tau^*)) =
\frac{1}{\sqrt{|G_\gamma|}}
\sum_{k \in G_{\gamma}} \langle k\rangle_{\gamma}^A F_{\gamma}(k)^B 
\,.
\end{equation}
\end{theorem}
\noindent 
This and the fact that $F_\gamma(k)$ are algebraic units 
~\cite{RW:stark} implies the first part of 
Conjecture~\ref{conj.tralg} for $\Omega_{A,B}$ (see 
Corollary~\ref{cor.AB}).
Equation~\eqref{trom} also implies that the left hand side 
is a periodic function of $A$, as was experimentally observed. 

Combining ~\eqref{omg} and~\eqref{trom}, we obtain 
\begin{equation}
\label{tromAB}
\sum_{\lambda \in \Lambda_{\Omega_{A,B},\gamma}} \lambda^n
 = \frac{1}{\sqrt{|G_{\gamma^n}|}}
\sum_{k \in G_{\gamma^n}} \langle k\rangle_{\gamma^n}^A F_{\gamma^n}
(k)^B, 
 \qquad (n \neq 0) \,.
\end{equation}
This identity where $\Lambda_{\Omega_{A,B},\gamma}$ is a finite set 
of units is by no means obvious since the summand of 
Equation~\eqref{trom} 
lies in $\BK_{\gamma^n}$, and as $n$ increases, so does the field. 
It does follow from the results of~\cite{RW:stark} that the right hand side of the above 
equation lies in a fixed number field independent of $n$. Even so, 
integrality of the right hand-side of the above equation is still
not obvious, and its shape is reminiscent of averages of 
the cyclic quantum dilogarithm (see~\cite[Eqn.(12)]{DG2} and 
\cite[Eqn.(7.6)]{GZ24}), now known to be algebraic 
integers~\cite[Cor.3.12]{GSWZ}.

The fact the set $\Lambda_{\Omega_{A,B},\gamma}$ consists of units is 
interesting but they seem to satisfy an even more interesting number 
theoretic property.
From the rank-one Stark conjectures~\cite{RW:stark}, it is known 
that a 
Galois automorphism $\sigma\in\mathrm{Gal}(\overline{\QQ}/\QQ)$ such 
that $\sigma(\eps)=\eps^{-1}$ satisfies $|\sigma(F_\gamma(u))|=1$ for 
$u\neq0$. This implies that the right hand side of~\eqref{trom} is
at most $\sqrt{|G_{\gamma^n}|}$ which in turn implies that
an upper bound on the size of the corresponding $\sigma(\lambda)$ 
given by $|\sigma(\lambda)| \leq \varepsilon^{1/2}$. Combining this
with the fact (see Remark~\ref{rem:recip}) that 
if $\lambda\in\Lambda_{\Omega_{A,B},\gamma}$, then $
\lambda^{-1}\in\Lambda_{\Omega_{A,B},\gamma}$ we conclude that 
absolute values that $\sigma(\lambda)$ lie in the critical strip
\begin{equation}
\label{eq.strip}
\varepsilon^{-1/2} \leq |\sigma(\lambda)| \leq \varepsilon^{1/2}\,.
\end{equation}
Experimenting further leads to the 
conjecture that in fact $|\sigma(\lambda)|=1$.

\begin{example}
Consider $A=1,B=2$ (an example associated to the figure eight
knot) and 
$\tau^*=\tfrac{3+\sqrt{5}}{2}$ where $\gamma=(\smat3{-1}
10)$.
We find that
\[
\Omega_{1,2,\big(\smat3{-1}10\big)}(\tfrac{3+\sqrt{5}}
{2})
\approx
\begin{pmatrix}
2.1781 + 0.13857i & -0.27714i\\
0.86909 - 0.081504i & 0.43993 - 0.13857i
\end{pmatrix}
\,.
\]
We could not recognise these numbers, however, from 
Theorem~\ref{thm.nahm}, we can explicitly compute the 
characteristic 
polynomial
\[
\lambda^2-\eps \lambda+1\,,\quad\text{with roots}\quad
\lambda=0.46431\cdots,\;2.1537\cdots\,.
\]
Exchanging $\eps\to\eps^{-1}$ gives
\[
\lambda^2-\eps^{-1}\lambda+1\,,\quad\text{with roots}\quad
\lambda=0.19098\cdots \pm 0.98159\cdots i\,.
\]
These have absolute value equal to $1$, which can be 
explicitly proved.
\end{example}

\begin{example}
Consider $A=4,B=1$ and $\tau^*=\tfrac{5+\sqrt{21}}{2}$ 
where $\gamma=(\smat5{-1}10)$.
Then $\Omega_{4,1,\gamma}(\tau^*)$ is (approximately) equal to
\[
\begin{pmatrix}
1.7797 + 0.48236i & 0.0056762 - 0.24307i & -0.28531 - 
0.25524i 
& -1.0837 + 1.1034i\\
1.5074 - 0.99185i & 0.21566 + 0.16402i & -0.056410 - 0.29246i 
& 
-0.59748 + 1.1575i\\
0.085692 - 0.46996i & -0.31378 - 0.024085i & 0.16794 - 
0.18752i 
& -0.33088 + 0.91402i\\
0.51746 + 0.68632i & 0.37349 - 0.045816i & 0.21516 + 
0.087157i 
& -0.26765 + 0.63560i
\end{pmatrix}\,.
\]
Again we can prove that the characteristic polynomial is given by
\[
\lambda^4+a\lambda^3-a\lambda^2+a\lambda+\ee(1/6)\,,
\quad\text{where }a\approx-1.8956 - 1.0945\,
\mathrm{i}\quad\text{and}\quad a^2 + (\eps - 1)a + \eps=0\,.
\]
The roots of this polynomial are
\[
\ee(1/21)\,,\; \ee(1/6)\,,\; \ee(4/21)\,,\; \ee(-5/21)\,.
\]
These are sent to other roots of unity under the action of a Galois 
automorphism. In this example the units $\lambda$ do lie in an
abelian extension of $\BQ$.
\end{example}

\begin{example}
Consider $A=4,B=1$ and $\tau^*=\tfrac{3+\sqrt{5}}{2}$ 
where $\gamma=(\smat3{-1}10)$.
Then $\Omega_{4,1,\gamma}(\tau^*)$ is (approximately) equal to
\[
\begin{pmatrix}
0.71134 + 0.23773i & -0.33237 + 0.035057i & -0.26503 + 0.26890i & 
0.33671 - 0.87669i\\
1.2549 + 0.13461i & 0.36464 - 0.12931i & -0.28795 - 0.14590i & 
0.30596 - 0.56211i\\
1.1139 - 0.24006i & -0.46853 - 0.050466i & 0.21590 - 0.23370i & 
0.37531 - 0.24302i\\
0.78389 - 0.22224i & 0.18558 + 0.55438i & -0.033550 + 0.34892i & 
0.32615 + 0.12528i
\end{pmatrix}\,.
\]
One can prove that the characteristic polynomial is given by
\[
\lambda^4-\tfrac{1+\sqrt{5}}{2}\lambda^3+\lambda^2-\tfrac{1+\sqrt{5}}
{2}\lambda+1\,.
\]
In this example the units $\lambda$ do not lie in an abelian 
extension of $\BQ(\tau^{*})$.
Their Galois conjugates satisfying 
\[
\lambda^4-\tfrac{1-\sqrt{5}}{2}\lambda^3+\lambda^2
-\tfrac{1-\sqrt{5}}{2}\lambda+1=0
\]
all lie on the unit circle and are approximately equal to
\[
    -0.67784 \pm 0.73521i\,,\quad  0.36882 \pm 0.92950i\,.
\]
\end{example}

\subsection{Traces from finite TQFTs}
\label{sub.Cgamma}

We already discussed a source of quantum modular cocycles, namely
the ones $\Omega_M$ that are conjecturally associated to 3-manifolds
$M$~\cite{GZ24}. The simplest invariant of such a cocycle 
is its rank, and in~\cite{GZ24} it was conjectured that the rank 
of $\Omega_M$ should equal to the rank $\Sk_M$ of its Kauffman 
bracket skein module introduced independently by Przytycki and 
Turaev~\cite{Przytycki,Tu:conway}, and proven to be of finite
rank by~\cite{Gunningham-Jordan}.
The latter is also conjectured to equal to the degree 0 part of the
$\SL_2(\BC)$ Floer cohomology $\HF(M)$ of 
Abouzaid--Manolescu~\cite{AM20}, altogether giving
\begin{equation}
\label{3e}
\rk \, \Omega_M \stackrel{?}= \rk \, \Sk_M \stackrel{?}
= \rk \, \HF^0(M) \,. 
\end{equation}

In this section we discuss an alternative and rather surprising 
way to compute the traces of $\Omega_M$ from invariants 
associated by a finite TQFT denoted $\ZTV_\gamma$ that depends on a 
hyperbolic element $\gamma$ of $\SL_2(\BZ)$. 
Conjecture~\ref{conj.TQFT} below
\begin{itemize}
\item
suggests a $\BZ/2\BZ$-graded extension of the 
cocycle $\Omega_{M,\gamma}$ and of the skein module
$\Sk_M$ whose rank is the even and odd rank of $\HF(M)$, both
with a $\BZ$-grading coming from the $\varepsilon^{1/2}$-weight of 
Equation~\eqref{wep}. 
\item 
gives a way to compute the ranks and the traces of the extended 
cocycle and skein module.
\item
the finite TQFTs take values in the field $\BK_{\gamma}$ 
and leads to surprising identities among algebraic numbers that
currently have not be explained by the explicit class field theory 
of~\cite{RW:stark}.
\end{itemize}

The above mentioned TQFT comes from fusion categories $\calC$
studied by Izumi~\cite{Izu01} whose Grothendieck ring has
objects $[g]$, $g \in G$ and $[\rho]$ satisfying the fusion rules 
\begin{equation}
\label{frules}
[g][h]=[g+h], \qquad [g][\rho]=[\rho][g]=[\rho], \qquad 
[\rho]^2=N[\rho]
+\sum_{g \in G} [g] 
\end{equation}
for certain metric abelian groups $G$. Izumi gave a construction of 
such categories from a solution to a system of quadratic and 
cubic equations~\cite[Thm.5.3]{Izu01} and questioned whether such solutions exist. 

In~\cite{RW}, two of the authors proved that the functions $F
_\gamma$
are a solution to that system of equations, thus proving the 
existence 
of such 
fusion categories $\calC_\ga$ parametrised by a hyperbolic elements 
$\ga \in \SL_2(\BZ)$ (i.e., $\tr(\gamma)=N_\gamma+2>2$), with corresponding
abelian group $G_\ga=\BZ^2/\BZ^2(\ga-I)$ of order $N_\gamma$.\footnote{We 
write vectors in row form, so 
$\BZ^2=\{(a,b) \, | \, a, b \in \BZ \}$.}

These categories, which depend only on the equivalence class of
a hyperbolic element $\gamma$ in $\SL_2(\BZ)$, realise all finite 
abelian groups $G$ that
are products of two abelian groups, and provide an infinite class of 
fusion categories
not defined over a cyclotomic field answering a question 
of~\cite[p. 591]{ENO05}. Morrison--Snyder \cite{Morrison} gave two 
examples of non-cyclotomic fusion categories. 

Using the calculus of Leavitt algebras, we compute the 
$6j$-symbols of the categories $\calC_\gamma$ (see 
Proposition~\ref{prop.6j} below) and use them to construct a 
Turaev--Viro type TQFT denoted $\ZTV_\gamma$ based on $\calC_\gamma$.
Recall that $\ZTV_\gamma(M)$ 
is the corresponding invariant of a 3-manifold $M$. We will denote 
$\ZTVred(M)=\ZTV_{\gamma}(M) \,\ZTV_{\gamma}(S^3)^{-1}$.

\begin{conjecture}
\label{conj.TQFT}
For all 3-manifolds $M$ there exists a $\ZZ/2\ZZ$-graded quantum modular cocycle $\Omega_M$ and for all hyperbolic elements $
\gamma\in\SL_2(\BZ)$, we have:
\begin{equation}
\ZTVred_{\gamma}(M) = \tr(\Omega_{M,\gamma}(\tau^*)) \,.
\end{equation}
This implies that there exist finite sets
of algebraic units $\Lambda_{M,\gamma}^\pm$ (whose cardinality is
independent of $\gamma$) so that for all integers $n$, 
we have:
\begin{equation}
\label{ztvl}
 \ZTVred_{\gamma^n}(M) = 
 \sum_{\lambda \in \Lambda^+_\gamma(M)}  \lambda^n 
 -\sum_{\lambda \in \Lambda^-_\gamma(M)}  \lambda^n\,.
\end{equation}
\end{conjecture}

Note that the invariant for the three sphere is 
$\ZTV_\gamma(S^3)
\varepsilon^2(\varepsilon-1)^{-2}(\varepsilon+1)^{-1}$.
Thus, if $\ZTV_\gamma(M)$ satisfies Equation~\eqref{ztvl}, so 
does $\ZTV_\gamma(M) \ZTV_\gamma(S^3)^{-1}$.

\begin{remark}
\noindent
1. Unlike Conjecture~\ref{conj.TQFT}, Equation~\eqref{ztvl} is 
explicitly
computable, numerically or exactly, and hence can be confirmed or 
falsified.

\noindent
2. 
If we think of $\gamma^n$ as some
version of Planck's constant, Equation~\eqref{ztvl} predicts that 
the numerical invariants $\ZTV_{\gamma^n}$ as $n$ varies are 
exact (i.e., there are power series in $1/n$ involved unlike 
the case of perturbative quantum field theory). 

\noindent
3.
Arithmetically, Equation~\eqref{ztvl} is similar to the Weil-type 
counts (of some realization of a not-yet-discovered motive) of the 
number of points of varieties over finite fields~\cite{Deligne:weil1}.

\noindent
4.
Equation~\eqref{ztvl} is by no means obvious even for the explicitly 
given 
invariant of the $4_1$  (see Equation~\eqref{trom}). 
\end{remark}

\subsection{The Drinfeld center of \texorpdfstring{$\calC_\gamma$}{Cgamma}}
\label{sub.ZCgamma}

We now discuss the Drinfeld center of the fusion category 
$\calC_\gamma$. Recall that 
every fusion category $\calC$ has a Drinfeld center $\calZ(\calC)$, 
which is a
braided monoidal category, \cite[Sec.7.13]{ENO05}, and in fact a 
modular category by 
a fundamental theorem of M\"uger~\cite{Muger}. 
The objects in $\calZ(\calC)$ consist of pairs $(X,e_X)$ of 
objects $X$ 
of $\calC$
and half-braidings $e_X$, i.e., collections of invertible maps
$e_X(Y) : X \otimes Y \stackrel{\simeq}{\to} Y \otimes X$ for every 
object 
$Y$ of $\calC$ that satisfy certain hexagon relations. 

Fusion and modular categories can be used to construct TQFTs via two 
general methods: 
\begin{itemize}
\item 
the Reshetikhin--Turaev approach uses as input a modular category and 
presents
manifolds by surgery on framed links in $S^3$~\cite{Turaev},
\item 
the Turaev--Viro approach uses as input a spherical category and 
presents 
manifolds via triangulations~\cite{Barrett-Westbury}.
\end{itemize}
The two approaches lead to the same invariants of 3-manifolds and in 
fact to
equivalent TQFTs if one uses a spherical category 
$\calC$ on the TV side and its Drinfeld center $\calZ(\calC)$ 
on the RT side.
For a detailed discussion, see Turaev--Virelizier
~\cite{Turaev-Virelizier} 
and also Balsam--Kirillov~\cite{Balsam-Kirillov}. An advantage of
having an explicit description of the Drinfeld center $
\calZ(\calC_\gamma)$ is that it leads to explicit representations
of the mapping class groups, and in particular to invariants of 
links, much in the spirit of Jones. 

Towards this task, 
Izumi introduced a second system of quadratic equations 
\cite[Eqns.(6.17)-(6.21)]{Izu01} over a solution to
his first system of quadratic and cubic equations, and using 
a solution of this
second system, he constructed the half braiding of the objects of the 
fusion 
categories that he considered. He questioned whether a solution 
to this 
second 
system of equations existed.

\begin{theorem}
\label{thm.second}
When $\tr\gamma$ is odd, the function $F_{\ga^2}$ solves the second system of equations 
~\eqref{eq:izumigamma_4}, \eqref{eq:izumigamma_1}, \eqref{eq:izumigamma_3}, and~\eqref{eq:izumigamma_2}
and constructs the modular category $\calZ(\calC_\gamma)$. 
\end{theorem}

As a corollary of Izumi's work and the above theorem, we find that 
$\calZ(\calC_\gamma)$ has four types 
of simple objects
$\mathfrak{a}_g, \mathfrak{b}_g, \mathfrak{c}_{g,h},\mathfrak{d}_j$
indexed by $g \in G_\gamma$ for the first two types, unordered pairs 
$(g,h)$
of different elements of $G_\gamma$ for the third type and the set 
$J_\gamma=(G_{\gamma^2}\sm G_\gamma)/\gamma$ of $N_\gamma(N_\gamma+3)/2$ 
elements. The images of the simple objects under the forgetful 
functor 
$U: \calZ(\calC_\ga) \to \calC_\ga$ mapping 
$(X,e_X)$ to $X$ 
are given by 
$
U(\mathfrak{a}_g) = [g], 
U(\mathfrak{b}_g) = [g] + [\rho], 
U(\mathfrak{c}_{g,h}) = [g] + [h] + [\rho],
U(\mathfrak{d}_j) = [\rho] \,.
$
Moreover, the $S$ and $T$ matrices of the $\SL_2(\BZ)$-representation 
of the modular category $\calZ(\calC_\gamma)$ are 
defined over $\BQ(\zeta_{N_\gamma(N_\gamma+4)})$ where $|G_{\gamma^2}|
=N_\gamma(N_\gamma+4)$. 
The $S$ and $T$ matrices have a block decomposition corresponding 
to the four types of simple modules of the Drinfeld center. The 
$T$-matrix is diagonal and all but the last block of the $S$ matrix
is expressed in terms of the bilinear forms on $G_\gamma$ and 
$G_{\gamma^2}$. 

With the notation of Equation~\eqref{Gconstants}, we have the 
following.

\begin{theorem}
\label{thm.S}
When $\tr\gamma$ is odd, 
the interesting part of the $S$ matrix is
\begin{equation}
\label{Sint}
S_{\fd_j,\fd_\ell}
=
2\frac{\Re\big(\ang{j;\ell}_{\gamma^2}\ang{\varpi(j);
\varpi(\ell)}
_{\gamma}^{-1/2}\big)\ang{\varpi(j);\varpi(\ell)}_\gamma^{1/2}}
{\sqrt{N_\gamma(N_\gamma+4)}}
=
\frac{\ang{j;\ell}_{\gamma^2}+\ang{j;\ell}_{\gamma^2}
^{-1}\ang{\varpi(j);\varpi(\ell)}_{\gamma}}
{\sqrt{N_\gamma(N_\gamma+4)}}\,,
\end{equation}
where $j, \ell \in J_\gamma$ and $\varpi:G_{\gamma^2}\to 
G_\gamma$ is the 
projection map induced from $\mathrm{im}(\gamma^2-
I)\subseteq\mathrm{im}
(\gamma-I)$.
\end{theorem}
For completeness, the full $S$ and $T$ matrices are given in
Remark~\ref{rem.ST}.

Finally, this gives rise to explicit representations of the 
braid groups discussed in Section~\ref{sub.Rmat} below.

\subsection{Evidence for Conjecture \ref{conj.TQFT}}

In this last section we discuss evidence for Conjecture 
\ref{conj.TQFT}.

\noindent$\bullet$
{\bf The Poincare homology sphere $S(2,3,5)$.}
Recall that the Brieskorn spheres $S(p,q,r)$ for positive integers 
$p,q,r$ are integer homology spheres that are Seifert fibered 
3-manifolds obtained by surgery on a Hopf link, with one of its 
components tripled, with appropriate framing. Their 
$\SL_2(\BC)$-character variety consists of the trivial
representation and 
$$
r_{p,q,r}:=\tfrac{1}{4}(p-1)(q-1)(r-1)
$$
irreducible representations~\cite{BC06}. The rank of the framed 
Floer cohomology $\HF$ is given by~\cite[Thm.1.6]{AM20}
\begin{equation}
\label{HFB}
\rk \, \HF^0(S(p,q,r)) = 1 + r_{p,q,r}, \qquad 
 \rk \, \HF^{-3}(S(p,q,r)) =  r_{p,q,r} \,.
\end{equation}
The rank of their $\SL_2(\BC)$-skein module is also $1 + r_{p,q,r}$
~\cite{DKS25} confirming the conjecture that
\begin{equation}
\label{conj.sk}
\rk \, \Sk(M) = \rk \, \HF^0(M) \,.
\end{equation}

Since Seifert fibered manifolds have a surgery presentation
in terms of a chain of Hopf links arranged in a tree, it follows that
their TQFT invariants can be computed in terms of the $S$ and $T$
matrices alone. In particular, for the Poincare homology
sphere $S(2,3,5)$, we have

\begin{equation}
\label{Z235}
\ZTV(S(2,3,5))=\sum_{m} \frac{(T_{mm})^{-1}}{S_{0m}} 
\Big( 
\sum_{a} (S T^2 S)_{m0}
\Big) \Big( 
\sum_{b} (S T^3 S)_{m0}
\Big) \Big( 
   \sum_{c} (S T^5 S)_{m0}
\Big) \,.
\end{equation}

Using this with $\gamma=\left(\smat{N+2}{-1}{1}{0}\right)$ and $N=1$
and the values of the $S$ and $T$ matrices we compute
for  $\ZTVred_{\gamma^n}(S(2,3,5))$ for $n=1,2,3,4,5$:
\begin{equation}
8 - 3 \varepsilon, \,\, -11 + 3 \varepsilon,  \,\, 
-55 + 24 \varepsilon,  \,\, 
610 - 231 \varepsilon,  \,\,
-3189 + 1210 \varepsilon \,.
\end{equation}
Using the first four values and the fact that the invariant is
real, we found out an interpolation of the form 
\begin{equation}
\label{eq.S235}
\ZTVred_{\gamma^n}(S(2,3,5)) = 1^n 
+ (\varepsilon^{1/2} \zeta_{10}^{-3})^n
+ (\varepsilon^{1/2} \zeta_{10}^{3})^n
- (\varepsilon^{-3/2} \zeta_{10}^{-3})^n
- (\varepsilon^{-3/2} \zeta_{10}^{3})^n
\end{equation}
for $n=1,2,3,4$, where $\zeta_m=\ee(1/m)$.
It was remarkable that this involved five $\lambda$, three coming
with positive sign and two with negative. Having found this
interpolation, we then computed the invariant for $n=5$ and found
an agreement with the interpolated answer. This gives us confidence 
that Equation~\eqref{eq.S235} is valid for all $n>0$. (Presumably,
this can be proven using formulas for Gauss sums.) The conclusion
is that the five $\lambda$ have $\BZ/2\BZ$-parity and integer weight 
(corresponding to the power of $\varepsilon^{1/2}$).  
Their number, with fixed parity, agrees with the rank of the even 
and odd $\HF$. 

Note that the size (i.e., number of rows) of the $S$ and $T$ 
matrices for $n=1,\dots,5$
are given by $5, 45, 320, 2205, 15125$, a sequence that grows 
exponentially with respect to $n$, and the sum~\eqref{eq.S235} 
requires $O(n^2)$-computations.

\noindent$\bullet$
{\bf Surface bundles over a circle.}
We next compute the invariants of the TQFT for trivial 
surface bundles over a circle. 
The formula suggests clearly the the need to consider 
$\BZ/2\BZ$-graded cocycles. 

\begin{proposition}
\label{prop.Sigmag}
Fix $\gamma$ hyperbolic. For every $g \geq 1$, there exist finite 
subsets $\Lambda_{\gamma,g}^+$ and $\Lambda_{\gamma,g}^-$ 
of $\varepsilon^\BZ$
such that for all nonzero integers $n$ we have 
\begin{equation}
\ZTVred_{\gamma^n}(S^1 \times \Sigma_g) = 
\sum_{\lambda \in \Lambda_{\gamma,g}^+} \lambda^n 
- \sum_{\lambda \in \Lambda_{\gamma,g}^-} \lambda^n \,.
\end{equation}
\end{proposition}

\noindent
Moreover, $|\Lambda_{\gamma,g}^+|=|\Lambda_{\gamma,g}^-|$ is 
independent of $\gamma$. For example, for $g=1,2,3$ and
with $s=\varepsilon^n$, we have:
\begin{equation}
\begin{aligned}
\ZTVred_{\gamma^n}(S^1 \times \Sigma_1) &= 
s^{-4} - 2 s^{-3} + s^{-1} + 1 - 2 s^2 + s^3, \\ 
\ZTVred_{\gamma^n}(S^1 \times \Sigma_2) &=  s^{-7} 
- 4 s^{-6} + 2 s^{-5} + 20 s^{-4} - 60 s^{-3} + 79 s^{-2} 
- 38 s^{-1} -38 + 79 s \\ & - 60 s^2 + 20 s^3 + 2 s^4 - 
  4 s^5 + s^6,  \\
\ZTVred_{\gamma^n}(S^1 \times \Sigma_3) &=
  s^{-11} - 7 s^{-10} + 16 s^{-9} + s^{-8} - 67 s^{-7}
  + 119 s^{-6} - 111 s^{-5} + 223 s^{-4} 
  - 596 s^{-3}  + 838 s^{-2} \\ & - 417 s^{-1} -417 + 838 s - 
  596 s^2 + 223 s^3 - 111 s^4 + 119 s^5 - 67 s^6 + s^7 + 16 s^8 - 
  7 s^9 + s^{10} \,.
 \end{aligned}
\end{equation}
  
It is not clear how the cardinality of the set 
$\Lambda_{\gamma,g}^{\pm 1}$ 
(which are independent of $\gamma$ and depend on $g$ alone) is related
to the rank $2^{2g+1}+2g-1$ of the $\SL_2(\BC)$-skein module on 
$S^1 \times \Sigma_g$; see Gilmer--Masbaum~\cite{GM19}. 
Experimentally, we have
$|\Lambda_{\gamma,g}^{\pm 1}|= \tfrac{5}{2}\cdot 2^{4g-4} + 
O(9.5^g)$, 
and the first few
values of $|\Lambda_{\gamma,g}^{\pm 1}|$ for $g=1,\dots,10$ are given 
by
\begin{equation}
4, \, 204, \, 2396, \, 30858, \, 459026, \, 6872676, \, 
98425386, \, 1487300558, \, 22405945978, \, 340039248620 \,.
\end{equation}

\noindent$\bullet$
{\bf Lens spaces.} We now discuss Lens spaces. Recall that they
are obtained by rational surgery on the unknot, and their TQFT 
invariants are determined by the $S$ and $T$ matrices. In particular,
for positive integer $p$, we have
\begin{equation}
\label{TQFTLens}
\ZTV(L(p,1)) = (S T^p S)_{0,0} \,.
\end{equation}

With the notation of Section~\ref{sub.finiteQDL}, for 
every hyperbolic element $\gamma$ and every $p>0$, we have:
\begin{equation}
\label{ZLens}
  \ZTVred_\gamma(L(p,1))
 =      \frac{1}{2}(1+\varepsilon^{-1})
 \big|G_\gamma[p]\big|
+ \frac{1}{2}(1-\varepsilon^{-1}) \cdot 
\frac{1}{\sqrt{|G_{\gamma^2}|}}
\sum_{j\in G_{\gamma^2}}
\frac{\ang{j}_{\gamma^2}^p}{\ang{\varpi j}_\gamma^p} 
\end{equation}
where $G[p]$ denotes the $p$-torsion subgroup of $G$.
Using this with $\gamma=\left(\smat{N+2}{-1}{1}{0}\right)$ with $N=1$
we compute  $\ZTVred_{\gamma^n}(L(2,1))$ for $n=1,\dots,9$:
\begin{equation}
3 - \varepsilon, \, 8 - 3 \varepsilon, \, 24 - 8 \varepsilon, 
\, 55 - 21 \varepsilon, \, 144 - 55 \varepsilon, 
\, 380 - 144 \varepsilon, \, 987 - 377 \varepsilon, \,
2584 - 987 \varepsilon, \, 6768 - 2584 \varepsilon \,.
\end{equation}
This data interpolates to
\begin{equation}
\ZTVred_{\gamma^n}(L(2,1)) = 1^n + (\zeta_3^{-1})^n 
+ (\zeta_3)^n + (\varepsilon^{-1})^n
\end{equation}
for $n=1,\dots,9$ suggesting that the even part of the conjectured
motive has rank 4. We repeated the computations and found that
\begin{equation}
\label{L21lambdas}
\Lambda^+_{\left(\smat{N+2}{-1}{1}{0} \right)}(L(2,1)) = 
\begin{cases} \{1, \, \zeta_3,  \, \zeta_3^{-1},  \,
\varepsilon^{-1} \}
& \text{if $N$ is odd} \\
 \{1,  \, 1,  \, -1,  \, \varepsilon^{-1} \}
& \text{if $N>0$ is even}
\end{cases}, \qquad
\Lambda^-_{\left(\smat{n+2}{-1}{1}{0} \right)}(L(2,1)) 
= \emptyset \,.
\end{equation}
Next for $L(3,1)$ and $\gamma=\left(\smat{N+2}{-1}{1}{0}\right)$ with 
$N=1$ we found the values
\begin{equation}
3 - \varepsilon, \, -6 + 3 \varepsilon, \,
21 - 8 \varepsilon, \, 171 - 63 \varepsilon, \, 
144 - 55 \varepsilon, \, -375 + 144 \varepsilon, 
\, 987 - 377 \varepsilon, \, 7758 - 2961 \varepsilon, 
6765 - 2584 \varepsilon
\end{equation}
for $n=1,\dots,9$ which interpolate to 
\begin{equation}
\ZTVred_{\gamma^n}(L(3,1)) = 1^n + 1^n +
(-1)^n +  (-1)^n + i^n + (-i)^n + (\varepsilon^{-1})^n 
+ (i \varepsilon^{-1})^n + (-i \varepsilon^{-1})^n
\end{equation}
for $n=1,\dots,9$, suggesting that the conjectured
motive has even rank 9 and odd rank 0. Going the the next 
example with $L(4,1)$ we find the values
\begin{equation}
1, \, 1, \, 136 - 48 \varepsilon, \, 1, \, 1, \, 2272 - 864 
\varepsilon, \, 1, \, 1, \, 40600 - 15504 \varepsilon, \, 1, \, 1
\end{equation}
for $n=1,\dots,11$ which interpolates to 
\begin{equation}
\ZTVred_{\gamma^n}(L(4,1)) = 6 \cdot 1^n 
+ 3 \cdot (\zeta_3^n+(\zeta_3^{-1})^n) + 2 \cdot 
((\varepsilon^{-1} \zeta_3)^n+(\varepsilon^{-1}\zeta_3^{-1})^n)
\end{equation}
for $n=1,\dots,11$, suggesting that the conjectured
motive has even rank 16 and odd rank 0. 

Summarising, the data for $L(2,1)$, $L(3,1)$ and $L(5,1)$ 
suggests that Equation~\eqref{ztvl} holds for $L(p,1)$ with 
$|\Lambda^+_{\gamma,L(p,1)}|=p^2$ and 
$|\Lambda^-_{\gamma,L(p,1)}|=0$.
This result differs from the rank of the 
framed $\SL_2(\BC)$ Floer Homology~\cite[Thm.1.5]{AM20} 
\begin{equation}
\rk \, \HF^0(L(p,q)) = 1 + \Big\lfloor \frac{p-1}{2} \Big\rfloor,
\qquad \rk \, \HF^{-2}(L(p,q)) = \Big\lfloor \frac{p-1}{2} \Big\rfloor
\end{equation}
and with the rank $1 + \lfloor \frac{p-1}{2}\rfloor$ of the
Kauffman bracket skein module; see Hoste--Przytycki~\cite{HP93}.

\noindent$\bullet$
{\bf The complement of the $4_1$ knot.}
Our next result is a proof of Conjecture~\ref{conj.TQFT} for 
the simplest
hyperbolic knot complement.

\begin{theorem}
\label{thm.41}
Conjecture~\ref{conj.TQFT} holds for the complement of the $4_1$ knot.
\end{theorem}

\subsection*{Acknowledgments} The authors wish to thank 
Dustin Clausen, Maxim Kontsevich, Doron Lubinsky, Ciprian Manolescu
and Don Zagier for enlightening conversations. C.W. was supported 
by the Institut des Hautes \'Etudes Scientifiques, France.
C.W. wishes to thank the Southern University of Science and 
Technology’s International Center for Mathematics in Shenzhen 
for their hospitality.
D.R. acknowledges funding by the European Union (ERC, FourIntExP, 
101078782).

\smallskip
\noindent \textbf{AI use disclosure.}
We did not make use of LLMs in the preparation and writing of this paper.

\section{\texorpdfstring{Solutions of $q$-difference equations at $|q|=1$}{Solutions of q-difference equations at |q|=1}}
\label{sec:|q|=1}

In this self-contained section, which may be of independent interest,
we construct solutions of linear $q$-difference equations
with coefficients in $\BQ[t^{\pm 1},q^{\pm 1}]$ when $q \in \BC$ 
satisfies $|q|=1$. We will use the existence of such solutions in
the next section to prove that the meromorphic functions 
$\Omega_{A,B,\gamma}$ indeed satisfy the cocycle property.

There is a vast literature on the subject of linear $q$-difference 
equations, and their formal power series solutions, and
their actual solutions when $|q|<1$. 

Let us recall some basic facts. All $q$-difference modules $M$ over 
$\CC^\times$ (that is, modules over the $q$-Weyl algebra 
$\BQ(q)[t^{\pm1}]\langle \sigma \rangle/(\sigma t - q t \sigma)$) 
have cyclic 
vectors, hence their solution space can be expressed in terms of 
solutions of a linear $q$-difference equation
\begin{equation}
\label{qdiffeq}
\sum_{j=0}^d c_j(t,q) f(q^j t, q)=0, \qquad c_j(t,q) \in \BZ[t^{\pm 
1},
q^{\pm 1}]  \quad j=0,\dots,d\,.
\end{equation}
The Newton polygon of the above $q$-difference module $M$ is the 
convex 
hull 
of the polynomial $\sum_{j=0}^d c_j(t,q) \sigma^j$ in the plane with 
coordinates 
$(\sigma,t)$. The edges of the lower (resp., upper) part of the 
Newton 
polygon
define two increasing filtrations on the formal solution space of $M$
at $t=0$ and $t=\infty$. 

\begin{definition}
We say that a $q$-difference module $M$ over $\CC^\times$ 
has quasi-unipotent monodromy if the roots of its edge polynomials 
are in $\mu_\infty q^{\BQ}$,
where $\mu_\infty$ denotes the set of complex roots of unity.
\end{definition}

This is a natural condition, reminiscent to the fact that local 
monodromy of
varieties defined over $\BQ$ is quasi-unipotent. In fact, all $q$-difference 
modules that appear naturally in quantum topology (see, e.g., 
~\cite{GL:qholo}) 
are conjectured to be quasi-unipotent. 

\begin{theorem}
\label{thm.q=1}
Every $q$-difference module over $\CC^\times$ has a fundamental 
basis of solutions whose 
power series have positive radius of convergence for all $q$ 
in the unit circle 
of Hausdorff dimension 0 and logarithmic dimension 2 
(in particular, almost
everywhere on the unit circle).
\end{theorem}

Some related results were also studied by di Vizio~\cite{Div09}.

\begin{proof}
The Frobenius method uses the Newton polygon of $M$ to 
construct a basis 
of formal solutions of the form
\begin{equation}
\label{ff}
   t^{\alpha \tfrac{\log t}{\log q}} t^{\log_q \rho} (\log q)^m 
   f(t,q), 
   \qquad f(t,q)=\sum_{k=0}^\infty a_k(q) t^k
\end{equation}
that depend on the slope of an edge of the Newton polygon, 
and a solution 
$\rho$ of its edge polynomial, where $m$ is strictly less 
than the multiplicity
of $\rho$ in its edge polynomial.

The method, described in detail for instance 
in~\cite[Sec.2.2]{GW22}, uses 
one edge at a time, which after shearing, we can assume that 
it is horizontal.
To a horizontal edge, one associates an edge polynomial, whose 
roots $\rho$ 
determine the leading terms of a formal power series~\eqref{ff} 
with $\alpha=0$.

For example, Equation~\eqref{qdiffeq} has a horizontal edge 
$[0,d] \times \{0\}$ 
when $c_j(t,q) \in \BZ[t,q]$ for all $j=0,\dots,d$ and
$c_0(0,q) c_d(0,q) \neq 0$.
The edge polynomial is 
$p_e(\sigma,q)=\sum_{j=0}^d c_j(0,q) \sigma^k=c(d)
\prod_{\rho'} (\sigma-\rho')$. 

Assume first that $\rho$ is a simple root. 
Substituting this ansatz in Equation~\eqref{qdiffeq},
we obtain a recursion 
relation for $a_k(q)$ of the form
\begin{equation}
\label{arec}
p_e(\rho q^k,q) a_k(q) 
= c'_{-1}(q^k,q) a_{k-1}(q) 
+ \dots + c'_{-s}(q^k,q) a_{k-s}(q)
\end{equation}
for some $s>0$ where $c'_j(t,q,q^k) \in \BQ[\rho][q,q^k,t]$ for $j=1,
\dots,s$.
Assuming that $p_e(\rho q^k) \neq 0$ for $k>0$, from this one can 
solve inductively to determine the coefficients of the formal power 
series solution in terms of $a_0(q)$. Solving 
inductively~\eqref{arec}, 
we obtain that 
\begin{equation}
\label{abounds}
 a_k(q) \prod_{j=1}^k p_e(\rho q^j,q) \in \BC[q], \qquad  | a_k(q) 
 \prod_{j=1}^k p_e(\rho q^j,q) |_1=O(C^k) \,.
\end{equation}
Factorising the edge polynomial into a product of linear forms,
combined with Lemma \ref{lem.dioph} below implies that 
the constructed formal power series solutions are convergent on a 
disk around $t=0$.

When $\rho$ is a root of the edge polynomial with multiplicity $m$, 
one
complements this method using $p_e(\rho q^{k+\varepsilon},q) + 
O(\varepsilon^m)$ to obtain $m$ formal power series solutions. Using 
Lemma~\ref{lem.eps} below,
we again conclude that the power series coefficients $a_k(q)$ satisfy 
the bounds
~\eqref{abounds}.

Finally, in case $p_e(\rho q^k,q)=0$ for some $k>0$, one replaces 
$k$ by 
$k+\eps$, uses the Lemma~\ref{lem.eps} and concludes the same 
bounds ~\eqref{abounds}. This concludes the proof of the theorem.
\end{proof}

We now give an $\varepsilon$-deformation of the quantum factorial.
Recall the $n$th polylogarithm $\Li_n(q)=\sum_{s=1}^\infty s^{-n} 
q^s$, 
which satisfies the property $\Li_{-\ell}(q)(1-q)^{\ell+1} \in 
\BZ[q]$ 
for $\ell >0$. For $k, \ell >0$, consider the $q$-harmonic sums 
$H_{k,\ell}(q)$
defined by
\begin{equation}
\label{Hlk}
H_{\ell,k}(q) = \sum_{s=1}^\infty s^{\ell-1} q^{js} 
= \sum_{j=1}^k \Li_{1-\ell}(q^j) \,.
\end{equation}
For example,
\begin{equation}
\qquad H_{1,k}(q) =  \sum_{j=1}^k \frac{q^j}{1-q^j} , 
\qquad  H_{2,k}(q) =  \sum_{j=1}^k \frac{q^j}{(1-q^j)^2}
\qquad H_{3,k}(q) = \sum_{j=1}^k \frac{q^j(1+q^j)}{(1-q^j)^3} \,.
\end{equation}

\begin{lemma}
\label{lem.eps}
For every $k \geq 0$, we have:
\begin{equation}
\label{qqk2}
\frac{(q;q)_k}{(q^{1+\ve};q)_k} =\exp \Big( 
\sum_{\ell=1}^\infty \frac{(\varepsilon \log q)^\ell}{\ell!} 
H_{\ell,k}(q) 
\Big) =
 \frac{(q;q)_k}{(q^{1+\ve};q)_k} = \sum_{\ell=0}^\infty 
\frac{P_{\ell,k}(q)}{(q;q)_k^{\ell}}(\log q)^\ell \ve^\ell
\end{equation}
for unique $P_{k,\ell}(q) \in \BQ[q]$.
Moreover, for every $\ell \geq 1$, there exist $C_\ell>0$ such that 
$|H_{\ell,k}(q)|_1 \leq C_\ell  2^{k \ell}$ for all $k \geq 0$.
\end{lemma}

For example, we have $P_{0,k}(q)=1$ and 
\begin{equation}
P_{1,k}(q) = \tilde H_{1,k}(q), \quad  
P_{2,k}(q) = \frac{1}{2} \tilde H_{1,k}(q)^2 
+ \frac{1}{2} \tilde H_{2,k}(q), \quad
P_{3,k}(q) = \frac{1}{6} \tilde H_{1,k}(q)^3 
+ \frac{1}{2} \tilde H_{1,k}(q) \tilde H_{2,k}(q)
+ \frac{1}{2} \tilde H_{3,k}(q)
\end{equation}
where $\tilde H_{\ell,k}(q)=(q;q)^\ell H_{\ell,k}(q)$.
For related formulas, see also~\cite[Sec.2.2]{GK17} 
and~\cite[Lem.A.4]{GZ24}.

\begin{proof}
We compute:
\begin{equation}
\begin{aligned}
-\log(q^{1+\varepsilon};q)_k + \log(q;q)_k  &= 
-\sum_{j=1}^k \log(1-q^{j+
\varepsilon})
+ \sum_{j=1}^k \log(1-q^{j}) = \sum_{j=1}^k \sum_{s=1}^\infty 
\frac{q^js (q^{\varepsilon s}-1)}{s} \\ &= 
\sum_{s=1}^\infty \frac{q^s(1-q^{ks})}
{1-q^s} \frac{\log(q^{\varepsilon s}-1)}{s} = \sum_{s=1}^\infty 
\frac{q^s(1-q^{ks})}{1-q^s} 
\sum_{\ell=1}^\infty \frac{(\varepsilon s \log 
q)^\ell}{s \, \ell!} \\
& =\sum_{\ell=1}^\infty \frac{(\varepsilon \log q)^\ell}{s \,\ell!} 
\sum_{s=1}^\infty s^\ell(q^s + \dots + q^{ks}) 
= \sum_{\ell=1}^\infty \frac{(\varepsilon \log q)^\ell}{\ell!} 
\sum_{j=1}^k 
\sum_{s=1}^\infty s^{\ell-1} q^{js} \\
&= \sum_{\ell=1}^\infty \frac{(\varepsilon \log q)^\ell}{\ell!} 
H_{\ell,k}(q) \,.
\end{aligned}    
\end{equation}
This proves~\eqref{qqk2}, which in turn implies the first equality 
in~\eqref{qqk2}.
This and the fact that 
$|(q;q)_n|_1 \leq \prod_{j=1}^n |1-q^j|_1=2^n$ implies 
the existence of the bound on their $\ell_1$-norm. 
\end{proof}

The next lemma is inspired by a theorem of 
Hardy--Littlewood~\cite{Hardy} who proved that
\begin{equation}
\liminf_n |(q;q)_n|^{1/n} = \liminf_n |1-q^n|^{1/n} 
\end{equation}
for $q$ in the complex unit circle $\mathbb{T}$, combined with
a result of Lubinsky~\cite{Lub99}
who proved that
\begin{equation}
\liminf_n |1-q^n|^{1/n}=1
\end{equation}
for a set of $q$ in the unit circle of Hausdorff dimension 0 and 
logarithmic dimension 2~\cite{Lub99}, in particular, almost 
everywhere on the unit circle. 

\begin{lemma}
\label{lem.dioph}
Fix an algebraic function $\alpha(q)$ of $q$. For almost all $q \in 
\mathbb{T}$, we have:
\begin{equation}
\label{liminf}
  \liminf_{n\to\infty} \left|(\alpha(q),q)\right|^{1/n} =
  \liminf_{n\to\infty} \left|\alpha(q)-q^n\right|^{1/n} =
  1 \,.
\end{equation}
\end{lemma}

Doron Lubinsky informs us that the above equation~\eqref{liminf}
actually holds for a set of $q$ in the unit circle of Hausdorff 
dimension 0 and logarithmic dimension 2~\cite{Lubinsky:personal}.

\begin{proof}
Let us prove first that
$\liminf_{n\to\infty} \left|\alpha(q)-q^n\right|^{1/n}$ almost 
everywhere on $\mathbb{T}$.
The upper bound is trivial since $\left|\alpha(q)-q^n\right| \leq M$
for all $q \in \mathbb{T}$.
For the lower bound, fix $q = \ee(x):= e^{2\pi i \theta}$,
If $\left|\alpha(q)\right| \neq 1$, then  
$2 \geq \left|q^n - \alpha(q)\right| \geq |1-|\alpha||>0$, so 
its $n$-th root
convergences to 1. 

If $\left|\alpha(q)\right| = 1$ (this happens on a finite union
$J$ of closed intervals in $S^1$), we use the Borel-Cantelli 
lemma, see e.g.,~\cite{Tao:measure}.
When $q \in J$, we write $\alpha(q) = e^{2\pi i \phi(\theta)}$
and denote $\|\cdot\|$ the distance of a real number to the 
nearest integer. 
Then,
\begin{equation}
  |\ee(n\theta)-\ee(\phi(\theta))| = 2 \sin \pi (n\theta-\phi(\theta))
  = 2 \sin \pi \|n\theta-\phi(\theta))\| = 2 \pi \|n\theta-
  \phi(\theta))\|
  + O(\|n\theta-\phi(\theta))\|^3) \,.
\end{equation}
If $\liminf < 1$, there exists some $\eps > 0$ such that the 
inequality
\[
\|n\theta - \phi(\theta)\| < \frac{1}{(1+\eps)^n}
\]
holds for infinitely many values of $n$. We define $E_n$ as the set of
$\theta \in [0,1)$ satisfying this inequality for a fixed $n$. 
Since $\phi(\theta)$ is a fixed piecewise smooth curve independent of 
$n$,
the function $g(\theta) = n\theta - \phi_j(\theta)$ winds around the 
circle $O(n)$
times. Thus, $E_n$ is a union of $O(n)$ intervals, each having length
$O\left(\frac{1}{n(1+\eps)^n}\right)$. The total Lebesgue 
measure $\mu$ of these
intervals satisfies:
\[
  \mu(E_n) = O(n) \times O\left(\frac{1}{n(1+\eps)^n}\right)
  = O\left(\frac{1}{(1+\eps)^n}\right) \,.
\]
Because the geometric series $\sum_{n=1}^\infty \mu(E_n)$ converges,
the first Borel-Cantelli Lemma implies that the set of $\theta$ 
belonging to
infinitely many $E_n$ has a measure $0$.

This proves the second equality in~\eqref{liminf}. For the first
equality, we use the identity
\begin{equation}
\label{Phi11}
\sum_{n=0}^\infty \frac{t^n}{(q \alpha(q);q)_n}
= (1-\alpha(q)) 
\left(\sum_{n=0}^\infty \frac{t^n}{(q;q)_n} \right)
\left(\sum_{n=0}^\infty q^{n(n+1)/2} \frac{(-t)^n}{(q;q)_n(1-
\alpha(q) q^n)} \right)
\end{equation}
of ~\cite[Lem.1.3]{Div09}. Together with
\begin{equation}
\liminf_n |\alpha(q)-q^n|^{1/n} = \liminf_n |(q;q)_n|^{1/n} =1
\end{equation}
almost everywhere, it implies that 
$1 \leq \liminf_n |(q;q)_n(\alpha(q)-q^n)|^{1/n} \leq 1$ almost
everywhere, hence the radius of convergence of the right hand 
side of~\eqref{Phi11} is $1$, equal to the radius of convergence of
the left hand side, concluding that 
$\liminf_n|(q \alpha(q);q)_n|^{1/n}
=1$ almost everywhere on $\mathbb{T}$.
\end{proof}

\begin{figure}
\centering
\begin{tikzpicture}
\draw[xshift=-5cm](0,0)--(3,0)--(-1,1)--(0,0);
\draw(0,0)--(3,0)--(2,1)--(0,0);
\draw[xshift=5cm](0,0)--(3,0)--(4,1)--(0,0);
\draw[xshift=-5cm](0,0) node[below] {$(0,0)$};
\draw[xshift=-5cm](3,0) node[below] {$(B,0)$};
\draw[xshift=-5cm](-1,1) node[above] {$(A,1)$};
\draw(0,0) node[below] {$(0,0)$};
\draw(3,0) node[below] {$(B,0)$};
\draw(2,1) node[above] {$(A,1)$};
\draw[xshift=5cm](0,0) node[below] {$(0,0)$};
\draw[xshift=5cm](3,0) node[below] {$(B,0)$};
\draw[xshift=5cm](4,1) node[above] {$(A,1)$};
\end{tikzpicture}
\caption{Shapes of the Newton polygon for 
$(1-\sigma)^B - (-1)^A t \sigma^A$.}
\label{fig:NPforAB}
\end{figure}
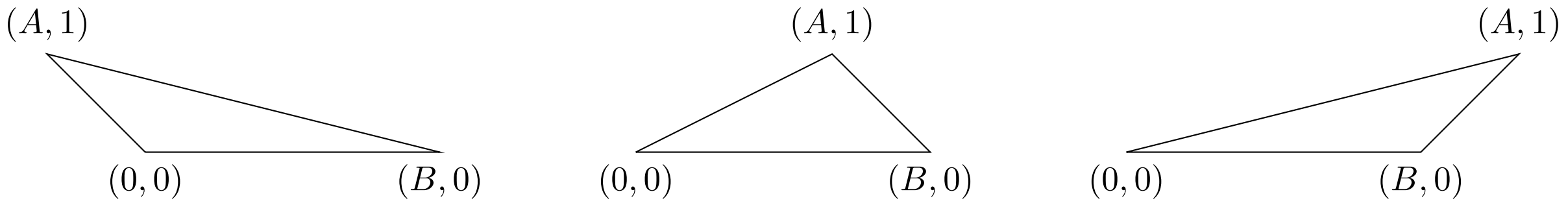

We now illustrate our theorem for the linear $q$-difference 
equation 

\begin{equation}
\label{nahmAB}
((1-\sigma)^B - (-1)^A t \sigma^A ) f(t,q)=0, \qquad (\sigma f)
(t,q)=f(qt,q) \,.
\end{equation}
associated to $\Omega_{A,B}$. When $B>0$, its Newton polygon is a 
triangle 
with vertices $(0,0)$, $(0,B)$ and $(1,A)$. In the three regions 
where $B>0$
is shown in Figure~\ref{fig:NPforAB}.

It is easy to see that it is quasi-unipotent.

When $B>A>0$, a fundamental solution to the $q$-difference 
equation~\eqref{nahmAB} is
\begin{equation}
t^{\ve} \sum_{k=0}^\infty 
\frac{q^{A(k(k-1)/2+A k \ve}}{(q^{1+\ve};q)_k^B} + O(\ve^B) \,.
\end{equation}

This can be expanded into a power series in $\ve+O(\ve^B)$ using 
Lemma~\ref{lem.eps} and the coefficients of the powers of 
$\ve$, arranged 
into a column vector of size $B$ gives a fundamental solution 
to the equation. 

On the other hand, when $B > A > 0$, to find the solution 
corresponding to the 
first edge of the upper Newton polygon of the $q$-difference 
equation~\eqref{nahmAB}, we need to shear it so that it 
becomes horizontal. This
is achieved by the change of variables
$f(t,q)= t^{\alpha \tfrac{\log t}{\log q}} g(t,q)$. Substituting 
into~\eqref{nahmAB} and using the fact that 
$\sigma( t^{\alpha \tfrac{\log t}{\log q}})=q^\alpha t^{2 \alpha}$,
we see that
with $\alpha=-1/(2A)$, $g$ satisfies the equation
\begin{equation}
\label{geqn}
  (  (1-q t^{-\tfrac{1}{A}}\sigma)^B -(-1)A q^A \sigma^A ) g(t,q) = 0 
\end{equation}
which whose upper Newton polygon has first segment 
$[0,A] \times \{0\}$.
The edge polynomial of the horizontal edge is 
$p_e(\sigma,q)=1-(-1)^A q^A \sigma^A$
whose roots are $\rho=(-1)^A \zeta_A q^{-1}$ where $\zeta^A=1$.

\section{State integrals}
\label{sec.state}

In this section, we give the detailed definition of the
quantum modular cocycle $\Omega_{A,B}$, prove its cocycle property,
and express its trace in terms of the finite quantum dilogarithm.

\subsection{Cocycles via the modular quantum dilogarithm}

It is explained in~\cite{GKZ} that entries of the matrix valued 
cocycles can be 
expressed in terms of integrals of Faddeev's modular quantum 
dilogarithm.
For
$D_{\gamma,m,n}(z;\tau)_{ij}=\ee(j(z+m\tau)+(i-1)(\frac{z}
{c\tau+d}+n\frac{a\tau+b}{c\tau+d}))$,
we define
\[
\Omega_{A,B,\gamma}(\tau)
=\!\!\!\!
\sum_{m\in\ZZ/c\ZZ}
\int_{\calC_{m+1,0}}\!\!\!\!\!\!\!\!
\Phi_{\gamma,m+1,0}(z;\tau)^B\;
\ee(\tfrac{Acz^2}{2(c\tau+d)}+A(m-\tfrac{1}{2(c\tau+d)}
-\tfrac{1}{2})z+A(\tfrac{m^2\tau}{2}-\tfrac{m\tau}{2}-\tfrac{m}
{2}))\;D_{\gamma,m,0}(z;\tau)\;P\;dz.
\]
where $\calC_{m,n}=\frac{n}{c}-m\frac{c\tau+d}{c}+i\RR-0$.
As defined, it is not obviously a quantum modular cocycle.
We will prove this without proving $q$-series or root of unity 
identities (the quadratic identities of~\cite{GZ23}).
The proof will use the extension of the cocycle to the positive reals.
Then the failure of the cocycle will satisfy a simple $q$-difference 
equation.
We will show that this equation has only one solution---satisfying 
certain 
growth conditions---that is given by a multiple of the identity 
matrix.
We will then compute a limiting case to identify this constant.
As a consequence, this will prove the quadratic relations needed for 
a 
proof 
using factorisation, which in turn prove that we have associated 
quantum modular forms.

Firstly, we will prove an elementary duality satisfied by the 
companion matrix 
of the linear $q$-difference equation~\eqref{nahmAB}. 
Let $C=C_{A,B}$ denote the companion matrix of~\eqref{nahmAB} 
and $P=P_{A,B}$ 
denote the gauge transformation matrix from Section~\ref{sub.13}. 
These are 
matrices of size $A$, $B$ and $B-A$ respectively when $A>B>0$, 
$B>A>0$ and $B>0>A$, respectively.

\begin{lemma}
\label{lem.CP}
For all $B>0$, we have: $C(t)^TP=PC(t)$.
\end{lemma}

\begin{proof}
The proof follows from showing that both sides are equal to 
the following matrices:
\begin{itemize}
\item for $A<0<B$ we have the block matrix with decomposition 
$B-A=1+(B-A-1)$
\[
\left(\!\!\begin{array}{c|c}
   (-1)^{A+B}t & 0 \\\hline
   0 & \begin{array}{ccccc}
    \cdots & -\binom{B}{3} & \binom{B}{2} & -\binom{B}{1} & 1
    \\
    \cdots & \binom{B}{2} & -\binom{B}{1} & 1 & 0\\
    \cdots & -\binom{B}{1} & 1 & 0 & 0\\
    \cdots & 1 & 0 & 0 & 0\\
    \vdots & \vdots & \vdots & \vdots & \vdots\\
    1 & 0 & 0 & \cdots & 0
    \end{array}
\end{array}\!\!\right)\!,
\]
\item 
for $0<A<B$ we have the block matrix with decomposition 
$B=(A+1)+(B-A-1)$
\[
\left(\!\!\begin{array}{c|c}
\begin{array}{ccccc}
0 & \cdots & 0 & 0 & (-1)^B\\
0 & \cdots & 0 & (-1)^B & (-1)^{B+1}\binom{B}{1}\!\!\!\!\\
0 & \cdots & (-1)^B & (-1)^{B+1}\binom{B}{1} & 
(-1)^{B}\binom{B}{2}\!\!\!\!\\
\vdots & \vdots & \vdots & \vdots & \vdots\\
(-1)^B & \cdots & (-1)^{B+A}\binom{B}{A-2} & 
(-1)^{B+A+1}\binom{B}{A-1} & \overset{(-1)^{B+A}\binom{B}
{A}}{+(-1)^At}\!\!\!\!
\end{array} & 0
\\ \hline
0 & \begin{array}{ccccc}
\!\!\!\!\cdots & -\binom{B}{3} & \binom{B}{2} & -\binom{B}
{1} & 1\\
\!\!\!\!\cdots & \binom{B}{2} & -\binom{B}{1} & 1 & 0\\
\!\!\!\!\cdots & -\binom{B}{1} & 1 & 0 & 0\\
\!\!\!\!\cdots & 1 & 0 & 0 & 0\\
\!\!\!\!\vdots & \vdots & \vdots & \vdots & \vdots\\
\!\!\!\!1 & 0 & 0 & \cdots & 0
\end{array}
\end{array}\!\!\right)\!,
\]
\item for $0<B<A$ it is a sum of two matrices
\[
(-1)^At^{-1}\left(\begin{array}{ccccc}
\binom{B}{0}\binom{B}{0} & -\binom{B}{0}\binom{B}{1} & 
\binom{B}{0}\binom{B}{2} & \cdots & 0\\
-\binom{B}{1}\binom{B}{0} & \binom{B}{1}\binom{B}{1} & -
\binom{B}{1}\binom{B}{2} & \cdots & 0\\
\binom{B}{2}\binom{B}{0} & -\binom{B}{2}\binom{B}{1} & 
\binom{B}{2}\binom{B}{2} & \cdots & 0\\
\vdots & \vdots & \vdots & \vdots & \vdots\\
0 & 0 & 0 & \cdots & 0\\
\end{array}\right)
+
\left(\!\!\begin{array}{ccccc}
0 & \cdots & 0 & 0 & 0\\
0 & \cdots & 0 & 0 & \binom{B}{0}\\
0 & \cdots & 0 & \binom{B}{0} & -\binom{B}{1}\\
\vdots & \vdots & \vdots & \vdots & \vdots\\
0 & \binom{B}{0} & -\binom{B}{1} & \cdots & 0\\
\end{array}\!\!\right)\!.
\]
\end{itemize}
\end{proof}

Now we will deform the cocycle by adding a Jacobi variable $w$ 
and 
shifts so 
that $D_{\gamma,m,n}^{(u)}(z;\tau)_{ij}=
\ee((j+u_1)(z+m\tau)+(i-u_2-1)(\frac{z}
{c\tau+d}+n\frac{a\tau+b}{c\tau+d}))$
and
\[
\begin{aligned}
\Omega_{A,B,\gamma}^{(u)}(w;\tau)&=\sum_{m\in\ZZ/c\ZZ}
\int_{\calC_{m+1,0}}\!\!
\Phi_{\gamma,m+1,0}(z;\tau)^B
\ee\big(\tfrac{Acz^2}{2(c\tau+d)}+A(m-\tfrac{1}{2(c\tau+d)}
-\tfrac{1}{2})z\\
&\qquad\qquad\qquad\qquad+A(\tfrac{m^2\tau}{2}-\tfrac{m\tau}{2}-
\tfrac{m}{2})+
\tfrac{(cz+(c\tau+d)m)w}{c\tau+d}\big)\;D_{\gamma,m,0}(z;\tau)\;P\;dz.
\end{aligned}
\]

Denoting $\Omega_{A,B,\gamma}(w;\tau)=
\Omega_{A,B,\gamma}^{(0,0)}(w;\tau)$
we have the following:

\begin{theorem}
\label{thm:cocycle}
The functions $\Omega_{A,B,\gamma}$ satisfy the cocycle 
condition, i.e., we have the following identity:
\[
\Omega_{A,B,\gamma\gamma'}(w;\tau)
=
\Omega_{A,B,\gamma}(\gamma'(w;\tau))\;
\Omega_{A,B,\gamma'}(w;\tau)\,.
\]
\end{theorem}

\noindent We will now prove several lemmas used to prove 
this identity before giving the proof.

Trivially, we find that
\begin{align*}
&\Omega_{A,B,\gamma}^{(u)}(w+a\tau+b;\tau) =
\Omega_{A,B,\gamma}^{(u+(a,1))}(w;\tau), 
&&\Omega_{A,B,\gamma}^{(u)}(w+\tau;\tau) = 
\Omega_{A,B,\gamma}^{(u+(1,d))}(w;\tau),\\
&\Omega_{A,B,\gamma}^{(u)}(w+c\tau+d;\tau) = 
\Omega_{A,B,\gamma}^{(u+(c,0))}(w;\tau), 
&&\Omega_{A,B,\gamma}^{(u)}(w+1;\tau) = 
\Omega_{A,B,\gamma}^{(u+(0,-c))}(w;\tau)\,.
\end{align*}
We also deduce that
\[
\Omega_{A,B,\gamma}^{(u+(1,0))}(w;\tau)
=
\Omega_{A,B,\gamma}^{(u)}(w;\tau)C(q^{u_1}\ee(w))\,,\qquad
\Omega_{A,B,\gamma}^{(u+(1,0))}(w;\tau)
=
C(\tq^{u_2}\ee(\tfrac{w}{c\tau+d}))^{-1}
\Omega_{A,B,\gamma}^{(u)}(w;\tau)\,,
\]
where we use the defining equation for $P$.
Notice that when $\ee(\tfrac{w}{c\tau+d})=\ee(w)$ and $q=\tq$, 
we find that 
for all $m\in\ZZ$
\[
\mathrm{tr}(\Omega_{A,B,\gamma}^{(u+(m,m))}(w;\tau))
=
\mathrm{tr}(\Omega_{A,B,\gamma}^{(u)}(w;\tau))\,.
\]

We will now fix $\gamma,\gamma'\in\SL_2(\ZZ)$ for the remainder 
of the section.
We define
\[
H(w;\tau)
=
P^{-1}(\Omega_{\gamma\gamma'}(w;\tau)^{-1}\Omega_{\gamma}(\gamma'(w;
\tau))\Omega_{\gamma'}(w;\tau))^{T}P\,.
\]

\begin{lemma}
The function $H(w;\tau)$ is meromorphic for $w\in \CC$ and 
$c\tau+d\notin\RR_{\leq0}$ and satisfies $H(w+1;\tau)=H(w;\tau)$ and 
$H(w+\tau;\tau)=C(\ee(w))H(w;\tau)C(\ee(w);\tau)^{-1}$.
\end{lemma}

Therefore, we find that $H(w;\tau)$ has a Fourier expansion and 
letting $t=\ee(w)$ 
we have $H=H(t;\tau)$. Therefore, when we fix $\tau$ to be a 
generic real 
number, we find that $H(t;\tau)$ is a holomorphic function 
of $t\in\CC^\times$,
since any pole would create a dense set of poles 
of a circle. We will now determine the behaviour of $H$ as 
$t^{\pm1}\to 0$.

\begin{lemma}
For generic $\tau\in\RR$, we find that 
$H(t;\tau)\in \GL_{r}(\CC[t^{\pm1}])$.
\end{lemma}

\begin{proof}
From Theorem~\ref{thm.q=1}, we have the existence of a 
fundamental matrix of 
solutions 
$G(t,\tau)$ that consist of Laurent series in $t$ 
times $t^{\alpha\frac{\log(t)}
{\log(q)}+\beta}$ 
for some $\alpha,\beta$. Therefore, the matrix 
$G(t;\tau)^{-1}H(t;\tau)G(t;\tau)$ 
is invariant under $t\mapsto qt$ and therefore constant.
Since $H(t,\tau)$ is holomorphic on $\CC^\times$ (which in 
particular means it has 
no monodromy around $0,\infty$), it means that $H$ must behave 
like a Laurent 
series near $0,\infty$.
Therefore, we have a meromorphic function on $\CC^\times$ with
possible poles 
only at $0,\infty$.
\end{proof}

We will now compute the behaviour of the function $H$ as we 
approach either 
$0,\infty$ in our different cases. In general, we will study
the asymptotics 
on the side which has a single associated slope of the Newton 
polygon.

\begin{lemma}
\label{lem:limiting.behaviour}
For $0<A<B$ the function $H(t,\tau)$ has asymptotics as $t\to 0$ 
given by
\[
H(t;\tau) = \mathrm{I} +\mathrm{O}(t)\,.
\]
\end{lemma}

\begin{proof}
We can compute the asymptotics of $\Omega_\gamma(w;\tau)$ as 
$w\to i\infty$ 
with $w\in i\RR$.
In particular, we see that moving the contours to the right 
will decrease their 
size to a constant times $|\ee(w \alpha)|$ for some $\alpha$.
At some point they will pass a pole of the integrand giving
\[
\Omega_{A,B,\gamma}(w;\tau)
=
\Res{z=0}\;
\Phi_{\gamma,1,0}(z;\tau)^B\;
\ee(\tfrac{Acz^2}{2(c\tau+d)}-A(\tfrac{1}{2(c\tau+d)}
+\tfrac{1}{2})z+\tfrac{czw}{c\tau+d})\;D_{\gamma,0,0}(z;\tau)\;P\;dz
+
O(\ee(w\alpha))\,.
\]
for some $\alpha>0$.
We can apply the cocycle condition to find
\[
\begin{aligned}
&\Res{z=0}\;
\frac{\Phi_{\gamma\gamma',1,1}(z;\tau)^B}{(1-\ee(\frac{z}
{c''\tau+d''}))^B}\;
\ee(\tfrac{Ac''z^2}{2(c''\tau+d'')}-A(\tfrac{1}{2(c''\tau+d'')}+
\tfrac{1}{2})z+
\tfrac{c''zw}{c''\tau+d''})\;D_{\gamma\gamma',0,0}(z;\tau)\;dz\\
&=\Res{z=0}\;
\frac{\Phi_{\gamma,1,1}(\gamma'(z;\tau))^B
\Phi_{\gamma',1,1}(z;\tau)^B}{(1-\ee(\frac{z}
{c''\tau+d''}))^B}\;
\ee(\tfrac{Ac''z^2}{2(c''\tau+d'')}-A(\tfrac{1}{2(c''\tau+d'')}+
\tfrac{1}{2})z+
\tfrac{c''zw}{c''\tau+d''})\;D_{\gamma\gamma',0,0}(z;\tau)\;dz\\
\end{aligned}
\]
This expresses the asymptotics of $\Omega_{A,B,\gamma\gamma'}$ 
as a bilinear 
combination of the asymptotics of 
$\Omega_{A,B,\gamma}(\gamma'(w;\tau))$ and 
$\Omega_{A,B,\gamma'}(w;\tau)$.
An explicit computation of the bilinear form produces the 
required identity.
\end{proof}

\begin{lemma}
\label{lem:limiting.behaviour.2}
For $0<B<A$ or $A<0<B$ the function $H(t,\tau)$ has asymptotics as 
$t\to\infty$ given by
\[
H(t;\tau) = \mathrm{I} +\mathrm{O}(t^{-1})\,.
\]
\end{lemma}
\begin{proof}
For $A>B$, applying Fourier invariance of the modular 
quantum dilogarithm, 
we find that
\[
\begin{aligned}
&\sum_{m\in\ZZ/c\ZZ}
\int_{\calC_{m+1,0}}
\Phi_{\gamma,m+1,0}(z;\tau)^B\;
\ee\big(\tfrac{Acz^2}{2(c\tau+d)}
+A(m-\tfrac{1}{2(c\tau+d)}-\tfrac{1}{2})z\\
&\qquad\qquad\qquad\qquad\qquad\qquad\qquad\qquad+A(\tfrac{m^2\tau}
{2}-\tfrac{m\tau}{2}-\tfrac{m}{2})+\tfrac{(cz+(c\tau+d)m)w}
{c\tau+d}\big)\;D_{\gamma,m,0}(z;
\tau)\;P\;dz\\
&=\mu_\gamma^{-B}(c\tau+d)^{B/2}\ee(B\tfrac{\gamma\tau-\tau}
{24})\sum_{m\in(\ZZ/
c\ZZ)^{B+1}}
\int_{\calC_{m+1,0}^{B+1}}
\prod_{i=1}^{B}\Phi_{\gamma,m_i+1,0}(z_i;\tau)\;
\ee\Big(\tfrac{Acz_0^2}{2(c\tau+d)}+A(m_0-\tfrac{1}{2(c\tau+d)}-
\tfrac{1}
{2})z_0\\
&\qquad+A(\tfrac{m_0^2\tau}{2}-\tfrac{m_0\tau}{2}-\tfrac{m_0}{2})
+\sum_{i=1}^B\tfrac{c(z_i+m_i(\tau+c^{-1}d))(z_i+(m_i-1)(\tau+c^{-1}
d))-
z_i+2(w+z_i+(\tau+c^{-1}d)m_i)(cz_0+(c\tau+d)m_0)}{2(c\tau+d)}\\
&\qquad\qquad\qquad\qquad\qquad\qquad\qquad\qquad+\ell 
(z_i+m_i\tau)\Big)D_{\gamma,m_0,0}(z_0;\tau)Pdz
\end{aligned}
\]
The integral in $z_0$ is decoupled from the quantum dilogarithms 
and we can compute 
it explicitly to produce Gaussian behaviour of the form 
$\ee(-Acw^2/2(c\tau+d))$.
Then shifting the remaining contour to the right produces 
a contribution from a 
residue and then an error of relative size $|\ee(w\alpha)|$.
An explicit computation of a similar form outlined in the previous 
lemma~\ref{lem:limiting.behaviour} produces the required 
asymptotic behaviour.
For $A<0<B$ a similar argument follows after applying the 
inversion relation before 
Fourier invariance, which produces Gaussian behaviour of the 
form $\ee((-A+B)cw^2/2(c\tau+d))$.
\end{proof}

\begin{lemma}
\label{lem:dans}
Let $A(t)$ be an $n\times n$ matrix of the following form
\[A(t)=
\begin{pmatrix}
0 & 1 & 0& \dots &0\\
0 & 0 &1 & \dots& 0\\
\vdots & \vdots& \vdots & \ddots & \vdots\\
0 & 0 & 0 &  \cdots & 1\\
-a_0(t) & -a_1(t) & -a_2(t) & \cdots & -a_{n-1}(t)
\end{pmatrix}\,,
\]
where either: (1) we assume $a_i(t)$, $i\ne k$ is a non-zero 
constant $\alpha_i$ and 
$a_k(t)=\alpha_k-\beta t$, where $k$ is a fixed number, $0\le k<n$, 
and $\beta\ne0$; or (2) we assume $a_i(t)=\beta_it^{-1}$ and 
$\beta_0\neq0$. Let $F(t)$ be an $n\times n$ matrix with entries in 
$\CC[t]$. 
Assume that for some $q$ with $|q|=1$ that is not a root of unity 
we have
\[
F(qt)=A(t)F(t)A(t)^{-1}\,.
\]
Then $F(t)=\lambda I_n$ for some constant $\lambda\in\CC$.
\end{lemma}

\begin{proof}
Suppose that $A$ is of the form (1).
Consider the $\CC$-linear mapping $X\colon \CC[t]^n\to \CC[t]^n$ 
given by 
$g(t)\mapsto g(qt)A(t)$ (where $g(t)=(g_0(t),\dots,g_{n-1}(t))$ 
is a row vector of 
polynomials in~$t$). Also consider $Y\colon \CC[t]^n\to \CC[t]^n$ 
given by $g(t)\mapsto 
g(t)F(t)$. The condition of the lemma implies that $X$ and $Y$ 
commute, but also $Y$ 
commutes with multiplication by $t$, whereas $Xt=qtX$. 
Let $e_0,\dots,e_{n-1}$ denote 
the standard basis of $\CC[t]^n$; then $e_i=X^ie_0$ for 
$i=0,\dots,n-1$. We also have 
$X^ne_0=Xe_{n-1}=e_{n-1}A(t)
=(\beta tX^k-\sum_{j=0}^{n-1}\alpha_jX^j)e_0$.

If $k>0$, then $X$ is invertible (as is easily seen since 
$A(t)^{-1}$ has polynomial 
entries) and the above computation of $X^ne_0$ implies that 
$\beta^{-1}q^kX^{-k}
P(X)e_0=te_0$, where $P(x)=x^n+\sum_{j=0}^{n-1} \alpha_jx^j$. 
By induction 
$t^me_0=P_m(X)e_0$ for some Laurent polynomial $P_m$, and then 
also $t^me_j=q^{-mj}
X^jP_{m}(X)e_0$. Thus $\CC[t]^n$ is a cyclic module over 
$\CC[X,X^{-1}]$, and since 
$\CC[t]^n$ is of infinite dimension as a $\CC$-vector space, 
it is also free (as a $
\CC[X,X^{-1}]$-module). Since $Y$ commutes with $X$, we must 
have $Y=Q(X)$ for some 
Laurent polynomial~$Q$ (that depends only on $Ye_0$). On the 
other hand, $Q(X)t = 
tQ(qX)$, and since~$Y$ commutes with $t$, we must have 
$Q(qz)=Q(z)$, which implies 
that $Q$ is constant.
	
If $k=0$, the same argument works with $\CC[X]$ in place of 
$\CC[X,X^{-1}]$. If $A$ is of the form (2), then simply replace 
$q\to q^{-1}$ take $A(t)^{-1}$ conjugated by the matrix with $1$ 
on the off-diagonal and $0$ elsewhere to reduce to the previous 
case.
\end{proof}

\begin{proof}[Proof of Theorem~\ref{thm:cocycle}]
Lemma~\ref{lem:dans} imply that $H(t;\tau)$ is a constant 
multiple of the identity.
Lemmas~\ref{lem:limiting.behaviour}~and
~\ref{lem:limiting.behaviour.2} then imply 
that this must simply be the identity.
This holds for generic and large enough $\tau$ on the real line.
Since $H$ is holomorphic on a cut plane, we conclude that this 
identity must 
hold there.
\end{proof}

\subsection{\texorpdfstring{Proof of 
Theorem~\ref{thm.nahm}}{Proof of Theorem 1.2}}

We will now use the same contour integration tricks in~\cite{RW} 
to prove the main theorem.

\begin{proof}[Proof of Theorem~\ref{thm.nahm}]
Recall that $\tau^*$ is a fixed point of $\gamma$ and $N=\mathrm{tr}(\gamma)-2$. 
Firstly, notice that for $z_{m,w}=\frac{m\tau^*+w}{\eps^{-1}-1}$ 
we find that
\[
\begin{aligned}
&\ee\Big(c\frac{z_{m,w}^2}{2\eps}+mz_{m,w}+\frac{m^2\tau^*}{2}
+\frac{z_{m,w}}{2\eps}-\frac{z_{m,w}}{2}-\frac{m\tau^*}{2}
-\frac{m}{2}\Big)\\
&=
\ee_{2N}(w^2c - bm^2 + (a-d)wm+Nwm+Nw+Nm)
=
\ang{(m,w)}_\gamma\,.
\end{aligned}
\]
Then one can compute that
\[
\mathrm{tr}((x^{i-1}y^j)_{ij}P)
=
\frac{(1-x)^By^A-(1-y)^Bx^A}{1-x/y}\,.
\]
Combining all of this we find that
\[
\begin{aligned}
&\frac{\eps}{1-\eps}\sum_{x_1=0}^{c-1}\sum_{x_2=0}^{N-1}
\ang{x}_\gamma^AF_{\gamma}^{+}(x)^B\\
&=\sum_{m=0}^{c-1}
\bigg(\int_{\calC_{m,0}}\!\!\!\!\!\!-\int_{\calC_{m,0}+
\eps-1}\bigg)\frac{\ee\big(\tfrac{Acz^2}{2\eps}+A(m+\tfrac{1}{2\eps}-
\tfrac{1}{2})z+A(\tfrac{m^2\tau^*}{2}-\tfrac{m\tau^*}{2}-\tfrac{m}
{2})+
\tfrac{(cz+\eps m)w}{\eps}\big)}{1-q^{-m}\ee(z(1-\eps)/
\eps)}\Phi_{\gamma,m,0}(z;\tau^*)^Bdz\\
&=\sum_{m,j=0}^{c-1}
\int_{\calC_{m,0}+j\frac{\eps-1}{c}}\frac{\ee\big(\tfrac{Acz^2}
{2\eps}+A(m+\tfrac{1}{2\eps}-\tfrac{1}{2})z+A(\tfrac{m^2\tau^*}{2}-
\tfrac{m\tau^*}{2}-\tfrac{m}{2})+\tfrac{(cz+\eps m)w}{\eps}\big)}{1-
q^{-
m}\ee(z(1-\eps)/\eps)}\Phi_{\gamma,m,0}(z;\tau^*)^Bdz\\
&\qquad\qquad-\sum_{m,j=0}^{c-1}
\int_{a\tau^*+b+\calC_{m,0}+j\frac{\eps-1}
{c}}\frac{\ee\big(\tfrac{Acz^2}{2\eps}+A(m+1-a+\tfrac{1}{2\eps}-
\tfrac{1}{2})z\big)}{1-q^{-m-1+a}\ee(z(1-\eps)/\eps)}\Phi_{\gamma,m+1-
a,0}(z;\tau^*)^B\\
&\qquad\qquad\qquad\qquad\qquad\qquad\qquad\qquad\times
\ee\big(A(\tfrac{(m+1-a)^2\tau^*}{2}-\tfrac{(m+1-a)\tau^*}{2}-
\tfrac{m+1-a}{2})+\tfrac{(cz+\eps(m+1-a))w}{\eps}\big)dz\\
&=\sum_{m,j=0}^{c-1}
\int_{\calC_{m,0}+j\frac{\eps-1}{c}}\frac{\ee\big(\tfrac{Acz^2}
{2\eps}+A(m+\tfrac{1}{2\eps}-\tfrac{1}{2})z+A(\tfrac{m^2\tau^*}{2}-
\tfrac{m\tau^*}{2}-\tfrac{m}{2})+\tfrac{(cz+\eps m)w}{\eps}\big)}{1-
q^{-
m}\ee(z(1-\eps)/\eps)}\Phi_{\gamma,m,0}(z;\tau^*)^B\\
&\qquad\qquad\qquad\qquad\qquad\qquad\qquad\qquad\times\Big(1-
\frac{(1-
\ee(z/\eps))^B}{(1-q^m\ee(z))^B}\ee(Az(1-\eps)/\eps\Big)\,dz\\
&=\sum_{m,j=0}^{c-1}
\int_{\calC_{m,0}+j\frac{\eps-1}{c}}\ee\big(\tfrac{Acz^2}{2\eps}+A(m-
\tfrac{1}{2\eps}-\tfrac{1}{2})z+A(\tfrac{m^2\tau^*}{2}-\tfrac{m\tau^*}
{2}-\tfrac{m}{2})+\tfrac{(cz+\eps m)w}{\eps}\big)\Phi_{\gamma,m+1,0}
(z;\tau^*)^B\\
&\qquad\qquad\qquad\qquad\qquad\qquad\qquad\qquad\times\frac{(1-
q^m\ee(z))^B\ee(Az/\eps)-(1-\ee(z/\eps))^Bq^{Am}\ee(Az)}{1-q^{-
m}\ee(z(1-\eps)/\eps)}\,dz\\
&=
c\sqrt{\eps}\;\mathrm{tr}\;\Omega_{A,B,\gamma}(\tau^*)\,.
\end{aligned}
\]
\end{proof}

\begin{remark}\label{rem:recip}
The cocycle $\Omega_{A,B}$ satisfies the extra symmetry that
$(\Omega_{A,B,\gamma})$ is conjugate to $ (\Omega_{A,B,\gamma^{-1}})^t$, which is a consequence of 
Lemma~\ref{lem.CP}.  It follows that 
$\Lambda_{\Omega_{A,B,\gamma}} = \Lambda_{\Omega_{A,B,\gamma}}^{-1}$.
\end{remark}

We end this section with a useful lemma. 

\begin{lemma}
\label{lem.vp}
Fix a positive integer $r$, distinct nonzero algebraic numbers 
$\lambda_i$
and nonzero algebraic numbers $c_i$ for $i=1,\dots,r$.
If $\sum_{i=1}^r c_i \lambda_i^n$ is an algebraic integer for 
all sufficiently
large positive integers $n$, then $\lambda_i$ is an algebraic 
integer for 
all $i=1,\dots,r$.
\end{lemma}

\begin{proof}
The generating series 
$G(z)=\sum_{n=0}^\infty (\sum_{i=1}^r c_i \lambda_i^n) z^n
=\sum_{i=1}^r \frac{c_i}{1-\lambda_i z}$  
is a rational function. For every prime $p$, let $|x|_p=p^{-v_p(x)}$ 
denote the $p$-adic norm. By assumption, the coefficients of $G(z)$ 
have
$p$-norm less than or equal to 1, hence the series is convergent for 
$|z|_p < 1$. Thus, $G(z)$ has no poles in the open unit disk 
$|z|_p<1$.
But $G(z)$ is a rational function with poles at $\lambda_i^{-1}$ for 
$i=1,\dots,r$. Thus, $|\lambda_i|_p \leq 1$ for all $i=1,\dots,r$ 
and all primes $p$. The result follows.
\end{proof}

Since 
\begin{equation}
\label{eq.Ng}
N_{\gamma^n} = \varepsilon^{n}+ \varepsilon^{-n}-2 \in \BZ_{\geq 0}
\end{equation}
for all nonzero integers $n$, we obtain.

\begin{corollary}
\label{cor.Ng}
For every hyperbolic $\gamma$ and every prime $p$ we have:
\begin{equation}
   \limsup_n  |N_{\gamma^n}^{1/n}|_p = 1 \,. 
\end{equation}
\end{corollary}

\begin{corollary}
\label{cor.AB}
Theorem~\ref{thm.nahm} implies the first part of Conjecture~\ref{conj.tralg} 
for $\Omega_{A,B}$.
\end{corollary}

\begin{proof}
Theorem~\ref{thm.nahm} and Equation~\eqref{omg} implies 
Equation~\eqref{tromAB} 
for a finite set $\Lambda_{\Omega_{A,B},\gamma}$ of nonzero 
algebraic
numbers. Now the summand 
$\langle k\rangle_{\gamma^n}^A F_{\gamma^n}(k)^B$
of the right hand side of~\eqref{tromAB} is an algebraic integer, 
in fact a
unit, as follows from~\cite{RW:stark}. This, together 
with Corollary~\ref{cor.Ng}
and Lemma~\ref{lem.vp} for $n >0$ conclude that the elements of 
$\Lambda_{\Omega_{A,B},\gamma}$ are integral, and for $n<0$, 
conclude that
their inverses are integral, hence they are units.
\end{proof}

\section{\texorpdfstring{Identities for finite quantum dilogarithms for 
$\gamma$ and $\gamma^2$}{Identities for finite quantum dilogarithms for 
gamma and gamma2}}
\label{sec.second}

In this section, we recall the defining equations of finite quantum 
dilogarithms, with a focus on those that are derived from a 
hyperbolic element $\gamma$ of $\SL_2(\BZ)$.
We then show that the finite quantum dilogarithms
constructed from $\gamma$ and $\gamma^2$ satisfy a system of
polynomial equations introduced by Izumi, thus giving a 
full description of the Drinfeld centre $\calZ(\calC_\gamma)$
of the fusion category $\calC_\gamma$ and its $S$-matrix.

\subsection{Finite quantum dilogarithms}
\label{sub.finiteQDL}

Let $(G,\ang{\cdot})$ be a finite metric group, that is,
a finite abelian group $G$ together with a choice of a function 
$\ang{\cdot} : G \to \CC^\times$ that satisfies $\ang{x}=\ang{-x}$
for all $x$ and such that
\begin{equation}
\ang{x;y} := \frac{\ang{x+y}}{\ang{x} \ang{y}}
\end{equation}
is a non-degenerate bicharacter on~$G$.
Following~\cite{RW}, we say that $E\colon G\to \CC$ is a \emph{finite quantum dilogarithm} if it satisfies
\[E(x)E(y)\ang{x;y} = -\widehat{E}(0) + \frac{1}{\sqrt{|G|}}\sum_{t\in G}E(t)E(x-t)E(y-t)\ang{t}\]
for all $x,y\in G$, together with the normalization $E(0)^2=\sqrt{|G|}E(0)+1$. This is a version of the more familiar non-commutative pentagon relation for quantum dilogarithm. Any finite quantum dilogarithm satisfies the reflection relation
\begin{equation} \label{eq:reflection}
E(x)E(-x)\ang{x} = 1+\sqrt{|G|}E(0)\delta(x)\,,
\end{equation}
and an identity for the Fourier transform 
\begin{equation} 
\label{eq:fourier}
\widehat{E}(x) = \mu \, E(x)\ang{x}\,,
\qquad
\mu = \widehat{E}(0)/E(0)\,,\qquad \mu^3 = \frac{1}{\sqrt{|G|}}\sum_{x\in G}\ang{x}^{-1}\,,
\end{equation}
see~\cite[Theorem 5]{RW}. From these it follows that the Fourier transform of $E^{-1}$ satisfies: 
\begin{equation}
\widehat{E^{-1}}(u)=\mu^{-1} E(u)-1 \,.
\end{equation}

We also have the following alternative versions of the pentagon relation (for the first one, see~\cite[Eq. (36)]{RW} and the second can be derived from the first using~\eqref{eq:reflection}):
\begin{equation} \label{eq:pentagon23}
\frac{1}{\sqrt{|G|}}\sum_{t\in G} \ang{t}E(t)E(v-t)\ang{t;u}=\frac{E(u)E(v)}{E(u+v)}+\sqrt{|G|}\frac{\delta(u+v)}{\ang{v}}\,,
\end{equation}
\begin{equation} \label{eq:ratiofourier}
\frac{1}{\sqrt{|G|}}\sum_{t\in G} \frac{E(v+t)}{E(t)}\ang{t;-u}=E(v)\Big(\frac{E(u)}{E(u+v)}-1\Big) +\sqrt{|G|}\frac{\delta(u+v)}{\ang{v}}\,.
\end{equation}
It is sometimes convenient to use a Fourier transform of the last relation:
\begin{equation} \label{eq:ratiofourier2}
\frac{E(v+t)}{E(t)}=\frac{E(v)}{\sqrt{|G|}}\sum_{u\in G}\Big(\frac{E(u)}{E(u+v)}-1\Big)\ang{u;t} +\frac{\ang{-v;t}}{\ang{v}}\,.
\end{equation}

Another useful property of finite quantum dilogarithms that is not stated in~\cite{RW} is the following
\begin{equation}\label{eq:correlation}
\frac{1}{\sqrt{|G|}}\sum_{t\in G}E(t)E(t+s) = \frac{\mu^{-1}}{\ang{s}} + \mu^2E(0)\,.	
\end{equation}
This follows from~\eqref{eq:fourier} and~\eqref{eq:reflection} together with Plancherel's identity $\sum_{t}f(t)g(t)=\sum_t\widehat{f}(t)\widehat{g}(-t)$.
We also have the following more complicated identity analogous to a relation discussed by Volkov in~\cite[\S5]{Vol05}.
\begin{lemma} \label{lem:volkov45}
Let $E$ be a finite quantum dilogarithm on~$G$. Then for all $u,v,s\in G$ we have
\begin{align*}
\frac{1}{\sqrt{|G|}}\sum_{t\in G}\frac{E(t+u)E(t+s+v)}{E(t)E(t+s)}
={}&\frac{\ang{u;s}E(v)E(u)E(s+v)}{\ang{u}E(u+v)E(s-u)}-
\frac{\ang{u;s}E(v)E(s)}{\ang{u}E(s-u)}
-\frac{E(u)E(s+v)}{E(s)}\\&+\sqrt{|G|}\frac{\ang{u;s}}{\ang{u}^2}\delta(u+v)+E(0).
\end{align*}
\end{lemma}

\begin{proof}
We only give an outline, leaving the remaining calculations to the reader. The idea is to apply~\eqref{eq:ratiofourier2} to rewrite $E(t+u)/E(t)$ and $E(t+s+v)/E(t+s)$. The main term that one needs to deal with is then
\[\frac{E(u)E(v)}{\sqrt{|G|}}\sum_{r\in G}\frac{E(r)E(-r)}{E(r+u)E(-r+v)}\ang{s;-r} = \frac{E(u)E(v)}{\sqrt{|G|}}\sum_{r\in G}\frac{1}{\ang{r}E(r+u)E(-r+v)}\ang{s;-r} + E(0)\,,\]
which after reflecting $\frac{1}{E(-r+v)}=E(r-v)\ang{r-v}-\sqrt{|G|}\delta(r-v)$ can be evaluated using~\eqref{eq:ratiofourier}.
\end{proof}

\subsection{\texorpdfstring{Finite quantum dilogarithms from hyperbolic 
$\SL_2(\ZZ)$ matrices}{Finite quantum dilogarithms from hyperbolic 
SL2(Z) matrices}}
\label{sub.gammaQDL}

In this section we recall the basic facts of the finite quantum 
dilogarithm of~\cite{RW}. Fix $\gamma \in \SL_2(\BZ)$ with $
\tr(\gamma)=2+N$ with $N>0$. Let 
\begin{equation}
\label{Ggdef}
G_\gamma=\BZ^2/\BZ^2(I-\gamma)
\end{equation}
denote the finite abelian group with $N$ elements. It is a 
metric group with 
\begin{equation}
\langle u \rangle=(-1)^{u_1+u_2+u_1 u_2} \ee_{2N}
(\omega(u\gamma,u)),
\qquad \langle u;v \rangle = \ee_N(\omega(u\gamma-u,v)) 
\end{equation}
where $u=(u_1,u_2) \in \BZ^2$. We define~\cite[Eqn.(2)]{RW}
\begin{equation}
E(u) = F_\gamma(u) \,.
\end{equation}
It is a finite quantum dilogarithm, in particular it satisfies all of the identities in Section~\ref{sub.finiteQDL} above. 
We additionally have $E(0)=\eps^{1/2}$, where $\eps=\eps_{\gamma}>1$ satisfies
\[\eps^2-(N+2)\eps+1=0\] 
and $\ol{E(x)}=\ang{x} E(x)$ (this property may not hold for arbitrary finite quantum dilogarithms).

For the metric group $G_\gamma$ we have
$\mu=\mu_\gamma$, where $\mu_\gamma$ (a $24$th root of unity) 
is the multiplier system of the Dedekind $\eta$-function. 
We define the constants
\begin{equation}
d=\frac{1}{2}(N-\sqrt{N(N+4)}), 
\qquad
\lambda=\frac{N(N+4)-N\sqrt{N(N+4)}}{2} \,.
\end{equation}
Note that $d$ and $\lambda$ have minimal polynomials
$ d^2-dN-N=0$ and $\lambda^2-N(N+4)\lambda+N^2(N+4)=0$. 
If we want to emphasise the dependence 
on $\gamma$, we will insert a $\gamma$
subscript to $\varepsilon$, $N$, $\langle u \rangle$, 
$\langle u;v \rangle$, $d$ or $\lambda$. Note that when 
$\varepsilon_\gamma>1$
\begin{equation}
\label{Gconstants}
N_\gamma=\eps^{-1}(\eps-1)^2
, \qquad 
d_\gamma = \varepsilon^{-1}-1, \qquad 
\sqrt{N_\gamma(N_\gamma+4)}=\varepsilon - \varepsilon^{-1},
\qquad
\lambda_\gamma=\varepsilon^{-2}(\varepsilon-1)^2(\varepsilon+1)
\end{equation}
and $\varepsilon_{\gamma^{-1}}=(\varepsilon_\gamma)^{-1}$.

Note finally that $G_\gamma$ and $G_{\gamma^2}$ are connected
by natural inclusion and projection maps
\begin{equation}
\label{2maps}
\gamma+I: G_{\gamma} \hookrightarrow G_{\gamma^2}, \qquad 
\varpi: G_{\gamma^2} \twoheadrightarrow G_{\gamma}
\end{equation}
with composition $G_{\gamma} \stackrel{2}{\to} G_{\gamma}$.

\subsection{A compatible pair of quantum dilogarithms}
\label{sub.Gpm}

Our goal for the rest of this section is to prove 
Theorem~\ref{thm.second}, which will require establishing several 
identities between finite quantum dilogarithms on the groups 
$G_{\gamma^2}$ and $G_{\gamma}$. It will be convenient to change 
notation and work with the group $G=(1-\gamma^2)^{-1}\ZZ^2/\ZZ^2$, 
isomorphic to $G_{\gamma^2}$. 
Define its subgroups $G_{\pm}$ by
\[G_{+}=\{x\in G\colon x\gamma=x\},
\qquad\qquad G_{-}=\{x\in G\colon x\gamma=-x\}\,.\]
(The group $G_{+}$ is canonically isomorphic to $G_{\gamma}$ via 
multiplication by $1-\gamma$.) We have finite quantum dilogarithms 
$E_2\colon G\to \CC$ and $E_1\colon G_{+}\to \CC$, related by 
$E_2(y)=E_1(y)^2$ for all $y\in G_+$, and in addition 
$E_2(x\gamma)=E_2(x)$ for all $x\in G$. Define Gaussians $
\ang{\cdot}\colon G\to \CC^{\times}$ and $\ang{\cdot}_{\pm}\colon 
G_{\pm}\to \CC^{\times}$ by
\begin{align*}
\ang{x} &= s(x(\gamma^2-1))\,\ee(\omega(x,x\gamma^2)/2)\,,\\ 
\ang{x}_{+} &= s(x(\gamma-1))\,\ee(\omega(x,x\gamma)/2)\,,\\ 
\ang{x}_{-} &= s(x(\gamma+1))\,\ee(\omega(x,-x\gamma)/2)\,,
\end{align*}
and use the same subscripts to denote the corresponding bicharacters 
on $G$, $G_{+}$, and $G_{-}$. Here $s\colon \ZZ^2\to \{\pm1\}$ is 
defined to take value~$1$ on $2\ZZ^2$ and $-1$ otherwise.
It is straightforward to verify that
\[\ang{y}=\ang{y}_{+}^2\,,\qquad \ang{z}=\ang{z}_{-}^2\,,\qquad 
\ang{x+x\gamma}_{+}\ang{x-x\gamma}_{-} = \ang{x}^2\,,\]
and
\[\ang{x;y} = \ang{x+x\gamma;y}_{+}\,,\qquad \ang{x;z} = \ang{x-
x\gamma;z}_{-}\,,\]
for $y\in G_{+}$, $z\in G_{-}$, $x\in G$.
We also denote $\mu=\frac{\widehat{E_1}(0)}{E_1(0)}$, so that $
\mu^2=\frac{\widehat{E_2}(0)}{E_2(0)}$. Note that $|G_{+}|=N$, $|
G_{-}|=N+4$, and $|G|=N(N+4)$. We also have $G_{\pm}^{\perp}=G_{\mp}$ 
(with respect to $\ang{\cdot;\cdot}$), but in general $G_{+}\cap 
G_{-}$ is the 2-torsion part of $G_{+}$. In particular, 
$G=G_{+}\oplus G_{-}$ if and only if~$N$ is odd. 

Using~\cite[Prop. 1]{RW} and the fact that $\Phi_{\gamma,1,1}(0;\tau)=\mu_\gamma^{-1} \tq^{\frac{1}{24}}q^{-\frac{1}{24}}$, we can show that for $z,\tau\in\RR$ 
\[
\frac{\Phi_{\gamma,1,1}(z;\tau)}{|\Phi_{\gamma,1,1}(z;\tau)|}
=
\mu_\gamma^{-1}\ee(Q_{\gamma,0,0}(z;\tau)/2)\,.
\]
The following useful calculation is an easy consequence of this formula:
\begin{equation}\label{eq:e2gminus}
E_2(z)=-\ang{z}_{-}^{-1}\,,\qquad z\in G_-\sm\{0\}\,.
\end{equation}
By Poisson summation together with~\eqref{eq:correlation} we have
\[\frac{1}{\sqrt{N+4}}\sum_{z\in G_{-}}\widehat{E_2}(z) = \frac{1}
{\sqrt{N}}\sum_{y\in G_{+}}E_2(y) = \frac{1}{\sqrt{N}}\sum_{y\in 
G_{+}}E_1(y)^2 = \mu^{-1}+\mu^{2}E_1(0)\,.\]
By~\eqref{eq:e2gminus} this gives an evaluation of the Gauss 
sum for $\ang{\cdot}_{-}$:
\begin{equation} \label{eq:gausssumminus}		
\frac{1}{\sqrt{N+4}}\sum_{z\in G_{-}}\ang{z}_{-}^{\pm1} 
= -\mu^{\mp3}\,.
\end{equation}

Also recall that $E_1(0)^2=\sqrt{N}E_1(0)+1$ and 
$E_2(0)^2=\sqrt{N(N+4)}E_2(0)+1=(N+2)E_2(0)-1$. We will also freely 
use the following simple corollaries:
\[\frac{1+E_2(0)^{-1}}{\sqrt{N+4}}=E_1(0)^{-1}\,,\qquad\qquad 
\frac{1+E_2(0)}{\sqrt{N+4}}=E_1(0)\,.\]

\subsection{\texorpdfstring{Identities between $E_2$ and $E_1$}{Identities between E2 and E1}}
Throughout this section, we assume that $|G|$ is odd, though 
this is only used in the proof of Equation~\eqref{eq:e2gauss}.

\begin{lemma} For any $x,y\in G$, we have
\begin{equation} \label{eq:e2beta}
\frac{1}{N+4}\sum_{z,w\in G_{-}}\ang{z;-w}_{-}\frac{E_2(x+z)E_2(y+w)}
{E_2(x+y+z+w)}= \frac{E_2(x)E_2(y)}{E_2(x+y\gamma)}\,.
\end{equation}
\end{lemma}
	
\begin{proof}
We have the following calculation using~\eqref{eq:pentagon23}:
\begin{align*}
\frac{1}{\sqrt{N(N+4)^3}}\sum_{\substack{z,w\in G_{-}\\t\in 
G}}\ang{z;-w}_{-}\ang{t}&E_2(t)E_2(y+w-t)\ang{t;x+z} = \frac{1}
{\sqrt{N(N+4)}}\sum_{t\in G}\ang{t}E_2(t)E_2(y-t\gamma)\ang{t;x}\\
&= \frac{E_2(y\gamma)E_2(x)}{E_2(x+y\gamma)}+
\sqrt{N(N+4)}\delta(x+y\gamma)\ang{x}^{-1}\,.
\end{align*}
Also note that
\begin{align*}
\frac{1}{N+4}\sum_{\substack{z,w\in G_{-}\\ z+w=-x-y}}
\ang{z;-w}
_{-}\ang{x+z}^{-1}
&= \frac{1}{N+4}\sum_{z\in G_{-}}\ang{z;x+y}_{-}\ang{z}
\ang{x+z}^{-1} 
\\
&= \frac{1}{N+4}\ang{x}^{-1}\sum_{z\in G_{-}}
\ang{z;x+y}_{-}\ang{-x;z} = \delta(x+y\gamma)\ang{x}^{-1}\,.
\end{align*}
(If $x+y\not\in G_-$, both sides vanish.)
Combining these two evaluations together 
with~\eqref{eq:pentagon23} gives~\eqref{eq:e2beta}.
\end{proof}

The key calculation is the following. 
\begin{lemma} For all $x\in G$ and $w\in G_{-}$ we have
\begin{equation} \label{eq:e2gauss}
\frac{1}{\sqrt{N+4}}\sum_{z\in G_{-}}\ang{z}_{-}E_2(x+z)
\ang{z;-w}_{-} = - \frac{E_1(x+x\gamma)}{\ang{w}_{-}E_2(x+w)}\,.
\end{equation}
\end{lemma}

\begin{proof} By~\eqref{eq:e2beta} we have
\[\frac{1}{N+4}\sum_{z,t\in G_{-}}
\ang{z;z-w-t}_{-}\frac{E_2(x+z)E_2(t-z)}{E_2(x+t)}
= \frac{E_2(x)E_2(-w)}{E_2(x+w)}\]
Writing $E_2(t-z)=-\ang{t-z}_{-}^{-1}+(E_2(0)+1)\delta(t-z)$ and 
$E_2(-w)=-\ang{-w}_{-}^{-1}+(E_2(0)+1)\delta(-w)$ gives
\[-\frac{1}{N+4}\sum_{z,t\in G_{-}}
\frac{\ang{z;z-w-t}_{-}}{\ang{t-z}_{-}}\frac{E_2(x+z)}{E_2(x+t)}
= -\frac{E_2(x)}{\ang{w}_{-}E_2(x+w)}\,.\]
Since $\frac{\ang{z;z-w-t}_{-}}{\ang{t-z}_{-}}=\frac{\ang{z}_{-}}
{\ang{t}_{-}}\ang{z;-w}$, we can rewrite the last identity as
\[\Big(\frac{1}{\sqrt{N+4}}\sum_{z\in G_{-}}\ang{z}_{-}
E_2(x+z)\ang{z;-w}_{-}\Big)\Big(\frac{1}{\sqrt{N+4}}\sum_{t\in 
G_{-}}\frac{1}{\ang{t}_{-}E_2(x+t)}\Big) = \frac{E_2(x)}{\ang{w}_{-}
E_2(x+w)}\,.\]
Writing $f_x(z)=\ang{z}_{-}E_2(x+z)$, we see that $\widehat{f_x}
(w)f_x(w)$ is constant in~$w$ (we write $\widehat{f}$ for the Fourier 
transform on $G_{-}$). On the other hand, we have
\[\widehat{f_x}(w)f_x(w) = \frac{1}{\sqrt{N+4}}\sum_{z\in 
G_{-}}\ang{z-w}_{-}E_2(x+z)E_2(x+w)=\widehat{f_{x+w}}(0)f_{x+w}(0)\,,
\]
so $\widehat{f_x}(w)f_x(w)$ depends only on the coset of $x$ modulo 
$G_{-}$. In other words $\widehat{f_x}(w)f_x(w)=-R(x+x\gamma)$ for 
some $R\colon G_{+}\to \CC$, and it suffices to show that $R=E_1$.

In view of~\eqref{eq:e2gminus}, we have $f_x(z)=-\frac{E_2(x+z)}
{E_2(z)}+(1+E_2(0)^{-1})E_2(x)\delta(z)$. 
Using~\eqref{eq:ratiofourier} together with Poisson summation on 
$G_{-}$ then gives
\[-\widehat{f_x}(w)+E_1(0)E_2(x) = \frac{1}{\sqrt{N}}\sum_{a\in 
G_+}\frac{E_2(x)E_2(w+a)}{E_2(x+w+a)}+\sqrt{N+4}\frac{[x+w\in G_{+}]}
{\ang{x}}\,.\]
Taking $w=0$ and $x\in G_{+}$, we get that
\[-\widehat{f_x}(0)+E_1(0)E_2(x) = \frac{1}{\sqrt{N}}\sum_{a\in 
G_+}\frac{E_2(x)E_2(a)}{E_2(x+a)}+\frac{\sqrt{N+4}}{\ang{x}}\,.\]
Since $x\in G_{+}$, in this case we have $E_2(a)=E_1(a)^2$ and 
$E_2(x+a)=E_1(x+a)^2$, and using Lemma~\ref{lem:volkov45} we can 
compute the sum on the right hand side as
\[\frac{1}{\sqrt{N}}\sum_{a\in G_+}\frac{E_1(x)^2E_1(a)^2}
{E_1(x+a)^2} = E_1(x)^2\Big(\frac{E_1(-x)^4}{E_1(-2x)}-E_1(-
x)^2(E_1(0)+E_1(0)^{-1})+\sqrt{N}\frac{\delta(2x)}{\ang{x}
^2}+E_1(0)\Big)\]
Therefore,
\[\frac{R(2x)}{E_2(x)} = E_1(x)^2\Big(\frac{E_1(-x)^4}{E_1(-2x)}-E_1(-
x)^2(E_1(0)+E_1(0)^{-1})+\sqrt{N}\frac{\delta(2x)}{\ang{x}^2}\Big)+
\frac{\sqrt{N+4}}{\ang{x}}\,.\]
For $x=0$ this easily implies $R(0)=E_1(0)$, so assume that $x\ne0$. 
Then
\begin{align*}
R(2x)&=\frac{\ang{x}^{-2}}{E_1(-2x)}-\ang{x}^{-1}
E_1(x)^2(E_1(0)+E_1(0)^{-1})+E_1(x)^4\sqrt{N}\frac{\delta(2x)}{\ang{x}
^2}+E_1(x)^2\frac{\sqrt{N+4}}{\ang{x}}\\
&=\frac{\ang{x}^{-2}}{E_1(-2x)}+E_1(x)^4\sqrt{N}\frac{\delta(2x)}
{\ang{x}^2}
= E_1(2x)+(E_1(x)^4-1)\sqrt{N}\frac{\delta(2x)}{\ang{x}^2} = 
E_1(2x)\,,
\end{align*}
since for $2x=0$, $x\ne 0$ we have $E_1(x)^4=E_1(x)^2E_1(-x)^2=\ang{x}
_{+}^{-2}=1$. This proves that $R=E_1$ on $2G_{+}$, proving the lemma 
since~$N$ is odd.
\end{proof}
Note that by Fourier inversion~\eqref{eq:e2gauss} also implies
\begin{equation} \label{eq:e2gaussfourier}
-\frac{E_1(x+x\gamma)}{\sqrt{N+4}}\sum_{t\in G_{-}}\frac{\ang{t;v}
_{-}}{\ang{t}_{-}E_2(x+t)} 
= \ang{v}_{-}E_2(x+v)\,.
\end{equation}

\begin{remark}
The above lemma implies that the function $f\colon G_{-}\to\CC$ 
given by $f(z)=E_2(x+z)\ang{z}_{-}$ satisfies $\widehat{f}(w)=-
E_1(x+x\gamma)/f(w)$. In particular, taking the unitary Galois 
conjugates of $E_1$ and $E_2$ gives examples of functions $g\colon 
G_{-}\to \CC$ with $|g(w)|=|\widehat{g}(w)|=1$ for all $w\in G_{-}$, 
with $\widehat{g}$ proportional to $g^{-1}=\overline{g}$.
\end{remark}

\begin{lemma} For all $u,v\in G$ we have
\begin{equation} \label{eq:e1e2beta}
\frac{1}{\sqrt{N+4}}\sum_{w\in G_{-}}
\ang{w}_{-}\frac{E_2(u+w)E_2(v+w)}{E_2(u+v+w)}
= -\frac{E_1(u+u\gamma)E_1(v+v\gamma)E_2(u+v\gamma)}
{E_1(u+v+u\gamma+v\gamma)E_2(u)E_2(v)}.
\end{equation}
\end{lemma}

\begin{proof}
Applying~\eqref{eq:e2gauss} to $E_1(u+u\gamma)/E_2(u)$ and 
$E_1(v+v\gamma)/E_2(v)$ gives that the right hand side is equal to
\[-\frac{E_2(u+v\gamma)}{(N+4)E_1(u+v+(u+v)\gamma)}\sum_{z_1,z_2\in 
G_{-}}\ang{z_1}_{-}\ang{z_2}_{-}E_2(u+z_1)E_2(v+z_2)\,.\]
Next, multiply and divide the inner sum by $\ang{z_1+z_2}_{-}
E_2(u+v+z_1+z_2)$ and apply~\eqref{eq:e2gaussfourier}. This gives
\[\frac{E_2(u+\eps v)}{(N+4)^{3/2}}\sum_{w\in G_{-}}\frac{1}{\ang{w}
_{-}E_2(u+v+w)}\sum_{z_1,z_2\in G_{-}}\frac{\ang{z_1;-z_2}
_{-}\ang{z_1+z_2;w}_{-}E_2(u+z_1)E_2(v+z_2)}{E_2(u+v+z_1+z_2)}\,.\]
Since $\ang{z_1;-z_2}_{-}\ang{z_1+z_2;w}_{-}=\ang{z_1-w;w-z_2}
_{-}\ang{w}_{-}^2$, applying~\eqref{eq:e2beta} with $x=u+w$, $y=v+w$ 
gives the claim.
\end{proof}

\begin{proposition} \label{prop:izumi4}
For $x\in G\sm G_{+}$, and $h,k\in G_{+}$ we have
\begin{equation} \label{eq:izumigamma_4}
\frac{1}{\sqrt{N}}\sum_{y\in G_{+}}\frac{E_1(y-h)E_1(y-k)}{E_2(x+y)} 
= E_1(0)+\frac{E_1(h+k+x+x\gamma)}{\ang{h-k}_{+}E_2(x+h)E_2(x+k)}\,.
\end{equation}
\end{proposition}

\begin{proof}
By shifting $x$, it is enough to prove 
\begin{equation} \label{eq:need1}
\frac{E_2(x)}{\sqrt{N}}\sum_{y\in G_{+}}\frac{E_1(y)E_1(y+a)}
{E_2(x+y)} = 
E_1(0)E_2(x)+\frac{E_1(x+x\gamma-a)}{\ang{a}_{+}E_2(x-a)}\,.
\end{equation}

First, since $E_2(y)=E_1(y)^2$, we can rewrite the sum on the 
left of~\eqref{eq:need1} as
\[\frac{E_2(x)}{\sqrt{N}}\sum_{y\in G_{+}}\frac{E_1(y+a)E_2(y)}
{E_1(y)E_2(x+y)}\,.\]
Since $x\not \in G_+$, by~\eqref{eq:ratiofourier} we have
\[\frac{E_2(y)}{E_2(x+y)} = 1+\frac{1}{\sqrt{|G|}}\sum_{t\in 
G}\frac{E_2(t+x)}{E_2(t)E_2(x)}\ang{t;-y}\,.\]
Plugging this into the above sum gives, by~\eqref{eq:ratiofourier}
\begin{align*}
&\frac{E_2(x)}{\sqrt{N}}\sum_{y\in G_{+}}\frac{E_1(y+a)}{E_1(y)}+
\frac{1}{\sqrt{N|G|}}\sum_{t\in G}\sum_{y\in G_{+}}\frac{E_1(y+a)}
{E_1(y)}\ang{t;-y}\frac{E_2(t+x)}{E_2(t)}\\
={}&\frac{E_2(x)}{\sqrt{N}}\sum_{y\in G_{+}}\frac{E_1(y+a)}{E_1(y)}+
\frac{1}{\sqrt{|G|}}\sum_{t\in G}\Big(E_1(a)\Big(\frac{E_1(t+t\gamma)}
{E_1(a+t+t\gamma)}-1\Big)+\sqrt{N}\delta(a+t+t\gamma)\ang{a}_{+}
^{-1}\Big)\frac{E_2(t+x)}{E_2(t)}\,.
\end{align*}
For $a=0$ the last expression is equal to
\[E_2(x)\sqrt{N}+\frac{1}{\sqrt{N+4}}\sum_{t\in G_{-}}\frac{E_2(t+x)}
{E_2(t)} = E_2(x)\Big(\sqrt{N}+\frac{E_2(0)^{-1}+1}{\sqrt{N+4}}\Big)-
\frac{1}{\sqrt{N+4}}\sum_{t\in G_{-}}E_2(t+x)\ang{t}_{-}\,,\]
which is equal to the right hand side of~\eqref{eq:need1} 
by~\eqref{eq:e2gauss}.

So we may assume that $a\ne 0$. Let us choose $b$ so that 
$a=b+b\gamma$ (note that $b\not\in G_{-}$ since $a\ne 0$). 
Using~\eqref{eq:ratiofourier} we calculate:
\[\frac{E_2(x)}{\sqrt{N}}\sum_{y\in G_{+}}\frac{E_1(y+a)}{E_1(y)} - 
\frac{1}{\sqrt{|G|}}\sum_{t\in G}E_1(a)\frac{E_2(t+x)}{E_2(t)} = 
E_2(x)(E_1(0)-E_1(a)) - E_1(a)(E_2(0)-E_2(x))\,.\]
Next, we directly rewrite
\[\frac{\ang{a}_{+}^{-1}}{\sqrt{N+4}}\sum_{t\in G}\frac{E_2(t+x)}
{E_2(t)}\delta(t+t\gamma+a) = \frac{\ang{a}_{+}^{-1}}
{\sqrt{N+4}}\sum_{w\in G_{-}}\frac{E_2(x+w-b)}{E_2(w-b)}\,.\]
Therefore, it remains to show that
\begin{equation} \label{eq:need2}
\frac{1}{\sqrt{|G|}}\sum_{t\in G}\frac{E_1(a)E_1(t+t\gamma)}
{E_1(a+t+t\gamma)}\frac{E_2(t+x)}{E_2(t)} = E_1(a)E_2(0) + 
\frac{E_1(x+x\gamma-a)}{\ang{a}_{+}E_2(x-a)} - \frac{\ang{a}_{+}^{-1}}
{\sqrt{N+4}}\sum_{w\in G_{-}}\frac{E_2(x+w-b)}{E_2(w-b)}
\end{equation}
Applying~\eqref{eq:e1e2beta} we can write
\begin{align*}
\frac{1}{\sqrt{|G|}}\sum_{t\in G}\frac{E_1(a)E_1(t+t\gamma)}
{E_1(a+t+t\gamma)}&\frac{E_2(t+x)}{E_2(t)}
= \frac{E_2(b)}{\sqrt{|G|}}\sum_{t\in G}\frac{E_1(a)E_1(t+t\gamma)}
{E_1(a+t+t\gamma)}\frac{E_2(t+b\gamma)}{E_2(b)E_2(t)}\frac{E_2(t+x)}
{E_2(t+b\gamma)}\\
= &-\frac{E_2(b)}{\sqrt{|G|(N+4)}}\sum_{t\in G}\frac{E_2(t+x)}
{E_2(t+b\gamma)}\sum_{w\in G_{-}}\ang{w}_{-}\frac{E_2(t+w)E_2(b+w)}
{E_2(t+b+w)}\\
= &-\frac{E_2(b)}{\sqrt{N+4}}
\sum_{w\in G_-}\ang{w}_{-}E_2(b+w)
\Big(\frac1{\sqrt{|G|}}\sum_{r\in G}\frac{E_2(r+x-b\gamma)E_2(r+w-
b\gamma)}{E_2(r)E_2(r+b+w-b\gamma)}\Big)\,.\end{align*}
Next, applying Lemma~\ref{lem:volkov45} with $u=w-b\gamma$, $v=x-b-
w$, $s=b-b\gamma+w$ shows that the inner sum is
\[\frac{\ang{b-b\gamma+w}E_2(x-b-w)}{\ang{b}E_2(b)}\Big(\frac{E_2(w-
b\gamma)E_2(x-b\gamma)}{E_2(x-b-b\gamma)}-E_2(b-b\gamma+w)\Big)
-\frac{E_2(w-b\gamma)E_2(x-b\gamma)}{E_2(b-b\gamma+w)}+E_2(0).\]
Substituting this expression for the inner sum, we expand and 
evaluate the four resulting sums. First, by~\eqref{eq:e2gauss} we 
have 
\[-\frac{E_2(b)}{\sqrt{N+4}}\sum_{w\in G_-}\ang{w}_{-}E_2(b+w)E_2(0) 
= E_2(0)E_1(a)\,.\]
Next, using $E_2(b+w)E_2(w-b\gamma) = \ang{b+w}^{-1}$ and 
\[E_2(s)^{\pm1} = -\ang{s}_{-}^{\mp1}+(E_2(0)^{\pm1}+1)\delta(s)\] 
we calculate
\begin{align*}
\frac{E_2(b)}{\sqrt{N+4}}&\sum_{w\in G_-}\ang{w}_{-}
E_2(b+w)\frac{E_2(w-b\gamma)E_2(x-b\gamma)}{E_2(b-b\gamma+w)} = 
\frac{E_2(b)E_2(x-b\gamma)}{\sqrt{N+4}}\sum_{w\in G_-}\ang{w}
_{-}\frac{\ang{b+w}^{-1}}{E_2(b-b\gamma+w)}\\
=&-\frac{E_2(b)E_2(x-b\gamma)}{\sqrt{N+4}}\sum_{w\in G_-}\ang{w}
_{-}\frac{\ang{b-b\gamma+w}_{-}}{\ang{b+w}}+\frac{E_2(b)E_2(x-
b\gamma)\ang{b\gamma-b}_{-}}{\ang{b}E_1(0)}\\
=&E_2(b)E_2(x-b\gamma)\frac{\ang{b\gamma-b}_{-}}{\ang{b}}\big(-
\sqrt{N+4}+E_1(0)^{-1}\big)
= -E_1(0)\frac{E_2(b)E_2(x-b\gamma)\ang{b\gamma-b}_{-}}{\ang{b}}
\end{align*}
and
\begin{align*}
\frac{1}{\sqrt{N+4}}\sum_{w\in G_-}&\ang{w}_{-}E_2(b+w)\frac{\ang{b-
b\gamma+w}E_2(b-b\gamma+w)E_2(x-b-w)}{\ang{b}}\\
=&-\frac{1}{\sqrt{N+4}}\sum_{w\in G_-}\ang{w}_{-}\frac{\ang{b-
b\gamma+w}_{-}E_2(x-b-w)}{\ang{b}\ang{w+b}E_2(-w-b)} + 
E_1(0)\frac{E_2(b)E_2(x-b\gamma)\ang{b\gamma-b}_{-}}{\ang{b}} \\
=&-\frac{\ang{a}_{+}^{-1}}{\sqrt{N+4}}\sum_{w\in G_-}\frac{E_2(x+w-b)}
{E_2(w-b)} + E_1(0)\frac{E_2(b)E_2(x-b\gamma)\ang{b\gamma-b}_{-}}
{\ang{b}}\,.
\end{align*}
Finally, we have
\begin{align*}
-\frac{1}{\sqrt{N+4}}&\frac{E_2(x-b\gamma)}{E_2(x-b-
b\gamma)}\sum_{w\in G_-}\ang{w}_{-}\frac{\ang{b-b\gamma+w}E_2(x-b-
w)E_2(b+w)E_2(w-b\gamma)}{\ang{b}} \\
&=-\frac{1}{\sqrt{N+4}}\frac{E_2(x-b\gamma)}{E_2(x-b-
b\gamma)}\sum_{w\in G_-}\ang{w}_{-}\ang{b+w;-b\gamma}E_2(x-b-w)\\
&=-\frac{\ang{b;-b\gamma}}{\sqrt{N+4}}\frac{E_2(x-b\gamma)}{E_2(x-b-
b\gamma)}\sum_{w\in G_-}\ang{w}_{-}\ang{w;b-b\gamma}_{-}E_2(x-b-w) = 
\frac{E_1(x+x\gamma-a)}{\ang{a}_{+}E_2(x-a)}\,,
\end{align*}
where we used~\eqref{eq:e2gauss} and $\frac{\ang{b;-b\gamma}}{\ang{b-
b\gamma}_{-}} = \frac{\ang{b}^2}{\ang{b+b\gamma}\ang{b-b\gamma}_{-}} 
= \frac{1}{\ang{b+b\gamma}_{+}} = \frac{1}{\ang{a}_{+}}$. Combining 
these evaluations together gives~\eqref{eq:need2} and hence the claim.
\end{proof}

It is easy to derive also the following identities.
\begin{proposition}
For $x\in G\sm G_{+}$ and $h\in G_{+}$ we have
\begin{equation} \label{eq:izumigamma_1}
\frac{\mu^2}{\sqrt{N}}\sum_{y\in G_{+}}\frac{1}{E_2(x+y)} = E_1(0)+
\mu^3\frac{\ang{x}^2}{\ang{x+x\gamma}_{+}}\,,
\end{equation}
\begin{equation} \label{eq:izumigamma_3}
\frac{\mu}{\sqrt{N}}\sum_{y\in G_{+}}\frac{E_1(y+h)}{E_2(x+y)} = 
E_1(0)+\mu\frac{\ang{x}^2}{\ang{x+x\gamma}_{+}}E_2(h-x\gamma)\,.
\end{equation}        
\end{proposition}

\begin{proof}
By $\widehat{E_2^{-1}}=\mu^{-2}E_2-1$ and Poisson summation we have
\[\frac{\mu^2}{\sqrt{N}}\sum_{y\in G_{+}}\frac{1}{E_2(x+y)} = \frac{1}
{\sqrt{N+4}}\sum_{z\in G_{-}}(E_2(z)-\mu^2)\ang{z;x-x\gamma}_{-} = 
E_1(0)-\sum_{z\in G_{-}}\ang{z}_{-}^{-1}\ang{z;x-x\gamma}_{-}\,,\]
which gives~\eqref{eq:izumigamma_1} by~\eqref{eq:gausssumminus}.
	
For~\eqref{eq:izumigamma_3}, by shifting~$x$ it suffices to prove it 
in the case $h=0$. In this case we can rewrite the identity as
\[\frac{\mu}{\sqrt{N}}\sum_{y\in G_{+}}E_1(y)^{-1}\frac{E_2(y)}
{E_2(x+y)} = E_1(0)+\mu\ang{x-\eps x}_{-}E_2(-x)\,.\]
Expanding $E_2(x)/E_2(x+y)$ as in the proof of 
Proposition~\ref{prop:izumi4} and using $\widehat{\mu E_1^{-1}}=E_1-
\mu$ gives that the sum on the left hand side is equal to
\[E_1(0)-\mu+\frac{1}{\sqrt{|G|}}\sum_{t\in G}(E_1(t+t\gamma)-
\mu)\frac{E_2(t+x)}{E_2(t)E_2(x)}
=E_1(0)-\mu\frac{E_2(0)}{E_2(x)}+\frac{E_2(x)^{-1}}{\sqrt{|
G|}}\sum_{t\in G}E_1(t+\eps t)\frac{E_2(t+x)}{E_2(t)}\,.\]
Next, using~\eqref{eq:e2gauss} and~\eqref{eq:correlation} we compute
\begin{align*}
\frac{1}{\sqrt{|G|}}\sum_{t\in G}E_1(t+t\gamma)\frac{E_2(t+x)}
{E_2(t)} 
&= -\frac{1}{(N+4)\sqrt{N}}\sum_{\substack{t\in G\\ z\in 
G_{-}}}\ang{z}_{-}E_2(t+z)E_2(t+x)\\
&= -\frac{1}{\sqrt{N+4}}\sum_{z\in G_{-}}\ang{z}_{-}(\mu^{-2}\ang{z-x}
^{-1}+\mu^4E_2(0))
= \mu\Big(E_2(0)+\frac{\ang{x-x\gamma}_{-}}{\ang{x}}\Big)\,,
\end{align*}
implying~\eqref{eq:izumigamma_3}, where we used
\[\ang{z}_{-}\ang{z-x}^{-1} = \ang{x}^{-1}\ang{z}_{-}^{-1}\ang{z;x-
x\gamma}_{-} = \frac{\ang{x-x\gamma}_{-}}{\ang{x}\ang{z-x+x\gamma}
_{-}}\,.\qedhere\]
\end{proof}

The remaining Izumi identity
\begin{equation} \label{eq:izumigamma_2}
\frac{1}{E_2(x+h)E_2(-x\gamma-h)} = \ang{x}\frac{\ang{h}
_{+}\ang{x+x\gamma+h}_{+}}{\ang{x+x\gamma}_{+}}\,\qquad x\not\in 
G_{+},\; h\in G_+
\end{equation}
follows from reflection for $E_2$, since 
$E_2(-x\gamma-h)=E_2(-x-h)$ and
\[\ang{x}\frac{\ang{h}_{+}\ang{x+x\gamma+h}_{+}}{\ang{x+x\gamma}_{+}} 
= \ang{x}\ang{h}\ang{x+x\gamma;h}_{+} = \ang{x+h}\,.\]
Together~\eqref{eq:izumigamma_4}, \eqref{eq:izumigamma_1}, 
\eqref{eq:izumigamma_3}, and~\eqref{eq:izumigamma_2} give solutions 
to Izumi's equations defining the Drinfeld center
of the fusion category $\calC_\gamma$ completing the proof of 
Theorem~\ref{thm.second}.
\qed

\begin{remark}
\label{rem.izumi}
To match with Izumi's notation of~\cite{Izu01} p.~637 
and specifically Equations (6.18)--(6.21), for $g\in G_+$ and 
$u\in G\sm G_+$ one can use
\begin{equation}
\label{omtau}
\tau_u=-u(\gamma+I)\,,\qquad \omega_u =
\frac{\ang{u}_{\gamma^2}}{\ang{\tau_u}_\gamma}
\end{equation}
and 
\begin{equation}
\begin{aligned}
\omega & =\omega_u, \quad
a(g)=\ang{g}_\gamma, \quad
b(g)=\frac{E_1(g)}{\sqrt{N}}, \quad
\hat{a}(0)=\mu_\gamma^{-3},\quad \\
-\frac{n}{\sqrt{d}}& =\sqrt{\eps},\quad
c=\mu_\gamma,\quad
\hat{\mu}(0)=\ang{\tau_u}_\gamma^{-1}\mu_\gamma^{-3},\quad
\xi_u(g)=E_2(u+g)^{-1}\mu_\gamma^2 \,.
\end{aligned}
\end{equation}
\end{remark}

\subsection{\texorpdfstring{The $S$-matrix}{The S-matrix}}
\label{sub.Smat}

We end this section with a proof of Theorem~\ref{thm.S}.
Following Izumi and the notation of Section~\ref{sub.Gpm}, 
we find that
\[
\begin{aligned}
S_{\fd_j,\fd_\ell}
&=
\frac{\omega_{j}\omega_\ell}{\lambda}
\sum_{g\in G_\gamma}
\ang{\varpi(j+\ell);g}_\gamma\ang{g;g}_\gamma^{-1}\\
&\quad+\frac{d_\lambda\omega_j\omega_\ell\ang{\varpi(j)}
_\gamma\ang{\varpi(\ell)}
_\gamma\mu_\gamma^2}{N\lambda}
\sum_{g,h\in G_\gamma}
\ang{\varpi(\ell-j);h-g}_\gamma
\ang{h-g;h-g}_\gamma F_{\gamma^2}(j+g(\gamma+I))F_{\gamma^2}
(\ell+h(\gamma+I))\\
&=
\frac{(1-\eps^{-1})^{-1}\ang{j}\ang{\ell}}{\ang{j(1+\eps)}_+\ang{\ell(1+\eps)}_+\sqrt{N(N+4)}}
\sum_{g\in G_+}
\ang{j+\ell;g}\ang{g}^{-1}\\
&\qquad-\frac{\ang{j}\ang{\ell}\mu^2}{N^{3/2}(N+4)^{1/2}}
\sum_{g,h\in G_+}
\ang{\ell-j;h-g}
\ang{h-g} E_2(j+g)E_2 (\ell+h) \,.
\end{aligned}
\]
We then compute
\[
\begin{aligned}
&-\frac{\mu^2\ang{j}\ang{\ell}}{N}\sum_{g,h\in G_+}
\ang{\ell-j;h-g}
\ang{h-g} E_2(j+g)E_2 (\ell+h)\\
&=
-\frac{\mu^2\ang{j}\ang{\ell}}{N(N+4)^2}
\sum_{k,k'\in G_-}
\sum_{g,h\in G}
\ang{k;g}
\ang{k';h}
\ang{\ell-j;h-g}
\ang{h-g} E_2(j+g)E_2 (\ell+h)\\
&=
-\frac{\ang{j;\ell}\ang{j}\mu^2}{N(N+4)^2}
\sum_{k,k'\in G_-}
\ang{k';\ell}^{-1}
\sum_{g\in G}
\ang{g}
\ang{g;j+k}
E_2(j+g)
\sum_{h\in G}
\ang{h}
\ang{k'-j-g;h}
E_2(h)\\
&=
-\frac{\ang{j;\ell}\ang{j}}{N^{1/2}(N+4)^{3/2}}
\sum_{k,k'\in G_-}
\ang{k';\ell}^{-1}
\sum_{g\in G}
\ang{g}
\ang{g;j+k}
E_2(j+g)
E_2(k'-j-g)\\
&=
-\frac{\ang{j;\ell}}{N^{1/2}(N+4)^{3/2}}
\sum_{k,k'\in G_-}
\ang{k';\ell}^{-1}
\ang{k;j}^{-1}
\sum_{g\in G}
\ang{g}
\ang{g;k}
E_2(g)
E_2(k'-g)\\
&=
-\frac{\ang{j;\ell}}{(N+4)}
\sum_{k,k'\in G_-}
\ang{k';\ell}^{-1}
\ang{k;j}^{-1}
\ang{k+k'}
E_2(k)E_2(k')E_2(-k-k')
\\
&\qquad
-
\frac{\ang{j;\ell}}{N+4}(\sqrt{N(N+4)}+\eps-\eps^3)\\
&=
\frac{\ang{j;\ell}}{(N+4)}
\sum_{k,k'\in G_-}
\ang{k';\ell}^{-1}
\ang{k;j}^{-1}
\ang{k;k'}_-
-(\eps+1)\frac{\ang{j;\ell}}{(N+4)}
\sum_{k\in G_-}
\ang{k;j}^{-1}\\
&\qquad-(\eps+1)\frac{\ang{j;\ell}}{(N+4)}
\sum_{k'\in G_-}
\ang{k';\ell}^{-1}
-(\eps+1)\frac{\ang{j;\ell}}{(N+4)}
\sum_{k\in G_-}
\ang{k;\ell-j}
\ang{k}_-^{-2}\\
&\qquad-\frac{\ang{j;\ell}}{(N+4)}(\sqrt{N(N+4)}-2\eps-2)\\
&=
\ang{j;\ell}
\ang{j(1-\eps);\ell(1-\eps)}_-^{-1}
-(\eps+1)\frac{\ang{j;\ell}}{(N+4)}
\sum_{k\in G_-}
\ang{k;\ell-j}
\ang{k}_-^{-2}+\ang{j;\ell}\\
&=
\frac{\ang{j(1+\eps);\ell(1+\eps)}_+}{\ang{j;\ell}}
-(\eps+1)\frac{\ang{j;\ell}}{(N+4)}
\sum_{k\in G_-}
\ang{k;\ell-j}
\ang{k}^{-1}+\ang{j;\ell} \,.
\end{aligned}
\]
Theorem~\ref{thm.S} then follows from the above 
together with the following identity
\[
(\eps+1)\frac{\ang{j;\ell}}{(N+4)}
\sum_{k\in G_-}
\ang{j-\ell;k}
\ang{k}^{-1}
=
\frac{(1-\eps^{-1})^{-1}\ang{j}\ang{\ell}}{\ang{j(1+\eps)}_+
\ang{\ell(1+\eps)}_+}
\sum_{g\in G_+}
\ang{j+\ell;g}\ang{g}^{-1}\,,
\]
which can be proved by computing the argument and the absolute values 
(where we note that the $2$-torsion of $G_+$ is the same as $G_-$).
\qed

\begin{remark}
\label{rem.ST} For completeness, we reproduce 
the $S$ and $T$ matrices of Izumi with the explicit 
relation to $\gamma^2$.
Noting that $S^t=S$ and $T$ is diagonal, we have
\[
\begin{aligned}
S_{\fa_g,\fa_h}&=\frac{\ang{g;h}_\gamma^2}{\lambda_\gamma}\quad &
S_{\fa_g,\fb_h}&=\frac{(d_\gamma+1)\ang{g;h}_\gamma^2}
{\lambda_\gamma}\quad &
S_{\fa_g,\fc_{h,k}}&=\frac{(d_\gamma+2)\ang{g;h+k}_\gamma}
{\lambda_\gamma}\\
S_{\fa_g,\fd_{j}}&=\frac{d_\gamma\ang{g;\varpi(j)}_\gamma}
{\lambda_\gamma}\quad &
S_{\fb_g,\fb_h}&=\frac{\ang{g;h}_\gamma^2}{\lambda_\gamma}\quad &
S_{\fb_g,\fc_{h+k}}&=\frac{\ang{g;h+k}_\gamma}{\lambda_\gamma}\\
S_{\fb_g,\fd_{j}}&=-\frac{d_\gamma\ang{g;\varpi(j)}_\gamma}
{\lambda_\gamma}\quad &
S_{\fc_{h+k},\fd_{j}}&=0\quad &
S_{\fd_j,\fd_\ell}&=\frac{\ang{j;\ell}_{\gamma^2}+\ang{j;\ell}
_{\gamma^2}^{-1}\ang{\varpi(j);\varpi(\ell)}_{\gamma}}
{\sqrt{N_\gamma(N_\gamma+4)}}\,.
\end{aligned}
\]
\[S_{\fc_{h+k},\fc_{h'+k'}}=\frac{(d_\gamma+2)}{\lambda_\gamma}
(\ang{h;k'}_\gamma\ang{k;h'}_\gamma+\ang{h;k}_\gamma
\ang{h';k'}_\gamma)
\]
and
\[
T_{\fa_g,\fa_g}=T_{\fb_g,\fb_g}=\ang{g;g}_\gamma\,,\qquad
T_{\fc_{h,k},\fc_{h,k}}=\ang{g;h}_\gamma\,,\qquad
T_{\fd_j,\fd_j}=\omega_j\,.
\]

\end{remark}

\section{TQFTs from the finite quantum dilogarithm}
\label{sec.6j}

An important aspect of the construction of Izumi's fusion categories 
$\calC_\gamma$ via endomorphisms of the Leavitt algebra is that 
this gives not 
only a monoidal structure and abstractly defined Hom spaces, 
but in addition 
it gives canonical bases of the Hom spaces (finite dimensional 
vector spaces
over $\BC$). This converts the tensor-valued $6j$-symbols (the 
building
blocks of the TV TQFT) into complex-valued $6j$-symbols, and 
the partition
function of the TV TQFT into a state-sum exactly as introduced 
by Turaev--Viro.

Thus, our task in this section is to compute the $6j$ symbols 
of $\calC_\gamma$.
To do so, we follow the excellent exposition of Barrett--Westbury
~\cite{Barrett-Westbury}.

\subsection{Izumi's fusion categories and Leavitt algebras}

We begin by a very brief discussion of how to define Izumi's 
fusion categories
from a collection of endomorphims of the Leavitt algebra. This is 
a subject
that has been treated comprehensively in ~\cite[Sec.5.2]{RW} 
following prior 
work of Evans--Gannon~\cite{EG17}. We need this discussion since 
we are 
using explicit formulas for the endomorphisms of the Leavitt 
algebra that 
lead to explicit bases of Homs and consequently to numerical 
values of 
the $6j$-symbols.

Fix a hyperbolic element $\gamma \in \SL_2(\BZ)$ and let 
$G=G_\gamma$ denote the finite abelian group of order $N$. 

Recall that the Leavitt algebra $L_{2N}$~\cite{Lea62} is the 
quotient of 
free associative algebra over $\BC$ with generators 
$\{S_g\}_{g\in G}$, $\{S_g'\}_{g\in G}$, 
$\{T_g\}_{g\in G}$, $\{T_g'\}_{g\in G}$, 
by the two-sided ideal generated by the relations
\begin{equation} 
\label{eq:leavitt2Na}
	S_g'S_h = \delta_{g,h}1\,,\qquad
	T_g'T_h = \delta_{g,h}1\,,\qquad
	S_g'T_h = T_h'S_g=0\, \qquad 
\sum_{g\in G}S_gS_g' + \sum_{g\in G}T_gT_g' = 1 \,.
\end{equation}
The category $\calC_\gamma$ has :

\noindent
$\bullet$
Set of simple objects $I=\{[g] \,| \, g \in G\} \cup \{[\rho]\}$ 
corresponding to the unital endomorphisms $\alpha_g, \rho$ 
of $L_{2N}$ defined below.

\noindent
$\bullet$
Homomorphisms: if $\alpha,\beta \in I$, 
\begin{equation}
    \Hom(\alpha,\beta)=\{x \in L_{2N} \, \, | \, \, x 
    \alpha(y)=\beta(y)x 
    \,\, \forall y \in L_{2N} \}
\end{equation}

\noindent
$\bullet$
Composition of homomorphisms is multiplication in $L_{2N}$:
\begin{equation}
    \Hom(\alpha,\beta) \times \Hom(\beta,\gamma) \to \Hom(\alpha,
    \gamma),
    \qquad (x,x') \mapsto x'x \,.
\end{equation}

\noindent
$\bullet$
Tensor product of simple objects is composition of endomorphisms of $L_{2N}$:
for $\alpha, \beta \in I$,
\begin{equation}
\alpha, \beta \in \End(L_{2N}), \qquad 
\alpha \otimes \beta = \alpha \circ \beta \in \End(L_{2N}) \,.
\end{equation}
\noindent
$\bullet$
Tensor product of homs:
\begin{equation}
r \in \Hom(\beta,\gamma), \,\, s \in \Hom(\sigma,\tau), \qquad
r \otimes s = r \beta(s)=\gamma(s) r \in 
\Hom(\beta \otimes \sigma, \gamma \otimes \tau) \,.
\end{equation}
The fursion rules are given in Equation~\eqref{frules}. 
We now define the endomorphisms $\alpha_g$ and $\rho$ of $L_{2N}$. 
The endomorphism $\alpha_{g}$ is given by:
\begin{equation}
\label{eq:alphag}
	\begin{aligned}
	\alpha_g(S_h) &= S_{g+h}\,, &&\alpha_g(T_h) = \ang{g;h}^{-1}T_{h}\,,
\\ 
	\alpha_g(S_h') &= S_{g+h}'\,, &&\alpha_g(T_h') = \ang{g;h}T_{h}'\,.
	\end{aligned}
\end{equation}

To define the endomorphism $\rho$, we need the auxilary 
representation
\begin{equation}
U\colon G\to L_{2N}^{\times}, \qquad
U_g = \sum_{h}\ang{g;h}^{-1}S_hS_h' + \sum_{h}T_{h-g}T_h'\,.
\end{equation}

The endomorphism associated to $\rho$ is given as follows:
\begin{align*}
\rho(S_g) &= U_g\Big(\frac{1}{d}\sum_{h}S_h + \frac{1}{\sqrt{d}}\sum_{h}
\ang{h}T_hT_{-h}\Big)U_{-g} = \Big(\frac{1}{d}
\sum_{h}\ang{g;h}^{-1}S_h + \frac{1}
{\sqrt{d}}\sum_{h}\ang{h}T_{h-g}T_{-h}\Big)U_{-g}\,,
\\ 
\rho(S_g') &= U_g\Big(\frac{1}{d}\sum_{h}S_h' + \frac{1}
{\sqrt{d}}\sum_{h}\ang{h}^{-1}
T_{-h}'T_{h}'\Big)U_{-g} = U_g\Big(\frac{1}{d}\sum_{h}\ang{g;h}S_h' + 
\frac{1}
{\sqrt{d}}\sum_{h}\ang{h}^{-1}T_{-h}'T_{h-g}'\Big)
\\
\rho(T_g) &=\frac{\lambda}{\sqrt{d}\sqrt{N}}\sum_{h,k\in G}\ang{h-g;k}
S_hT_k'
+\frac{\lambda^{-1}}{\ang{g}\sqrt{N}}\sum_{h,k\in G}\ang{g+h;k}^{-1}
T_kS_hS_h'
\\
&\qquad+\frac{1}{\sqrt{N}}\sum_{h,k\in G}\ang{h}E(h+g)\ang{g;k}^{-1}
T_{h+k}T_{-h}T_k'
\\
\rho(T_g') &=\frac{\lambda^{-1}}{\sqrt{d}\sqrt{N}}\sum_{h,k\in G}
\ang{h-g;k}^{-1}T_kS_h'+\frac{\lambda \ang{g}}{\sqrt{N}}\sum_{h,k\in 
G}\ang{g+h;k}S_hS_h'T_k'
\\
&\qquad+\frac{1}{\sqrt{N}}\sum_{h,k\in G}\frac{\ang{g+h}}{\ang{h}}E(-h-
g)\ang{g;k}T_kT_{-h}'T_{h+k}'
\end{align*}
where $\lambda$ and $E$ are as in Section~\ref{sub.finiteQDL}.

\subsection{Bases for hom spaces}

The endomorphisms corresponding to the simple objects of 
$\calC_\gamma$ satisfy the following properties

\begin{equation}
\label{rho2}
\begin{aligned}
\alpha_g \alpha_h   &= \alpha_{g+h}, \qquad
\alpha_g \rho = \rho, \qquad
\rho \alpha_g (x) = U_g \rho(x) U_{-g}, \qquad
\alpha_g U_h = \langle g;h \rangle U_h \\
\rho^2(x) & = \sum_{g\in G}S_g\alpha_g(x)S_g' + \sum_{h\in G}
T_h\rho(x)T_h'
\end{aligned}
\end{equation}
for all $g,h \in G$ and $x \in L_{2N}$. The identities in the first 
line of~\eqref{rho2} are elementary. The last one however, 
is equivalent to the pentagon relation for the finite quantum dilogarithm.

It follows from Equation~\eqref{rho2} and that the non-zero 
hom-spaces are based as follows:

\begin{equation}
\label{homs}
\begin{aligned}
\mathrm{Hom}([g+h],[g]\otimes[h])&=\CC 1\,,&\qquad&
\mathrm{Hom}([\rho],[g]\otimes[\rho])&\!\!\!\!=\;&\CC 1\,,\\
\mathrm{Hom}([\rho],[\rho]\otimes[g])&=\CC U_g\,,&\qquad&
\mathrm{Hom}([g],[\rho]\otimes[\rho])&\!\!\!\!=\;&\CC S_g\,,\\
\mathrm{Hom}([\rho],[\rho]\otimes[\rho])&=\oplus_{g\in G}\CC T_g\,.
\end{aligned}
\end{equation}
For example,
\begin{equation}
\begin{aligned}
x \in \Hom([\rho],[\rho] \otimes [g]) & \iff x \rho(y) = 
[\rho] \otimes [g](y)x
\iff x \rho(y) = \rho(\alpha_g(y))x \\ &
\iff x \rho(y) = U_g \rho(y) U_{-g} x
\end{aligned}
\end{equation}
for all $y \in L_{2N}$, which is equivalent to $U_{-g}x$ commutes
with $\rho(y)$ for all $y \in L_{2N}$. But by a property of $\rho$, 
it follows that $U_{-g}x \in \BC 1$, so $x \in \BC U_g$.

We denote the canonical bases in the hom spaces of~\eqref{homs} 
by $\bullet$ 
for the first four cases of (where at most two $\rho$ appear), 
and by $T_g$, $g \in G$ for the last case. 

Equations~\eqref{homs} are complemented with the identification
$\Hom(A,B) = \Hom(B,A)^*$, given explicitly by
\begin{equation}
\label{homs2}
\begin{aligned}
\mathrm{Hom}([g]\otimes[h],[g+h])&=\CC 1\,,&\qquad&
\mathrm{Hom}([g]\otimes[\rho],[\rho])&\!\!\!\!=\;&\CC 1\,,\\
\mathrm{Hom}([\rho]\otimes[g],[\rho])&=\CC U_{-g}\,,&\qquad&
\mathrm{Hom}([\rho]\otimes[\rho],[g])&\!\!\!\!=\;&\CC S_g'\,,\\
\mathrm{Hom}([\rho]\otimes[\rho],[\rho])&=\oplus_{g\in G}\CC T_g'\,.
\end{aligned}
\end{equation}

We define the dual (both left and right) of $[g]$ to be $[-g]$ 
and set both 
evaluation and co-evaluation morphisms to be 
$e_{\alpha_g}=b_{\alpha_g}=1$ (since 
the relevant hom-space is $\mathrm{Hom}(1,1)=\CC 1$).
The object $\rho$ is self-dual, and since $S_0\in\Hom(1,\rho^2)$ 
and $S_0'\in 
\Hom(\rho^2,1)$, the evaluation and coevaluation are given by 
$e_{\rho}=dS_0'$, $b_{\rho}=S_0$.

\subsection{\texorpdfstring{Thetas and $6j$-symbols}{Thetas and 6j-symbols}}

The quantum dimension of the simple objects (defined as the left 
or right trace of the identity endomorphism) is given by 
\begin{equation}
\dimq([g])=1, \qquad \dimq([\rho])=d, 
\end{equation}
where $d$ is as in~\eqref{Gconstants}. This follows from the fact 
that $e_\rho b_\rho=d$.

We will now compute the standard isomorphisms between 
$\mathrm{Hom}([\rho],[\rho]\otimes[\rho])$ and itself. 
This comes from the diagram

\begin{figure}[htpb!]
\begin{center}
\begin{tikzpicture}[scale=0.8]
\draw[->](1,4)--(1,3);
\draw(1,3.2) node[right] {$S_0$};
\draw[->](3,1)--(3,0);
\draw(3,0.8) node[right] {$dS_0'$};
\draw[->](4,4) to[out angle=-90, in angle=90,curve through = {(4,2) 
(3,1) (2,2) (1,1)}] (1,0);
\draw(1.9,2.1) node[right] {$T_g$};
\draw[->](2,2) to[out angle=90, in angle=90,curve through = {(1,3)}] 
(0,0);
\end{tikzpicture}   
\end{center}
\caption{$\theta$-tensors.}
\label{fig:theta}
\end{figure}
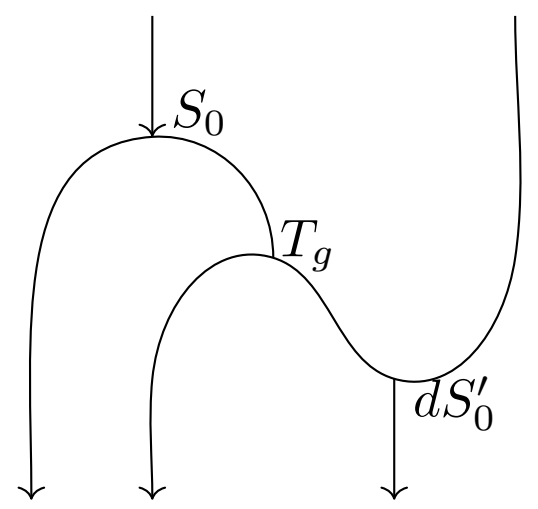

We find the isomorphism sends $T_g$ to 
\[
(1\otimes 1\otimes dS_0')(1\otimes T_g\otimes 1)(S_0\otimes 1)
=
d\rho^2(S_0')\rho(T_g)S_0
=
\frac{\lambda^{-1}}{\ang{g}\sqrt{N}}
\sum_{k\in G}\ang{g;k}^{-1}T_k\,.
\]
This isomorphism has order $3$.
The fact that these isomorphisms are not trivial means that the 
symmetries of the $6j$-symbols are not simple relabellings.

We now come to the main building block of the TV TQFT, namely the 
$6j$-symbols.
We will follow the exposition of 
Barrett--Westbury~\cite{Barrett-Westbury}.
We use ordered tetrahedra, which are either positive or negative, 
and we 
label their vertices, edges and faces as follows.

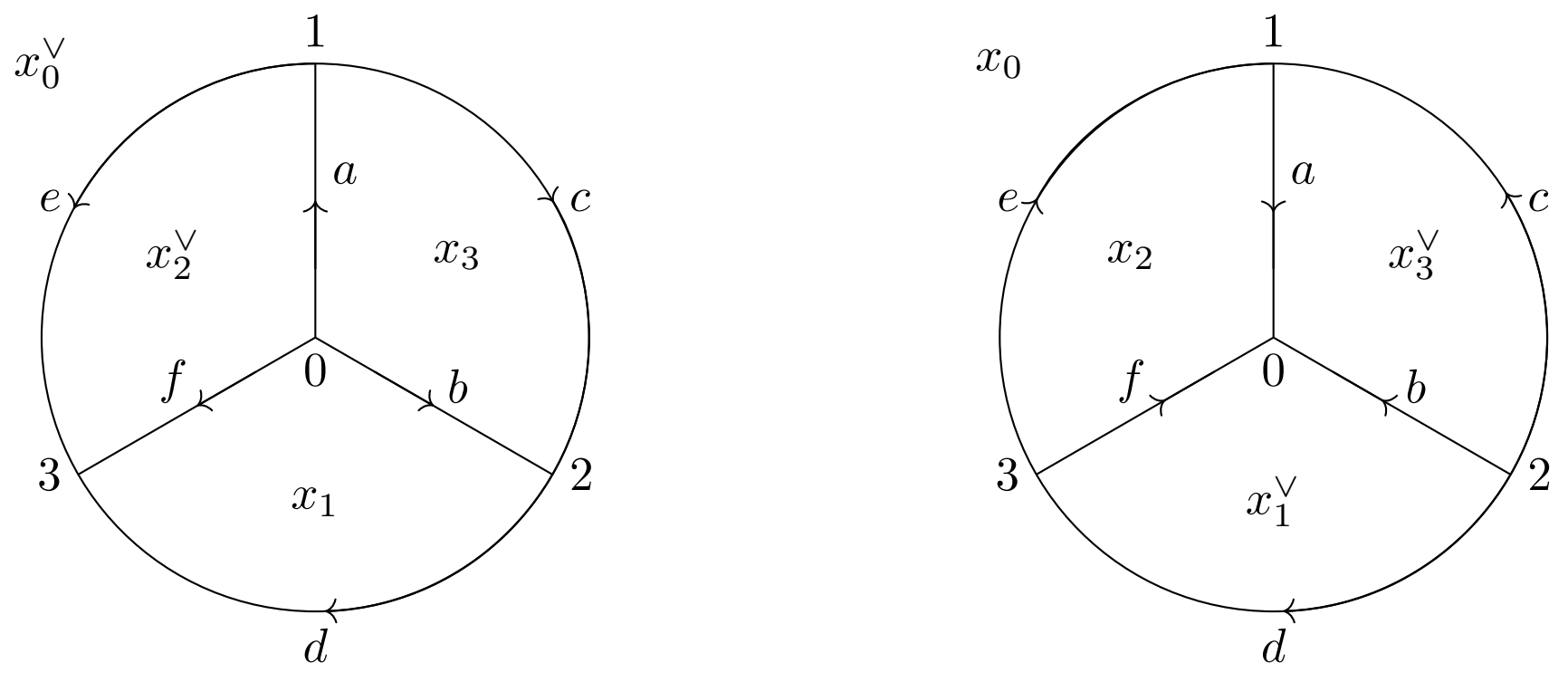
\begin{figure}[htpb!]
\begin{center}
\begin{tikzpicture}
\draw(0,0) circle (2cm);
\draw[->](0,0)--(0,1);
\draw(0,0.5)--(0,2);
\draw[->](0,0)--(-30:1);
\draw(-30:0.5)--(-30:2);
\draw[->](0,0)--(-150:1);
\draw(-150:0.5)--(-150:2);
\draw[-<](2,0) arc (0:2)--(30:2);
\draw[->](0,2) arc (90:2)--(150:2);
\draw[->](2,0) arc (0:2)--(-90:2);
\draw(0,0) node[below] {$0$};
\draw(0,2) node[above] {$1$};
\draw(-30:2) node[right] {$2$};
\draw(-150:2) node[left] {$3$};
\draw(0,1.2) node[right] {$a$};
\draw(-30:1.2) node[above] {$b$};
\draw(30:2) node[right] {$c$};
\draw(0,-2) node[below] {$d$};
\draw(150:2) node[left] {$e$};
\draw(-150:1.2) node[above] {$f$};
\draw(-2,2) node {$x_0^\vee$};
\draw(-90:1.2) node {$x_1$};
\draw(150:1.2) node {$x_2^\vee$};
\draw(30:1.2) node {$x_3$};

\draw[xshift=7cm](0,0) circle (2cm);
\draw[xshift=7cm,-<](0,0)--(0,1);
\draw[xshift=7cm](0,0.5)--(0,2);
\draw[xshift=7cm,-<](0,0)--(-30:1);
\draw[xshift=7cm](-30:0.5)--(-30:2);
\draw[xshift=7cm,-<](0,0)--(-150:1);
\draw[xshift=7cm](-150:0.5)--(-150:2);
\draw[xshift=7cm,->](2,0) arc (0:2)--(30:2);
\draw[xshift=7cm,-<](0,2) arc (90:2)--(150:2);
\draw[xshift=7cm,->](2,0) arc (0:2)--(-90:2);
\draw[xshift=7cm](0,0) node[below] {$0$};
\draw[xshift=7cm](0,2) node[above] {$1$};
\draw[xshift=7cm](-30:2) node[right] {$2$};
\draw[xshift=7cm](-150:2) node[left] {$3$};
\draw[xshift=7cm](0,1.2) node[right] {$a$};
\draw[xshift=7cm](-30:1.2) node[above] {$b$};
\draw[xshift=7cm](30:2) node[right] {$c$};
\draw[xshift=7cm](0,-2) node[below] {$d$};
\draw[xshift=7cm](150:2) node[left] {$e$};
\draw[xshift=7cm](-150:1.2) node[above] {$f$};
\draw[xshift=7cm](-2,2) node {$x_0$};
\draw[xshift=7cm](-90:1.2) node {$x_1^\vee$};
\draw[xshift=7cm](30:1.2) node {$x_3^\vee$};
\draw[xshift=7cm](150:1.2) node {$x_2$};
\end{tikzpicture}
\end{center}
\caption{Positive (on the left) and negative 
(on the right) tetrahedron.}
\label{fig:2tet}
\end{figure}

Here, $a,b,c,d,e,f$ are simple objects labelling the edges of a 
tetrahedron and 
\begin{equation}
\label{abcd2}
x_0\in\mathrm{Hom}(d\otimes c,e)\,,\quad
x_1\in\mathrm{Hom}(f,d\otimes b)\,,\quad
x_2\in\mathrm{Hom}(e\otimes a,f)\,,\quad
x_3\in\mathrm{Hom}(b,c\otimes a)\,,
\end{equation}
labelling the faces of a tetrahedron. Using the linear 
functional 
\begin{equation}
\label{delta+}
x_0 \otimes x_1 \otimes x_2 \otimes x_3 \mapsto 
\mathrm{tr}(x_2(x_0\otimes\mathrm{id})(\mathrm{id}\otimes x_3)x_1)
\end{equation}
on $
\mathrm{Hom}(d\otimes c,e) \otimes \mathrm{Hom}(f,d\otimes b)
\otimes \mathrm{Hom}(e\otimes a,f) \otimes \mathrm{Hom}(b,c\otimes a) 
\to \BC 
$,
we define~\cite[Sec.3]{Barrett-Westbury}
\[
\Big\{\begin{array}{ccc|cc}
a & b & c & x_0 & x_2\\
d & e & f & x_1 & x_3
\end{array}\Big\}_+
=
\mathrm{tr}(x_2(x_0\otimes\mathrm{id})(\mathrm{id}\otimes x_3)x_1)
=
\mathrm{ev}((x_2(x_0\otimes\mathrm{id})(\mathrm{id}\otimes 
x_3)x_1)\otimes\mathrm{id})\mathrm{coev}\,.
\]

Using the non-degeneracy of $\tr$ as a bilinear pairing between hom-
spaces, this
defines a map~\cite[Defn.1.5]{Barrett-Westbury}
\[
\Big\{\begin{array}{ccc}
    a & b & c\\
    d & e & f
\end{array}\Big\}_+
:
\mathrm{Hom}(f,d\otimes b)\otimes\mathrm{Hom}(b,c\otimes a)
\to
\mathrm{Hom}(e,d\otimes c)\otimes
\mathrm{Hom}(f,e\otimes a)\,.
\]

The morphism can be pictured as follows and, upon taking the trace, 
this recovers the dual graph of the tetrahedron (which happens to be 
again a tetraherdron):

\begin{figure}[htpb!]
\begin{center}
\begin{tikzpicture}[scale=0.6]
\draw[<-](0,4)--(0,5) node[right] {$f$};
\draw(0,4)--(-1,3);
\draw(0.5,3.5) node[right] {$b$};
\draw[->](0,4)--(1,3);
\draw[->](-1,3) to[out angle=-135, in angle=135,curve through = 
{(-1.5,2)}] (-1,1);
\draw(-1.5,2) node[left] {$d$};
\draw(1,3)--(0,2) node[right] {$c$};
\draw[->](1,3) to[out angle=-45, in angle=45,curve through = {(1.5,2)}] 
(0,0);
\draw(1.5,2) node[right] {$a$};
\draw[->](0,2)--(-1,1);
\draw(-0.5,0.5) node[right] {$e$};
\draw[->](-1,1)--(0,0);
\draw[->](0,0)--(0,-1) node[right] {$f$};
\filldraw(0,4) node[right] {$x_1$};
\filldraw(1,3)  node[right] {$x_3$};
\filldraw(-1,1) node[right] {$x_0$};
\filldraw(0,0) node[right] {$x_2$};
\end{tikzpicture}    
\end{center}
\caption{Composition of maps defining the quantum $6j$-symbols 
of a positive tetrahedron.}
\label{fig:tet6j}
\end{figure}
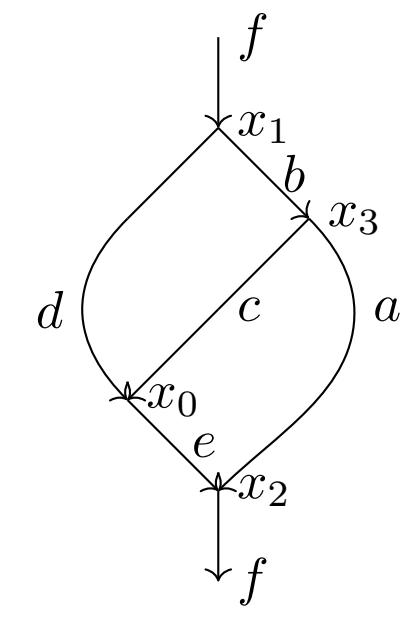

Likewise, for a negative tetrahedron, we use the linear 
functional
\begin{equation}
\label{delta-}
x_0 \otimes x_1 \otimes x_2 \otimes x_3 \mapsto 
\mathrm{tr}(x_1(\mathrm{id}\otimes x_3)(x_0 \otimes \mathrm{id})x_2)
\end{equation}
on $
\mathrm{Hom}(e, d \otimes c) \otimes \mathrm{Hom}(d\otimes b,f)
\otimes \mathrm{Hom}(f,e\otimes a) \otimes \mathrm{Hom}(c\otimes a,b) 
\to \BC 
$,
we define~\cite[Sec.3]{Barrett-Westbury}
\[
\Big\{\begin{array}{ccc|cc}
a & b & c & x_0 & x_2\\
d & e & f & x_1 & x_3
\end{array}\Big\}_-
=
\mathrm{tr}(x_1(\mathrm{id}\otimes x_3)(x_0 \otimes \mathrm{id})x_2)
=
\mathrm{ev}(x_1(\mathrm{id}\otimes x_3)(x_0 \otimes \mathrm{id})x_2)\otimes\mathrm{id})\mathrm{coev}\,.
\]

Using the non-degeneracy of $\tr$ as a bilinear pairing between hom-spaces, this
defines a map~\cite[Defn.1.5]{Barrett-Westbury}
\[
\Big\{\begin{array}{ccc}
    a & b & c\\
    d & e & f
\end{array}\Big\}_-
:
\mathrm{Hom}(e,d\otimes c)\otimes\mathrm{Hom}(f,e\otimes a)
\to
\mathrm{Hom}(f,d\otimes b)\otimes
\mathrm{Hom}(b,c\otimes a)\,.
\]

\begin{proposition}
\label{prop.6j}
The $6j$-symbols of a positive tetrahedron are given by:
\begin{equation}
\label{6jval:6}
\begin{aligned}
\Big\{\begin{array}{ccc|cc}
    \rho & \rho & \rho & g & h\\
    \rho & \rho & \rho & k & j
\end{array}\Big\}_+
&=\frac{\ang{h}}{d\sqrt{N}\ang{j;k}}E(j-h)\delta_{g,k-h}\,,
\end{aligned}
\end{equation}

\begin{equation}
\label{6jval:5}
\begin{array}{rlrl}
\Big\{\begin{array}{ccc|cc}
    a & \rho & \rho & g & \bullet\\
    \rho & \rho & \rho & k & \bullet
\end{array}\Big\}_+
&\!\!\!\!=\;\frac{\delta_{g,k+a}\ang{k+a}}{d\ang{k}}, 
\qquad
&\Big\{\begin{array}{ccc|cc}
    \rho & b & \rho & g & h\\
    \rho & \rho & \rho & \bullet & \bullet
\end{array}\Big\}_+
&\!\!\!\!=\;d^{-3/2}\delta_{g+h+b}\ang{h},
\\[0.5cm]
\Big\{\begin{array}{ccc|cc}
    \rho & \rho & c & \bullet & h\\
    \rho & \rho & \rho & k & \bullet
\end{array}\Big\}_+
&\!\!\!\!=\;d^{-1}\delta_{h,k+c},
\qquad
&\Big\{\begin{array}{ccc|cc}
    \rho & \rho & \rho & \bullet & h\\
    d & \rho & \rho & \bullet & j
\end{array}\Big\}_+
&\!\!\!\!=\;\frac{\delta_{h,j}}{d\ang{d;j}}\\[0.5cm]
\Big\{\begin{array}{ccc|cc}
    \rho & \rho & \rho & \bullet & \bullet\\
    \rho & e & \rho & k & j
\end{array}\Big\}_+
&\!\!\!\!=\;\frac{\lambda\ang{e-j;k}}{d^{3/2}\sqrt{N}},
\qquad
&\Big\{\begin{array}{ccc|cc}
    \rho & \rho & \rho & g & \bullet\\
    \rho & \rho & f & \bullet & j
\end{array}\Big\}_+
&\!\!\!\!=\;\frac{\lambda^{-1}\ang{j+f;g}}{\ang{j}\sqrt{N}}
\end{array}
\end{equation}

\begin{equation}
\label{6jval:4}
\Big\{\begin{array}{ccc|cc}
    a & \rho & \rho & \bullet & \bullet\\
    d & \rho & \rho & \bullet & \bullet
\end{array}\Big\}_+
=d^{-1}\ang{a;d},
\quad
\Big\{\begin{array}{ccc|cc}
    \rho & b & \rho & \bullet & \bullet\\
    \rho & e & \rho & \bullet & \bullet
\end{array}\Big\}_+
=\frac{1}{d^2\ang{b;e}},\quad
\Big\{\begin{array}{ccc|cc}
    \rho & \rho & c & \bullet & \bullet\\
    \rho & \rho & f & \bullet & \bullet
\end{array}\Big\}_+
=\ang{c;f}
\end{equation}

\begin{equation}
\label{6jval:3}
\begin{aligned}
\Big\{\begin{array}{ccc|cc}
    a & b & c & \bullet & \bullet\\
    \rho & \rho & \rho & \bullet & \bullet
\end{array}\Big\}_+
&=d^{-1}\delta_{b,a+c}, 
\qquad
\Big\{\begin{array}{ccc|cc}
    a & \rho & \rho & \bullet & \bullet\\
    \rho & e & f & \bullet & \bullet
\end{array}\Big\}_+
&\!\!\!\!=\,\delta_{a+e,f}, \\
\Big\{\begin{array}{ccc|cc}
    \rho & b & \rho & \bullet & \bullet\\
    d & \rho & f & \bullet & \bullet
\end{array}\Big\}_+
&=\delta_{b+d,f}, 
\qquad
\Big\{\begin{array}{ccc|cc}
    \rho & \rho & c & \bullet & \bullet\\
    d & e & \rho & \bullet & \bullet
\end{array}\Big\}_+
&\!\!\!\!=\,d^{-1}\delta_{c+d,e},
\end{aligned}
\end{equation}

\begin{equation}
\label{6jval:1}
\begin{aligned}
\Big\{\begin{array}{ccc|cc}
    a & b & c & \bullet & \bullet\\
    d & e & f & \bullet & \bullet
\end{array}\Big\}_+
&=
\delta_{f,b+d}\delta_{b,a+c}\delta_{c+d,e}\delta_{a+e,f}
\end{aligned}
\end{equation}
Replacing positive tetrahedra by negative replaces the right hand 
side of the above equations by complex-conjugation of roots of unity and $E(j)\to \ang{j}E(-j)$
\end{proposition}

In $6j$-symbols above interpolate from the presence of all 
$\rho$, that take the most arithmetically interesting values, 
to the simplest case of no $\rho$ that essentially describes
a Dijkgraaf--Witten TQFT with abelian group. 

\begin{remark}
\label{rem.kinker}
Kashaev's kinematical kernel is 
\begin{equation}
K^{x_0,x_2,x_4}_{x_1,x_3} = 
\langle x_0;x_4 \rangle \delta(x_0-x_1+x_2) \delta(x_2-x_3+x_4)
\end{equation}
where $x_i$ are variables in the $i$-th face of the positive 
tetrahedron and $x_4$ is a variable in the tetrahedron. 
It is a distribution that satisfies
the ordered pentagon identity for each locally compact 
abelian group~\cite{Kas4d}. 
The first $6j$-symbol in~\eqref{6jval:6} is of this form with  
$g=x_0$, $k=x_1$, $h=x_2$, $j=x_3$ and $z=x_3-x_2$:
\begin{equation}
\label{6jval2}
\begin{aligned}
\Big\{\begin{array}{ccc|cc}
    \rho & \rho & \rho & x_0 & x_2 \\
    \rho & \rho & \rho & x_1 & x_3
\end{array}\Big\}_+
&=\frac{d}{\sqrt{N}}E(x_4) K^{x_0,x_2,x_4}_{x_1,x_3} \,.
\end{aligned}
\end{equation}
The 5-tensor $K^{x_0,x_2,x_4}_{x_1,x_3}$ satisfies a system of 
cubic equations related to the (3,3)-Pachner moves in dimension 4, 
given explicitly in Eqn.(3.10) of~\cite{Kas4d} (see also Eqn.(3.22)). 
Here, the interpretation is that $x_i$ for $i=0,\dots,4$ is the 
tetrahedron opposite to face $i$ of the pentachoron 
(4-dimensional simplex).
\end{remark}

\begin{remark}
In Equation~\eqref{6jval:1}, one can remove any one
of the delta functions.
\end{remark}

\begin{proof}[Proof of Proposition~\ref{prop.6j}]
Using
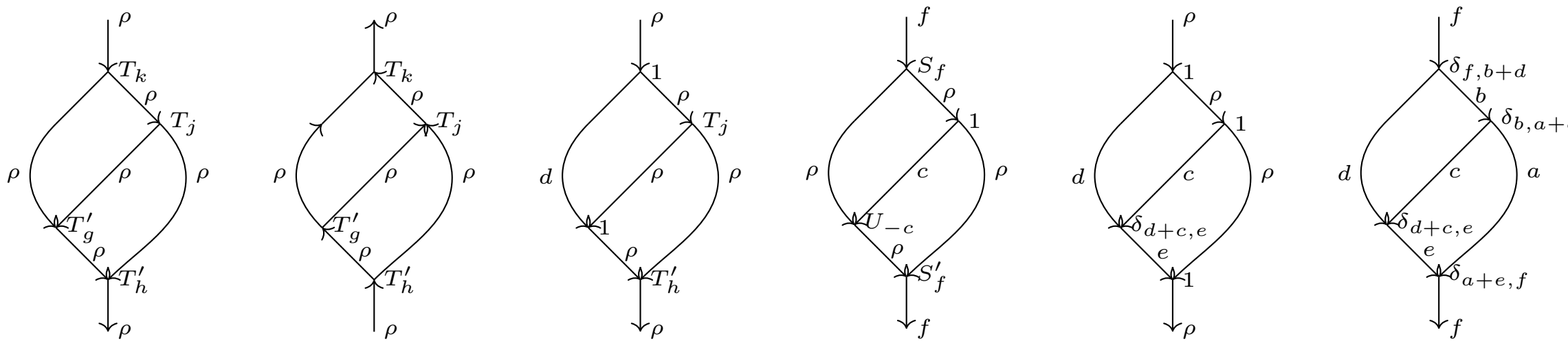
\begin{tiny}
\begin{center}
\begin{tikzpicture}[scale=0.5]
\draw[<-](0,4)--(0,5) node[right] {$\rho$};
\draw(0,4)--(-1,3);
\draw(0.5,3.5) node[right] {$\rho$};
\draw[->](0,4)--(1,3);
\draw[->](-1,3) to[out angle=-135, in angle=135,curve through = 
{(-1.5,2)}] (-1,1);
\draw(-1.5,2) node[left] {$\rho$};
\draw(1,3)--(0,2) node[right] {$\rho$};
\draw[->](1,3) to[out angle=-45, in angle=45,curve through = 
{(1.5,2)}] 
(0,0);
\draw(1.5,2) node[right] {$\rho$};
\draw[->](0,2)--(-1,1);
\draw(-0.5,0.5) node[right] {$\rho$};
\draw[->](-1,1)--(0,0);
\draw[->](0,0)--(0,-1) node[right] {$\rho$};
\filldraw(0,4) node[right] {$T_k$};
\filldraw(1,3)  node[right] {$T_j$};
\filldraw(-1,1) node[right] {$T_g'$};
\filldraw(0,0) node[right] {$T_h'$};
\end{tikzpicture}
$\quad$
\begin{tikzpicture}[scale=0.5]
\draw[->](0,4)--(0,5) node[right] {$\rho$};
\draw(0,4)--(-1,3);
\draw(0.5,3.5) node[right] {$\rho$};
\draw[<-](0,4)--(1,3);
\draw[<-](-1,3) to[out angle=-135, in angle=135,curve through = 
{(-1.5,2)}] (-1,1);
\draw(-1.5,2) node[left] {$\rho$};
\draw[<-](1,3)--(0,2) node[right] {$\rho$};
\draw[<-](1,3) to[out angle=-45, in angle=45,curve through = 
{(1.5,2)}] 
(0,0);
\draw(1.5,2) node[right] {$\rho$};
\draw(0,2)--(-1,1);
\draw(-0.5,0.5) node[right] {$\rho$};
\draw[<-](-1,1)--(0,0);
\draw[<-](0,0)--(0,-1) node[right] {$\rho$};
\filldraw(0,4) node[right] {$T_k$};
\filldraw(1,3)  node[right] {$T_j$};
\filldraw(-1,1) node[right] {$T_g'$};
\filldraw(0,0) node[right] {$T_h'$};
\end{tikzpicture}
$\quad$
\begin{tikzpicture}[scale=0.5]
\draw[<-](0,4)--(0,5) node[right] {$\rho$};
\draw(0,4)--(-1,3);
\draw(0.5,3.5) node[right] {$\rho$};
\draw[->](0,4)--(1,3);
\draw[->](-1,3) to[out angle=-135, in angle=135,curve through = 
{(-1.5,2)}] (-1,1);
\draw(-1.5,2) node[left] {$d$};
\draw(1,3)--(0,2) node[right] {$\rho$};
\draw[->](1,3) to[out angle=-45, in angle=45,curve through = 
{(1.5,2)}] 
(0,0);
\draw(1.5,2) node[right] {$\rho$};
\draw[->](0,2)--(-1,1);
\draw(-0.5,0.5) node[right] {$\rho$};
\draw[->](-1,1)--(0,0);
\draw[->](0,0)--(0,-1) node[right] {$\rho$};
\filldraw(0,4) node[right] {$1$};
\filldraw(1,3)  node[right] {$T_j$};
\filldraw(-1,1) node[right] {$1$};
\filldraw(0,0) node[right] {$T_h'$};
\end{tikzpicture}
$\quad$
\begin{tikzpicture}[scale=0.5]
\draw[<-](0,4)--(0,5) node[right] {$f$};
\draw(0,4)--(-1,3);
\draw(0.5,3.5) node[right] {$\rho$};
\draw[->](0,4)--(1,3);
\draw[->](-1,3) to[out angle=-135, in angle=135,curve through = 
{(-1.5,2)}] (-1,1);
\draw(-1.5,2) node[left] {$\rho$};
\draw(1,3)--(0,2) node[right] {$c$};
\draw[->](1,3) to[out angle=-45, in angle=45,curve through = 
{(1.5,2)}] 
(0,0);
\draw(1.5,2) node[right] {$\rho$};
\draw[->](0,2)--(-1,1);
\draw(-0.5,0.5) node[right] {$\rho$};
\draw[->](-1,1)--(0,0);
\draw[->](0,0)--(0,-1) node[right] {$f$};
\filldraw(0,4) node[right] {$S_f$};
\filldraw(1,3)  node[right] {$1$};
\filldraw(-1,1) node[right] {$U_{-c}$};
\filldraw(0,0) node[right] {$S_f'$};
\end{tikzpicture}
$\quad$
\begin{tikzpicture}[scale=0.5]
\draw[<-](0,4)--(0,5) node[right] {$\rho$};
\draw(0,4)--(-1,3);
\draw(0.5,3.5) node[right] {$\rho$};
\draw[->](0,4)--(1,3);
\draw[->](-1,3) to[out angle=-135, in angle=135,curve through = 
{(-1.5,2)}] (-1,1);
\draw(-1.5,2) node[left] {$d$};
\draw(1,3)--(0,2) node[right] {$c$};
\draw[->](1,3) to[out angle=-45, in angle=45,curve through = 
{(1.5,2)}] 
(0,0);
\draw(1.5,2) node[right] {$\rho$};
\draw[->](0,2)--(-1,1);
\draw(-0.5,0.5) node[right] {$e$};
\draw[->](-1,1)--(0,0);
\draw[->](0,0)--(0,-1) node[right] {$\rho$};
\filldraw(0,4) node[right] {$1$};
\filldraw(1,3)  node[right] {$1$};
\filldraw(-1,1) node[right] {$\delta_{d+c,e}$};
\filldraw(0,0) node[right] {$1$};
\end{tikzpicture}
$\quad$
\begin{tikzpicture}[scale=0.5]
\draw[<-](0,4)--(0,5) node[right] {$f$};
\draw(0,4)--(-1,3);
\draw(0.5,3.5) node[right] {$b$};
\draw[->](0,4)--(1,3);
\draw[->](-1,3) to[out angle=-135, in angle=135,curve through = 
{(-1.5,2)}] (-1,1);
\draw(-1.5,2) node[left] {$d$};
\draw(1,3)--(0,2) node[right] {$c$};
\draw[->](1,3) to[out angle=-45, in angle=45,curve through = 
{(1.5,2)}] 
(0,0);
\draw(1.5,2) node[right] {$a$};
\draw[->](0,2)--(-1,1);
\draw(-0.5,0.5) node[right] {$e$};
\draw[->](-1,1)--(0,0);
\draw[->](0,0)--(0,-1) node[right] {$f$};
\filldraw(0,4) node[right] {$\delta_{f,b+d}$};
\filldraw(1,3)  node[right] {$\delta_{b,a+c}$};
\filldraw(-1,1) node[right] {$\delta_{d+c,e}$};
\filldraw(0,0) node[right] {$\delta_{a+e,f}$};
\end{tikzpicture} 
\end{center}
\end{tiny}
and going though each of the six pictures, we compute:
\[
\begin{aligned}
\Big\{\begin{array}{ccc|cc}
    \rho & \rho & \rho & g & h\\
    \rho & \rho & \rho & k & j
\end{array}\Big\}_+
&\!\!=
\frac{S_0'}{d}T_h'(T_g'\otimes 1)(1\otimes T_j)T_kS_0
=
\frac{S_0'}{d}T_h'T_g'\rho(T_j)T_kS_0
=
\frac{\ang{h}E(j-h)}{d\sqrt{N}\ang{j;k}}\delta_{g,k-h}\,,\\
\Big\{\begin{array}{ccc|cc}
    \rho & \rho & \rho & g & h\\
    \rho & \rho & \rho & k & j
\end{array}\Big\}_-
&\!\!=
\frac{S_0'}{d}T_k'(1\otimes T_j')(T_g\otimes 1)T_hS_0
=
\frac{S_0'}{d}T_k'\rho(T_j')T_gT_hS_0
=
\frac{\ang{j-h}E(h-j)}{d\sqrt{N}\ang{h}\ang{j;k}^{-1}}\delta_{j,h-k}\,,\\
\Big\{\begin{array}{ccc|cc}
    \rho & \rho & \rho & \bullet & h\\
    d & \rho & \rho & \bullet & j
\end{array}\Big\}_+
&\!\!=
\frac{S_0'}{d}T_h'(1\otimes 1)(1\otimes T_j)1
S_0
=
\frac{S_0'}{d}T_h'\alpha_d(T_j)S_0
=\frac{\delta_{h,j+d}}{d\ang{d;j}}\,,\\
\Big\{\begin{array}{ccc|cc}
    \rho & \rho & c & \bullet & \bullet\\
    \rho & \rho & f & \bullet & \bullet
\end{array}\Big\}_+
&\!\!=
1S_f'(U_{-c}\otimes 1)(1\otimes 1)S_f1
=
S_f'U_{-c}S_f
=
\ang{c;f}\,,\\
\Big\{\begin{array}{ccc|cc}
    \rho & \rho & c & \bullet & \bullet\\
    d & e & \rho & \bullet & \bullet
\end{array}\Big\}_+
&\!\!=
\delta_{d+c,e}\frac{S_0'}{d}1(1\otimes 1)(1\otimes 1)1S_0
=
d^{-1}\delta_{d+c,e}\,,\\
\Big\{\begin{array}{ccc|cc}
    a & b & c & \bullet & \bullet\\
    d & e & f & \bullet & \bullet
\end{array}\Big\}_+
&\!\!=
\delta_{f,b+d}\delta_{b,a+c}\delta_{c+d,e}\delta_{a+e,f}\,.
\end{aligned}
\]
\end{proof}

The $6j$-symbols satisfy an inversion and a pentagon identity; 
see~\cite[Prop.5.2]{Barrett-Westbury}. We now relate their 
pentagon identity with
the one of the finite quantum dilogarithm~\cite{RW}.

\begin{lemma}
The pentagon identity for the $6j$-symbols where all external 
edges are 
labelled by $\rho$ is equivalent to the pentagon identity of 
the finite 
quantum dilogarithm.
\end{lemma}

\begin{figure}[htpb!]
\begin{center}
\begin{tikzpicture}[scale=0.6]
\draw(0.4,-0.2)--(0,2.2) node[above] {$3$};
\draw(0,2.2)--(2,0.5) node[right] {$2$};
\draw(0,2.2)--(-2,0.5) node[left] {$0$};
\draw(0.4,-0.2)--(-2,0.5);
\draw(0.4,-0.2)--(2,0.5);
\draw[dashed](-2,0.5)--(2,0.5);
\draw(2,0.5)--(0,-2) node[below] {$1$};
\draw(0,-2)--(-2,0.5);
\draw(0,-2)--(0.4,-0.2);
\draw[<->](3.2,0.5)--(4.8,0.5);
\draw(0.6,-0.4) node {$4$};
\draw[xshift=8cm](0.4,-0.2)--(0,2.2) node[above] {$3$};
\draw[xshift=8cm](0,2.2)--(2,0.5) node[right] {$2$};
\draw[xshift=8cm](0,2.2)--(-2,0.5) node[left] {$0$};
\draw[xshift=8cm](0.4,-0.2)--(-2,0.5);
\draw[xshift=8cm](0.4,-0.2)--(2,0.5);
\draw[xshift=8cm,dashed](-2,0.5)--(2,0.5);
\draw[xshift=8cm](2,0.5)--(0,-2) node[below] {$1$};
\draw[xshift=8cm](0,-2)--(-2,0.5);
\draw[xshift=8cm](0,-2)--(0.4,-0.2);
\draw[xshift=8cm,dashed] (0,2.2)--(0,-2);
\draw[xshift=8cm](0.6,-0.4) node {$4$};
\end{tikzpicture}
\end{center}
\caption{The 2--3 Pachner move.}
\label{fig:5tet}
\end{figure}
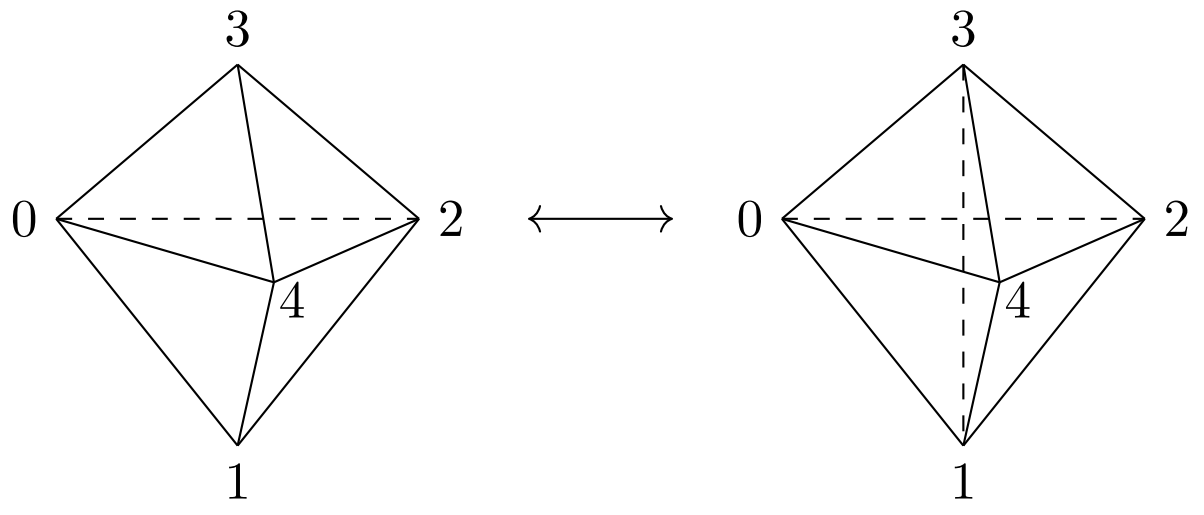

\begin{proof}
The pentagon relation for the 5 ordered tetrahedra of the 
above drawn 5-simplex where
all edges labelled by $\rho$ states that
\[
\begin{aligned}
&\sum_{c,b'}
\delta_{c,b'}
\Big\{\begin{array}{ccc|cc}
\rho & \rho & \rho & a & c\\
\rho & \rho & \rho & b & d
\end{array}\Big\}_+
\Big\{\begin{array}{ccc|cc}
\rho & \rho & \rho & a' & c'\\
\rho & \rho & \rho & b' & d'
\end{array}\Big\}_+ \\
&=
d \sum_{x_0,x_3,x_1',x_3',x_0'',x_2''}
\delta_{x_0,x_1'}
\delta_{x_3,x_2''}
\delta_{x_3',x_0''}
\Big\{\begin{array}{ccc|cc}
\rho & \rho & \rho & x_0 & x_2\\
\rho & \rho & \rho & x_1 & x_3
\end{array}\Big\}_+
\Big\{\begin{array}{ccc|cc}
\rho & \rho & \rho & x_0' & x_2'\\
\rho & \rho & \rho & x_1' & x_3'
\end{array}\Big\}_+
\Big\{\begin{array}{ccc|cc}
\rho & \rho & \rho & x_0'' & x_2''\\
\rho & \rho & \rho & x_1'' & x_3''
\end{array}\Big\}_+\\
&\qquad+
\sum_{g\in G}
\Big\{\begin{array}{ccc|cc}
\rho & \rho & g & \bullet & x_2\\
\rho & \rho & \rho & x_1 & \bullet
\end{array}\Big\}
\Big\{\begin{array}{ccc|cc}
\rho & g & \rho & x_0' & x_2'\\
\rho & \rho & \rho & \bullet & \bullet
\end{array}\Big\}
\Big\{\begin{array}{ccc|cc}
\rho & \rho & \rho & \bullet & \bullet\\
\rho & g & \rho & x_1'' & x_3''
\end{array}\Big\}
\end{aligned}
\]
where
\[
a=x_0'\,,\quad
b=x_1\,,\quad
d=x_1''\,,\quad
a'=x_2'\,,\quad
c'=x_2\,,\quad
d'=x_3''\,.
\]
Using the values of the quantum $6j$-symbols from 
Proposition~\ref{prop.6j}, one confirms that the above identity 
is equivalent to the pentagon identity for the finite 
quantum dilogarithm.
\end{proof}

\begin{example}
\label{ex.S3}
As a consistency check, we use the $6j$ symbols given above to compute
the TQFT invariant of $S^3$. 
$S^3$ has an ordered triangulation with one positive tetrahedron and 
face pairings $0 \leftrightarrow 3$
and $1 \leftrightarrow 2$ and two distinct edges $03$ and all other 5 
identified (which form a trefoil):

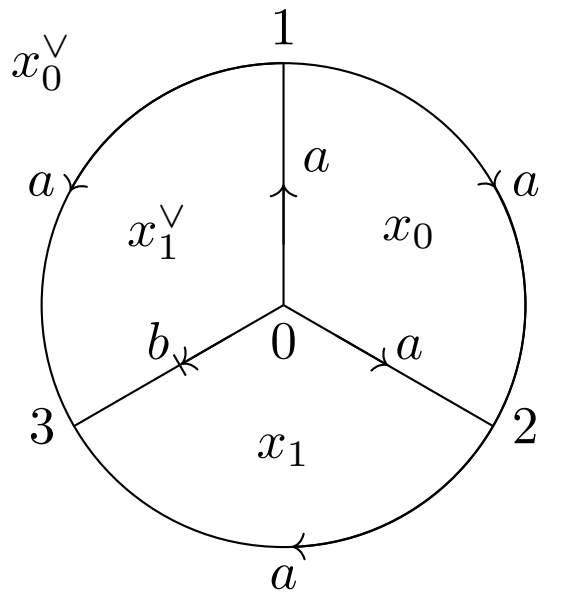
\begin{figure}[htbp!]
  \centering
  \begin{minipage}[c]{0.40\textwidth}
\centering
\begin{tikzpicture}[scale=0.8]
\draw(0,0) circle (2cm);
\draw[->](0,0)--(0,1);
\draw(0,0.5)--(0,2);
\draw[->](0,0)--(-30:1);
\draw(-30:0.5)--(-30:2);
\draw[->|](0,0)--(-150:1);
\draw(-150:0.5)--(-150:2);
\draw[-<](2,0) arc (0:2)--(30:2);
\draw[->](0,2) arc (90:2)--(150:2);
\draw[->](2,0) arc (0:2)--(-90:2);
\draw(0,0) node[below] {$0$};
\draw(0,2) node[above] {$1$};
\draw(-30:2) node[right] {$2$};
\draw(-150:2) node[left] {$3$};
\draw(0,1.2) node[right] {$a$};
\draw(-30:1.2) node[above] {$a$};
\draw(30:2) node[right] {$a$};
\draw(0,-2) node[below] {$a$};
\draw(150:2) node[left] {$a$};
\draw(-150:1.2) node[above] {$b$};
\draw(-2,2) node {$x_0^\vee$};
\draw(-90:1.2) node {$x_1$};
\draw(150:1.2) node {$x_1^\vee$};
\draw(30:1.2) node {$x_0$};
\end{tikzpicture}
  \end{minipage}
  \hfill 
  \begin{minipage}[c]{0.55\textwidth}
\begin{equation*}
\begin{array}{c|cccccc|}
\hbox{\diagbox[width=2.5em]{\tiny{tet}}{\tiny{edge}}} &  01 & 02 & 03 & 12 & 13 & 23 \\
\hline
    0  &            0 &  0 &  1 &  0 &  0 &  0 \\
\end{array}
\end{equation*}

\vspace{1em} 

\begin{equation*} 
\begin{array}{c|cccc|} 
\hbox{\diagbox[width=2.5em]{\tiny{tet}}{\tiny{face}}} & 012 & 013 & 023 & 123 \\
\hline    
    0  &             0 &  1 &  1 &  0  \\
\end{array}
\end{equation*}
  \end{minipage}
\caption{A triangulation of $S^3$.}
\label{fig:S3}
\end{figure}

Labelling edges $0$ and $1$ by $a$ and $b$, the invariant of
$S^3$ is given by

\begin{equation}
\label{ZS3}
\ZTVred_\gamma(S^3)=
\sum_{\substack{a,b \\ x_0, x_1}} 
\dimq(a) \dimq(b)
  \Big\{\begin{array}{ccc|cc}
    a & a & a & x_0 & x_1 \\
    a & a & b & x_1 & x_0
\end{array}\Big\}_+  \,.
\end{equation}
The bulk of the sum comes from $a=b=\rho$. Then the sum involves a 
delta function imposing $x_0=0$ resulting to
\begin{equation}
d^2 \sum_{\substack{a,b \\ x_0, x_1}} 
  \Big\{\begin{array}{ccc|cc}
    \rho & \rho & \rho & x_0 & x_1 \\
    \rho & \rho & \rho & x_1 & x_0
\end{array}\Big\}_+  = 
\frac{d}{\sqrt{N}}\sum_{x_1 \in G} \langle x_1 \rangle E(-x_1)
= d \varepsilon^{1/2} \mu_\gamma^{-1}
\end{equation}
The remaining sum involves either $b \in G$ and all other five edges equal to $\rho$:
\begin{equation}
   d  \sum_{\substack{b, x_0}} 
  \Big\{\begin{array}{ccc|cc}
    \rho & \rho & \rho & x_0 & \bullet \\
    \rho & \rho & b & \bullet & x_0
\end{array}\Big\}_+ = \frac{\mu_\gamma^{-1}}{\sqrt{N}} \sum_{x_0,b \in G}
\langle x_0 \rangle \langle b;x_0 \rangle =  N^{1/2} \mu_\gamma^{-1}
\end{equation}
or $a,b \in G$ giving:
\begin{equation}
\sum_{\substack{a,b \in G}} 
\Big\{\begin{array}{ccc|cc}
a & a & a & \bullet & \bullet \\
a & a & b & \bullet & \bullet
\end{array}\Big\}_+  = 1 \,.
\end{equation}
So, we obtain 
\begin{equation}
\ZTVred(S^3) =  
d \, \varepsilon^{1/2} \mu_\gamma^{-1} + N^{1/2} \mu_\gamma^{-1} + 1 =1 
\end{equation}
since $d \, \varepsilon^{1/2} + N^{1/2}=0$.
\end{example}

\begin{example}
\label{ex.41}
A second and substantially more interesting example is 
the finite TQFT invariant of the complement of the $4_1$ knot, 
which we compute and match it with the trace
of the quantum modular cocycle from~\cite{GZ23},
thus confirming Conjecture~\ref{conj.TQFT} for this example and 
closing a full circle of ideas.

Let $\calT$ denote the ideal triangulation of the complement 
of the figure-eight knot in $S^3$ with one positive tetrahedron 
$T_0$ and one negative tetrahedron $T_1$ with the edge and face 
identifications given in Figure~\ref{fig:4_1}.

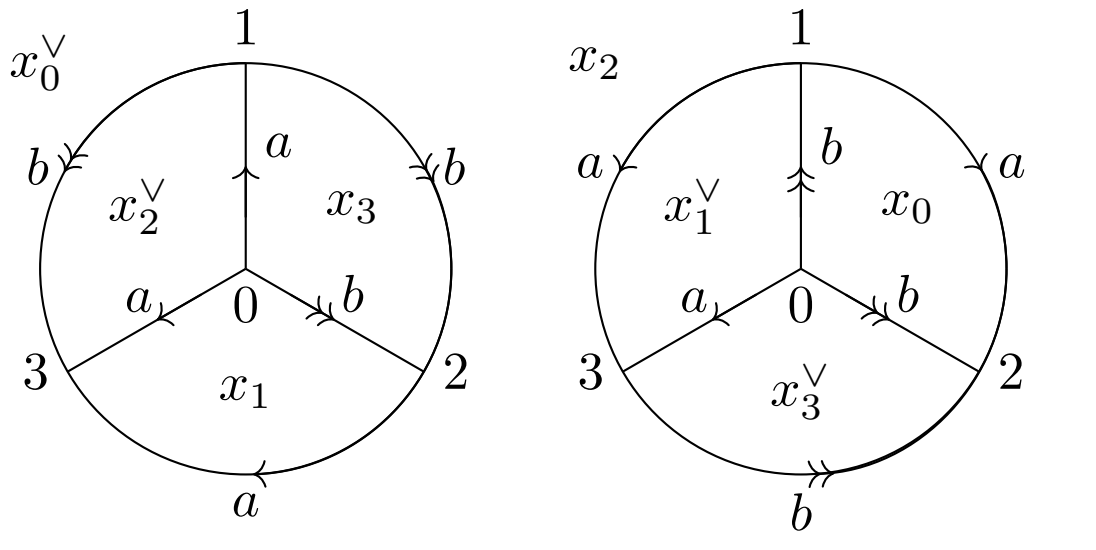
\begin{figure}[htbp!]
  \centering
  \begin{minipage}[c]{0.50\textwidth}
\centering
\begin{tikzpicture}[scale=0.68]
\draw(0,0) circle (2cm);
\draw[->](0,0)--(0,1);
\draw(0,0.5)--(0,2);
\draw[->>](0,0)--(-30:1);
\draw(-30:0.5)--(-30:2);
\draw[->](0,0)--(-150:1);
\draw(-150:0.5)--(-150:2);
\draw[-<<](2,0) arc (0:2)--(30:2);
\draw[->>](0,2) arc (90:2)--(150:2);
\draw[->](2,0) arc (0:2)--(-90:2);
\draw(0,0) node[below] {$0$};
\draw(0,2) node[above] {$1$};
\draw(-30:2) node[right] {$2$};
\draw(-150:2) node[left] {$3$};
\draw(0,1.2) node[right] {$a$};
\draw(-30:1.2) node[above] {$b$};
\draw(30:2) node[right] {$b$};
\draw(0,-2) node[below] {$a$};
\draw(150:2) node[left] {$b$};
\draw(-150:1.2) node[above] {$a$};
\draw(-2,2) node {$x_0^\vee$};
\draw(-90:1.2) node {$x_1$};
\draw(150:1.2) node {$x_2^\vee$};
\draw(30:1.2) node {$x_3$};

\draw[xshift=5.4cm](0,0) circle (2cm);
\draw[xshift=5.4cm,->>](0,0)--(0,1);
\draw[xshift=5.4cm](0,0.5)--(0,2);
\draw[xshift=5.4cm,->>](0,0)--(-30:1);
\draw[xshift=5.4cm](-30:0.5)--(-30:2);
\draw[xshift=5.4cm,->](0,0)--(-150:1);
\draw[xshift=5.4cm](-150:0.5)--(-150:2);
\draw[xshift=5.4cm,-<](2,0) arc (0:2)--(30:2);
\draw[xshift=5.4cm,->](0,2) arc (90:2)--(150:2);
\draw[xshift=5.4cm,-<<](2,0) arc (0:2)--(-90:2);
\draw[xshift=5.4cm](0,0) node[below] {$0$};
\draw[xshift=5.4cm](0,2) node[above] {$1$};
\draw[xshift=5.4cm](-30:2) node[right] {$2$};
\draw[xshift=5.4cm](-150:2) node[left] {$3$};
\draw[xshift=5.4cm](0,1.2) node[right] {$b$};
\draw[xshift=5.4cm](-30:1.2) node[above] {$b$};
\draw[xshift=5.4cm](30:2) node[right] {$a$};
\draw[xshift=5.4cm](0,-2) node[below] {$b$};
\draw[xshift=5.4cm](150:2) node[left] {$a$};
\draw[xshift=5.4cm](-150:1.2) node[above] {$a$};
\draw[xshift=5.4cm](-2,2) node {$x_2$};
\draw[xshift=5.4cm](-90:1.2) node {$x_3^\vee$};
\draw[xshift=5.4cm](30:1.2) node {$x_0$};
\draw[xshift=5.4cm](150:1.2) node {$x_1^\vee$};
\end{tikzpicture}
  \end{minipage}
  \hfill 
  \begin{minipage}[c]{0.47\textwidth}
\small 
\begin{equation*}
\begin{array}{c|cccccc|}
\hbox{\diagbox[width=2.5em]{\tiny{tet}}{\tiny{edge}}} &  01 & 02 & 03 & 12 & 13 & 23 \\
\hline
    0  &            0 &  1 &  0 &  1 &  1 &  0 \\
    1  &            1 &  0 &  1 &  0 &  0 &  1 \\
\end{array}
\end{equation*}

\vspace{0.4em} 

\begin{equation*} 
\begin{array}{c|cccc|} 
\hbox{\diagbox[width=2.5em]{\tiny{tet}}{\tiny{face}}} & 012 & 013 & 023 & 123 \\
\hline
    0  &             0 &  1 &  2 &  3  \\
    1  &             2 &  3 &  0 &  1  \\
\end{array}
\end{equation*}
  \end{minipage}

  \vspace{1em}
\caption{A triangulation of the complement of the $4_1$ knot.}
\label{fig:4_1}
\end{figure}

Labelling edges $0$ and $1$ by $a$ and $b$, we obtain that the invariant of
the complement of the $4_1$ knot is given by

\begin{equation}
\label{41}
\ZTVred_\gamma(S^3\sm 4_1)=
\sum_{\substack{a,b \\ x_0, x_1, x_2, x_3}} 
\dimq(a) \dimq(b)
  \Big\{\begin{array}{ccc|cc}
    a & b & b & x_0 & x_2 \\
    a & b & a & x_1 & x_3
\end{array}\Big\}_+  
   \Big\{\begin{array}{ccc|cc}
    b & a & a & x_2 & x_1 \\
    b & a & b & x_3 & x_0
\end{array}\Big\}_-   \,.
\end{equation}
The bulk of the sum comes from $a=b=\rho$ and 
$x_0=g$, $x_1=h$, $x_2=k$, $x_3=j$ with $g,h,k,j \in G$.
Then, the sum involves two delta functions imposing 
$g-h+k=0$ and $g-j+k=0$, with solutions $j=h$, $k=h-g$, 
resulting to

\begin{equation}
\frac{1}{N} \sum_{g,h} 
\frac{\langle h-g \rangle \langle g\rangle^{-1}}{
\langle h;h\rangle \langle h;h \rangle^{-1}}
E(g) \ang{h-g}E(g-h) \,.
\end{equation}
Since $|\langle h;h\rangle|=1$, after shifting
shifting $h \mapsto h+g$, 
the above integral decouples and equals to:
\begin{equation}
\frac{1}{N}\Big(\sum_g \langle g\rangle^{-1} 
E(g)\Big)\Big(\sum_g \langle g\rangle^2 E(g)\Big) 
= \tr(\Omega_{-1,1,\gamma}(\tau^*))
\tr(\Omega_{2,1,\gamma}(\tau^*))
\end{equation}
where the last equality follows from Theorem~\ref{thm.nahm}.
Compare this computation with~\cite[p.920]{AK14a}. 

The remaining part of the sum comes from tetrahedra where 
$a,b \in G$ and $x_0=x_1=x_2=x_3=\bullet$, and the 
delta functions equations at each of the triangle faces imply that 
$a=b=0$
and the contribution is $1$. 
\end{example}

\subsection{\texorpdfstring{$R$-matrices and representations of braid groups}{R-matrices and representations of braid groups}}
\label{sub.Rmat}

We now describe the $R$-matrices and the corresponding representations of the 
braid groups whose strands are colored by simple objects of $\calZ(\calC_\gamma)$.
Recall that 
$\calZ(\calC_\gamma)$ has four types of simple objects
$\mathfrak{a}_g,\mathfrak{b}_g, \mathfrak{c}_{g,h},\mathfrak{d}_u$
indexed by $g \in G_\gamma$ for the first two types, unordered pairs $(g,h)$
of different elements of $G_\gamma$ for the third type and the 
$J_\gamma=(G_{\gamma^2}\sm G_\gamma)/\gamma$, a set of $N(N+3)/2$ 
elements for the fourth type. The images of the simple objects under 
the forgetful functor $U: \calZ(\calC_\ga) \to \calC_\ga$ mapping 
$(X,e_X)$ to $X$
are given by 
$
U(\mathfrak{a}_g) = [g], 
U(\mathfrak{b}_g) = [g] + [\rho], 
U(\mathfrak{c}_{g,h}) = [g] + [h] + [\rho],
U(\mathfrak{d}_j) = [\rho] \,.
$
Recall that 
\begin{equation}
\Hom_{\calZ(\calC_\gamma)}((X,e_X), (Y,e_Y)) \subset 
\Hom_{\calC_\gamma}(X,Y)
\subset L_{2N}
\end{equation}

Let $u \in J_\gamma$ and color the strings of the braid group in
two strands by $\mathfrak{d}_u$. Using~\cite[Thm.6.6]{Izu01} we obtain that
the representation is given by
$B_2 \to \Aut(\mathfrak{d}_u^{\otimes 2})$ where $\sigma \mapsto e_{\mathfrak{d}_u}$
with
\begin{equation}
\label{sB1}
   e_{\mathfrak{d}_u} = \omega_u^{-1} \ang{\tau_u}_\gamma^{-1}\mu_\gamma^{-3}
   \Big(\sum_g S_g S'_g + E_2(u+g)^{-1} \mu_\gamma^2 T_g T'_g \Big)
   \frac{1}{\sqrt{N}} \Big( \sum_g 
   \langle \tau_u ; g \rangle \langle g \rangle^{-1}
   U_g \Big) \,.
\end{equation}
Likewise, recall $B_3=\langle \sigma_1, \sigma_2 \, | \sigma_1 \sigma_2 \sigma_1 = \sigma_1 \sigma_2 \sigma_1 \rangle$. 
The representation of the braid group $B_3$ with two strands is given by
\begin{equation}
\label{B3}
B_3 \to \Aut(\mathfrak{d}_u^{\otimes 3}), \qquad \sigma_1 \mapsto 
   s_1:= e_{\mathfrak{d}_u} \otimes 1 = e_{\mathfrak{d}_u}, \qquad 
\sigma_2 \mapsto 
  s_2= 1 \otimes  e_{\mathfrak{d}_u} =\rho(e_{\mathfrak{d}_u}) \,.
\end{equation}

Using the relations of $E_{\gamma^2}$ and $E_\gamma$, one can check the
Yang--Baxter equation:

\begin{equation}
\label{YB}
s_1 s_2 s_1 = s_2 s_1 s_2 \in L_{2N} \,. 
\end{equation}

If, on the other hand, we color the strings of $B_r$ by $\mathfrak{a}_g$,
the representation $B_r \to \Aut(\mathfrak{a}_g^{\otimes r})$ is given by
$\beta \mapsto \langle g \rangle^{2 \text{writhe}(\beta)}$ where writhe is
the signed sum of the crossings of $\beta$.

\begin{proof}[Proof of Proposition \ref{prop.Sigmag}]
The TQFT invariants for $S^1 \times \Sigma_g$, where $\Sigma_g$
is a closed surface of genus $g$, are given by the dimensions of
the vector spaces associated to $\Sigma_g$ by $\calZ(\calC_\gamma)$, and those are computed by the Verlinde
formula, which in our case is given in~\cite[Cor.8.9]{Izu01}:
\begin{equation}
\begin{aligned}
\ZTV_\gamma(S^1 \times \Sigma_g) & = \frac{N^{2 g - 1} (N - 1)}{2} 
+ \frac{N^g (N + 3) (N + 4)^{g - 1}}{2} \\ & + 
 N^{2 g - 1} (N + 4)^{g - 1} 
 \Big( \Big( \frac{\sqrt{N+4} + \sqrt{N}}{2} \Big)^{2 g - 2} 
 + \Big( \frac{\sqrt{N+4} - \sqrt{N}}{2} \Big)^{2 g - 2} 
 \Big) 
\end{aligned}
\end{equation}
with the first few values for $g=1,2,3$ given by
\begin{equation}
3 N + N^2, \,\, 6 N^2 + 11 N^3 + 7 N^4 + N^5, \,\,
24 N^3 + 20 N^4 + 37 N^5 + 81 N^6 + 50 N^7 + 12 N^8 + N^9 \,.
\end{equation}
where $N=\tr(\gamma)-2=\varepsilon+\varepsilon^{-1}-2$. Expanding the
binomial, it follows that $Z_{\gamma,S^1 \times \Sigma_g} \in 
\BZ_{\geq 0}[N]$
is a polynomial of degree $4g-3$ and lowest degree $g$ for $g \geq 2$.
It follows that
the invariant for $\gamma^n$ is given by replacing $N$ by  
$\varepsilon^n+\varepsilon^{-n}-2$ in the above formula.
Since $Z_{\gamma,S^1 \times \Sigma_g} \in \BZ_{\geq 0}[N]$ 
vanishes when $N=0$, the above substitution concludes the proof of 
Proposition~\ref{prop.Sigmag}.
\end{proof}

\begin{proof}[Proof of Equation~\eqref{ZLens}]
From Izumi's $S$ matrix and Equation~\eqref{TQFTLens} we can 
compute the invariants associated to Lens spaces $L(p,1)$:

\[
\begin{aligned}
\lambda_\gamma \ZTVred_\gamma(L(p,1))
&=
\sum_{g\in G_\gamma}\ang{g;g}_\gamma^{-p}
+
(d_\gamma+1)^2
\sum_{g\in G_\gamma}\ang{g;g}_\gamma^{-p}
+
(d_\gamma+2)^2
\sum_{g<h\in G_\gamma}\ang{g;h}_\gamma^{-p}
+
d_\gamma^2
\sum_{j\in (G_{\gamma^2}\sm G_\gamma)/\gamma}
\frac{\ang{j}_{\gamma^2}^p}{\ang{\varpi j}_\gamma^p}\\
&=
\sum_{g\in G_\gamma}\ang{g;g}_\gamma^{-p}
+
(d_\gamma+1)^2
\sum_{g\in G_\gamma}\ang{g;g}_\gamma^{-p}
+
\frac{(d_\gamma+2)^2}{2}
\Big(\sum_{g,h\in G_\gamma}\ang{g;h}_\gamma^{-p}-\sum_{g\in 
G_\gamma}\ang{g;g}_\gamma^{-p}\Big)\\
&\qquad+
\frac{d_\gamma^2}{2}
\Big(\sum_{j\in G_{\gamma^2}}
\frac{\ang{j}_{\gamma^2}^p}{\ang{\varpi j}_\gamma^p}-\sum_{g\in 
G_{\gamma}}
\ang{g}_\gamma^{-2p}\Big)\\
&=
\frac{(d_\gamma+2)^2}{2}
N_\gamma\big|G_\gamma[p]\big|
+
\frac{d_\gamma^2}{2}
\sum_{j\in G_{\gamma^2}}
\frac{\ang{j}_{\gamma^2}^p}{\ang{\varpi j}_\gamma^p}\,,
\end{aligned}
\]
where we use the fact that $\ang{ g(\mathrm{I}+\gamma)}
_{\gamma^2}\ang{g(\mathrm{I}+\gamma)}_\gamma^{-1}
=\ang{g}_\gamma^{-2}$.
Using 
\begin{equation}
\frac{(d_\gamma+2)^2 \lambda_\gamma}{2} N_\gamma 
= \frac{1}{2}(1+\varepsilon^{-1}),
\qquad
\frac{d_\gamma^2}{2 \lambda_\gamma} \sqrt{N_\gamma(N_\gamma+4)} 
= \frac{1}{2}(1-\varepsilon^{-1})
\end{equation}
the result follows. 
\end{proof}

\bibliographystyle{plain}
\bibliography{biblio}

\end{document}